\documentclass[leqno,12pt]{article}
\usepackage{enumerate,amsmath,amsthm,amssymb,ifthen}
\swapnumbers
\usepackage[numeric]{amsrefs}

\newtheorem{thm}[equation]{Theorem}
\newtheorem{cor}[equation]{Corollary}
\newtheorem{lem}[equation]{Lemma}
\newtheorem{prop}[equation]{Proposition}
\newtheorem{remark}[equation]{Remark}
\newtheorem{definition}[equation]{Definition}

\makeatletter
\renewcommand{\@seccntformat}[1]{\bf \S{\csname the#1\endcsname\enspace}}
\makeatother

\newcommand{\thmref}[1]{theorem~\ref{#1}}
\newcommand{\lemref}[1]{lemma~\ref{#1}}
\newcommand{\propref}[1]{proposition~\ref{#1}}
\newcommand{\corref}[1]{corollary~\ref{#1}}
\newcommand{\defref}[1]{definition~\ref{#1}}

\DeclareMathOperator{\sg}{sgn}

\numberwithin{equation}{section}

\newcommand{\gobble}[1]{}
  \newcommand{\rangeref}[2]{%
    \ref{#1}--\afterassignment\gobble\fam 0\ref{#2}%
  }

\renewcommand\d{\delta}

\renewcommand\l{\lambda}

\newcommand\D{\Delta}
\newcommand\G{\Gamma}
\newcommand\f{\frac}
\newcommand{\N}{{\mathbb{N}}}
\newcommand{\Z}{{\mathbb{Z}}}
\newcommand{\R}{{\mathbb{R}}}
\newcommand{\RP}{{\mathbb{RP}}}
\newcommand{\C}{{\mathbb{C}}}
\renewcommand{\D}{{\mathcal{D}}}

\newcommand{\Sch}{{\mathcal{S}}}

\renewcommand\O{{\mathcal O}}

\newcommand\I{{\mathcal I}}

\newcommand{\ttwo}[4]{\left(\begin{array}{cc}
{#1} & {#2} \\ {#3} & {#4} \end{array} \right)}

\renewcommand\({\left(}         
\renewcommand\){\right)}

\begin{document}

\title{Distributions and Analytic Continuation of Dirichlet Series}
\author{Stephen D. Miller\thanks{Supported by NSF grant
DMS-0122799 and an NSF post-doctoral fellowship}~~and Wilfried
Schmid\thanks{Supported in part by NSF grant DMS-0070714}}
\date{June 8, 2003}

\maketitle

\section{Introduction}\label{introduction}

Dirichlet series and Fourier series can both be used to encode sequences of complex numbers $\,a_n\,$, $n\in\N$. Dirichlet series do so in a manner adapted to the multiplicative structure of $\,\N$, whereas Fourier series reflect the additive structure of $\,\N$. Formally at least, the Mellin transform relates these two ways of representing sequences. In this paper, we make sense of the Mellin transform of periodic distributions and other tempered distributions, as a tool for the analytic continuation of various $L$-functions and the derivation of functional equations.

To illustrate what we mean, we sketch a heuristic argument for the functional equation of the Riemann zeta function. We let $\,\delta_n(x)$ denote the Dirac delta function at the point $\,n\in\Z$. The sum $\,\sum_{n\in\Z} \delta_n(x)$ is a tempered distribution; as such, it has a Fourier transform:
\begin{equation}
\label{intro1}
\mathcal F\bigl( \,\, {\sum}_{n\in\Z}\, \delta_n(x)\,\bigl) \ = \ {\sum}_{n\in\Z}\ e(nx)\qquad \bigl(\, e(x)\ =_{\text{def}}\ e^{2\pi i x}\,\bigr)\,.
\end{equation}
Here we are using L. Schwartz' normalization of the Fourier transform,
\begin{equation}
\label{intro2}
\widehat f(y) \ = \ \int_\R f(x)\ e(-xy)\,dy\,,
\end{equation}
for $\,f\in\Sch(\R)\,$\,=\,\,Schwartz space of $\,\R$. The Poisson summation formula is equivalent to the identity $\,\mathcal F(\sum_{n\in\Z} \delta_n(x))=\sum_{n\in\Z} \delta_n(x)$, hence
\begin{equation}
\label{intro3}
{\sum}_{n\in\Z}\, \delta_n(x) \ = \ {\sum}_{n\in\Z}\ e(nx)\,,
\end{equation}
as an equality of tempered distributions. We now formally integrate both sides against the ``even Mellin kernel" $\,|x|^{s-1}\,dx$, without worrying about convergence. Both $\,\delta_0(x)$ and the constant function $1$ have Mellin transform zero, in a sense that can be made precise. Neglecting these two terms and computing formally, one finds
\begin{equation}
\label{intro4}
\begin{aligned}
2\,\zeta(1-s)\ &= \ {\sum}_{n\neq 0}\,\int_\R \delta_n(x)\,|x|^{s-1}\,dx
\\
&= \ {\sum}_{n\neq 0}\ \int_\R e(nx)\,|x|^{s-1}\,dx\ = \ G_0(s)\,\zeta(s)\,,
\end{aligned}
\end{equation}
with $\,G_0(s)=\int_\R e(x)|x|^{s-1}dx=2(2\pi)^{-s}\Gamma(s)\cos(\pi s/2)$. That is Riemann's functional equation.

Our results on the Mellin transform of tempered distributions make the preceding formal argument perfectly rigorous. The presence of the constant function $1$ on the right hand side of (\ref{intro3}) accounts for the simple pole of $\,\zeta(s)$ at $\,s=1$, even though this term does not affect the functional equation itself. The same reasoning also gives the analytic continuation and functional equation for Dirichlet $L$-functions; details can be found in section \ref{examplesec}.

The use of distributions in the proof of the functional equation of $\,\zeta(s)$ might seem like a mere curiosity. However, distributions come up naturally in the study of automorphic forms on reductive groups. Classical modular forms and Maass forms on $GL(2,\R)$ have distribution boundary values. These are so-called automorphic distributions, periodic distributions $\,\tau\in C^{-\infty}(\R)$ which satisfy an equation of the type
\begin{equation}
\label{intro5}
\tau(x) \ = \ (\sg x)^\delta\,|x|^{\mu-1}\,\tau(1/x)\,.
\end{equation}
The $L$-function of the modular form or Maass form $\,F$ is the Dirichlet series formed from the Fourier coefficients -- suitably renormalized -- of the automorphic distribution that corresponds to $\,F$. To prove the analytic continuation and functional equation of such an $L$-function, we apply the Mellin transform to both sides of the equation (\ref{intro5}). This, too, is worked out in section \ref{examplesec}.

These are all known results, of course, and new proofs would not justify the writing of this paper. Its real purpose is to provide the analytic tools for our program, begun in \cite{MS2}, to study automorphic representations of higher rank groups from the point of view of automorphic distributions.

The central idea is a notion of distributions vanishing to a certain order along a submanifold of the manifold on which the distributions are defined. We introduce this notion in the next section, where we also deduce the most immediate consequences. For most of the rest of the paper, the compacti\-fied real line will play the role of the ambient manifold. In section \ref{onevariablesec}, we establish a number of equivalent criteria for the vanishing of a distribution of one variable at a point. We also define the signed Mellin transform of tempered distributions which vanish to sufficiently high order at the origin and the point at infinity, and we show that the Fourier transform $\,\widehat\sigma$ of a tempered distribution $\,\sigma$ vanishes to order $\,k$ at infinity if $\,\sigma$ vanishes to order $\,k$ at the origin. The properties of the Mellin transform of tempered distributions, in particular the interaction between the Mellin transform and the Fourier transform, are the subject of section \ref{4nmelsec}. The examples we mentioned earlier -- the Riemann zeta function, Dirichlet $L$-functions, and $L$-functions of automorphic forms on $GL(2)$ -- are worked out in section \ref{examplesec}. In section \ref{talphasec} we prove certain auxiliary statements for the Voronoi summation formula for $GL(3)$, which is the main result of \cite{MS2}. We return to the multi-variable case in the final section, where we discuss the summation and integration of distributions which vanish along a submanifold; these results are needed for the proof of the (known) converse theorem for $GL(3)$ in \cite{MS2}.

\section{Definitions and basic properties}\label{definitions}

In the following, $M$ will denote a $C^\infty$ manifold and $S \subset M$ a locally closed submanifold. We follow the convention of defining the space of distributions $C^{-\infty}(M)$ as the dual of the space of compactly supported, smooth measures. Functions and distributions take values in $\C$ unless we say otherwise. By means of the integration pairing between functions and measures, every $C^\infty$ function, and more generally every locally integrable function, can be regarded as a distribution:
\begin{equation}
\label{distribution1}
C^{\infty}(M) \ \subset \ L^1_{\text{loc}}(M)\ \subset \ C^{-\infty}(M)\,.
\end{equation}
Informally speaking, ``distributions transform like functions". We write the tautological pairing between distributions and compactly supported smooth measures as integration, since it extends the integration pairing between functions and measures. More generally, a distribution $\,\sigma$ can be paired against a smooth measure $\,dm\,$ if the intersection of their supports is compact. These tautological pairings make sense globally on the manifold $M$, quite independently of a choice of specific coordinate neighborhoods, but can be reduced to the analogous pairing on Euclidean space, by means of a suitable partition of unity.

We use the notation $\mathcal C^\infty_M$ for the sheaf of $C^\infty$ functions, $\mathcal C^{-\infty}_M$ for the sheaf of distributions, and $\I_S \subset \mathcal C^\infty_M$ for the ideal sheaf of the submanifold $S\subset M$. The term ``differential operator" will serve as shorthand for ``linear differential operator with $C^\infty$ coefficients". The differential operators constitute a sheaf of algebras $\D_M$ over the sheaf of rings $\mathcal C^\infty_M$. One calls a differential operator $D$ {\em tangential to} $S$ if
\begin{equation}
\label{tangential1}
D\,\I_S^{\,\,k} \ \subset \ \I_S^{\,\,k}\ \ \text{for every $k\in \N$}\,.
\end{equation}
If $D$ happens to be a vector field, this notion agrees with the usual, geometric notion of tangentiality: a vector field is tangential to $S$ if its values at all the points of $S$ lie in the tangent bundle $TS$. The differential operators which are tangential to $S$ constitute a sheaf of subalgebras of $\D_M$, which is generated over $\mathcal C^\infty_M$ by the sheaf of vector fields which are tangential to $S$; one can verify this assertion by a computation using suitably chosen local coordinates. We observe:
\begin{equation}
\label{tangential2}
\begin{aligned}
&\text{if a differential operator $D$ is tangential to $S$,}
\\
&\qquad\ \ \ \text{then so is its formal adjoint $D^*$},
\end{aligned}
\end{equation}
when the formal adjoint is defined relative to any particular Riemannian metric on $M$. Because of what was just said, it suffices to establish (\ref{tangential2}) for vector fields, which is a simple matter.

\begin{definition}\label{inf_order} A distribution $\,\sigma\in C^{-\infty}(M)\,$ vanishes to order $k\geq 0$ along the submanifold $S$ if every point $p\in S$ has an open neighborhood $U_p$ in $M$ with the following property: there exist differential operators $D_j$ on $U_p$ which are tangential to $S \cap U_p$, measurable locally bounded functions $h_j\in L^\infty_{\text{\rm{loc}}}(U_p)$, and $C^\infty$ functions $f_j\in C^\infty(U_p)$ which vanish to order $k$ on $S$, all indexed by $1\leq j \leq N$, such that
\begin{equation*}
\sigma \ = \  {\sum}_{1 \leq j \leq N} \ f_j \,D_j\, h_j\,,
\end{equation*}
as an identity between distributions on $U_p$. The distribution $\sigma$ vanishes to infinite order along $S$ if it vanishes to order $k$ for every $k\geq 0$.
\end{definition}

Let us record some formal consequences. If $\,0\leq k_1\leq k$, vanishing to order $k$ implies vanishing to order $k_1$. Since the definition does not involve a choice of coordinates, the notion of vanishing to order $k\leq\infty$ along a submanifold is preserved by diffeomorphisms. Vanishing to order $k\leq\infty$ along $S$ is a local condition. Also, if $\,\sigma,\tau\in C^{-\infty}(M)\,$ have this property, then so do $\sigma + \tau$ and the product $f\sigma$ with any $f\in C^\infty(M)$. To put it more succinctly, the distributions which vanish to order $k\leq\infty$ along $S$ constitute a subsheaf of $\mathcal C^{-\infty}_M$, viewed as sheaf of $\mathcal C^{\infty}_M$-modules.

\begin{lem}\label{differentiation}
If $D$ is a differential operator of degree $d$, and if $\,\sigma\in C^{-\infty}(M)\,$ vanishes to order $k\geq d$ along the locally closed submanifold $S\subset M$, the distribution $D\sigma$ vanishes to order $k-d$ along $S$. The distributions which vanish to infinite order along $S$ constitute a subsheaf of $\mathcal C^{-\infty}_M$, as sheaf of modules over the sheaf of differential operators $\mathcal D_M$.
\end{lem}

\begin{proof}
We may as well suppose $d=1$, and that $D$ is a vector field. When we express $\sigma$ as in \defref{inf_order},
\begin{equation}
\label{inf_order0}
D\,\sigma \ = \  {\sum}_{1 \leq j \leq N} \ (Df_j) \,D_j\, h_j\ + \ {\sum}_{1 \leq j \leq N} \ f_j\, D \,D_j\, h_j\,.
\end{equation}
The functions $Df_j$ vanish to order $k-1$ along $S$. The differential operators $DD_j$ may not be tangential to $S$, but can be made tangential by multiplication with $C^\infty$ functions which vanish on $S$. We shrink the neighborhood $U_p$, if necessary, so that each $f_j$ can be factored as a product of $k$ functions which vanish to order one. We take one of these factors to make $DD_j$ tangential to $S$, thereby reducing the order of vanishing of $f_j$ by one. This establishes the first assertion. The second follows formally.
\end{proof}

For some purposes, it is preferable to have presentations of a distribution $\,\sigma$ as in definition \ref{inf_order}, but with continuous functions $h_j$. One can accomplish that, at the expense of reducing the integer $k$:

\begin{remark}\label{cont_inf_order}
If the distribution $\,\sigma\in C^{-\infty}(M)$ vanishes to order $k_1$ along $S$, and if $0\leq k \leq k_1 - 2[\operatorname{dim}M/4]-2$, there exist differential operators $D_j$ on $U_p$ tangential to $S \cap U_p$, continuous functions $h_j\in C(U_p)$, and $C^\infty$ functions $f_j\in C^\infty(U_p)$ which vanish to order $k$ on $S \cap U_p$, $1\leq j\leq N$, such that
\begin{equation*}
\sigma \ = \  {\sum}_{1 \leq j \leq N} \ f_j \,D_j\, h_j\,.
\end{equation*}
If $\,\sigma$ vanishes to infinite order along $S$, an expression of this type exists for every $k>0$.
\end{remark}

Indeed, any $L^2$ function $h$ on an $m$-dimensional Riemannian manifold can be expressed locally as $h=\Delta^r \tilde h$, in terms of a continuous function $\tilde h$ and the Laplace operator $\Delta$, raised to the power $r=[m/4]+1$. From here on one can argue as in the proof of \lemref{differentiation}, moving the differential operator $\Delta^r$ from right to left instead of left to right.

\begin{lem}\label{support}
A distribution $\,\sigma$ which vanishes to order $k\geq 0$ along a locally closed submanifold $S\subset M$ of codimension at least one cannot have its support contained in $S$ unless $\sigma = 0$.
\end{lem}

We should note that vanishing to order $k\geq 1$ along an open submanifold implies vanishing on the submanifold. In that case, the codimension one hypothesis can be dropped. A similar lemma plays a crucial role in the proof, by Atiyah and the second named author, of Harish-Chandra's regularity theorem for invariant eigendistributions \cite{Atiyah}

\begin{proof}
This is a local problem. We may as well replace $M$ by an open neighborhood $U_p$ of some $p\in S$ on which $\sigma$ can be expressed as in \defref{inf_order}. Shrinking $U_p$, if necessary, we may suppose that there exist local coordinates $\{x_1,\,\dots,\,x_s,\,y_1,\,\dots,\,y_{n-s}\}$ on $U_p$ adapted to $S$, in the sense that
\begin{equation}
\label{inf_order1}
U_p \cap S \ = \ \{\,q\in U_p \, \mid \, x_1(q)= \dots = x_s(q)=0\,\}\,,
\end{equation}
and that the fibres of the map $\,U_p \to U_p \cap S\,$, $\,(x_i,y_j)\mapsto (0,y_j)$, correspond to balls centered at the origin. Assuming that the support of $\sigma$ is contained in $U_p \cap S$, we must show that $\int_{U_p}\! \sigma\,\psi \,dx_1\dots \,dx_s\, dy_1\dots \,dy_{n-s} = 0$\,,\, for all test functions $\,\psi\! \in\! C^{\infty}_c(U_p)$. Since $\psi$ has compact support in  $U_p$, we can choose $\,\phi\in C^{\infty}_c(U_p)$ such that
\begin{equation}
\label{inf_order2}
\phi \ \equiv \ 1 \ \ \text{on some neighborhood of}\,\ S\cap \operatorname{supp}\psi\,.
\end{equation}
Our hypotheses on the coordinate system ensure that the family of functions
\begin{equation}
\label{inf_order3}
\phi_t \in C^{\infty}_c(U_p)\,,\qquad  \phi_t(x_i,\,y_j)\ = \ \phi(t^{-1}\, x_i,\,y_j) \qquad(\,0 < t \leq 1\,)
\end{equation}
is well defined. Each $\phi_t$ inherits the property (\ref{inf_order2}) from $\phi$, which implies
\begin{equation}
\label{inf_order4}
\begin{aligned}
&\int_{U_p}\! \sigma\,\psi\,dx_1\dots \,dx_s\, dy_1\dots \,dy_{n-s} = \int_{U_p}\! \sigma\,\phi_t\psi\,dx_1\dots \,dx_s\, dy_1\dots \,dy_{n-s}
\\
&\qquad\qquad   = \ {\sum}_{1\leq j\leq N}\, \int_{U_p}\! h_j\,D^*_j\left(f_j\,\phi_t\,\psi\right)dx_1\dots \,dx_s\, dy_1\dots \,dy_{n-s}\,;
\end{aligned}
\end{equation}
the second step uses the expression for $\sigma$ given in \defref{inf_order}, and $D^*_j$ denotes the formal adjoint of $D_j$ with respect to the Euclidean metric. At this point, it suffices to show that $\int_{U_p}\! h_j\,D_j^*(f_j\phi_t\psi)\,dx_1\dots \,dx_s\, dy_1\dots \,dy_{n-s}\to 0$ as $t\to 0$, or more specifically, that
\begin{equation}
\label{inf_order5}
\operatorname{vol}\bigl(\operatorname{supp}\phi_t\bigr)\,\operatorname{sup}\left|h_j D^*_j\left(f_j\,\phi_t\,\psi\right)\right|\ \rightarrow \ 0\ \ \ \text{as}\ \ t \to 0\,,
\end{equation}
for $1\leq j\leq N$.

According to (\ref{tangential2}), $D^*_j$ is tangential to $U_p\cap S$. Thus, when $D^*_j$ is expressed as a linear combination $\sum_{I,J} a_{I,J}\,\frac{\partial^{|I|}\ }{\partial x^I}\,\frac{\partial^{|J|}\ }{\partial y^J}$ of monomials in the $\frac{\partial\ }{\partial x_i}$ and $\frac{\partial\ }{\partial y_j}$ with $C^\infty$ coefficients, each $a_{I,J}$ must vanish along $S$ to order equal to the total degree $|I|$ of normal derivatives. When a normal derivative of order $\ell$ is applied to $\phi_t$, the result is $t^{-\ell}$ times a bounded function, but partial derivatives of $\,\phi_t$ in directions tangential to $S$ are bounded indepentdently of $t$. As $t\to 0^+$, the diameter of the support of $\phi_t$ in the fibres of $\,U_p \to U_p \cap S\,$ shrinks down to $0$, linearly in $t$, hence
\begin{equation}
\label{inf_order6}
\operatorname{sup}\{\,|a_{I,J}(q)|\, \mid \, q \in \operatorname{supp}\phi_t\} = O(t^{|I|})\,,
\end{equation}
whereas $\,\psi$, the $h_j$ and $f_j$, and all the partial derivatives of $\,\psi$ and of the $f_j$ are uniformly bounded on the support of $\,\phi_t$. This bounds $|h_j D^*_j\left(f_j\,\phi_t\,\psi\right)|$ independently of $t$. Since the volume of the support of $\phi_t$ tends to $0$ in proportion to $t^s$, $s=\operatorname{codim}S\geq 1$, the estimate (\ref{inf_order5}) follows.
\end{proof}

If the functions $\,h_j\,$ in (\rangeref{inf_order4}{inf_order5}) are only locally $\,L^1$, one can still bound $|D^*_j\left(f_j\,\phi_t\,\psi\right)|$ independently of $t$. The supports of the $\phi_t$ shrink down to $S\cap\operatorname{supp}\phi$, which has volume $0$, so the integrals of the $|h_j|$ over the support of $\,\phi_t$ tend to $0$. Thus, even when $\,h_j\in L^1_{\text{loc}}(U_p)$, $\,\sigma$ must still vanish on $\,U_p$. For future reference we record this slight improvement of the lemma:

\begin{remark}\label{improved_support}
Suppose $S\subset M$ is a closed submanifold of codimension at least one. If $\,\sigma\in C^{-\infty}(M)$ can be represented, locally near each $p\in S$, as
\[
\sigma \ = \ {\sum}_{1\leq j \leq N}\ f_j\,D_j\,h_j\,,
\]
in terms of locally $L^1$ functions $h_j$, $C^\infty$ functions $f_j$, and differential operators $D_j$ which are tangential to $S$, then $\,\sigma$ cannot have its support contained in $S$ unless $\,\sigma = 0$.
\end{remark}

When the submanifold $S\subset M$ is not only locally closed but closed, one can restrict distributions from $M$ to $M-S$. In that situation, according to the lemma, a distribution $\sigma \in C^{-\infty}(M)$ which vanishes to order $k\geq 0$ along $S$ is completely determined by its restriction to $M-S$. This observation motivates the following terminology:

\begin{definition}\label{canonical}
A distribution $\,\tau\,$ defined on the complement $M-S$ of a closed submanifold $S\subset M$ has a canonical extension across $S$ if there exists a -- necessarily unique -- distribution $\,\sigma \in C^{-\infty}(M)$ that vanishes to infinite order along $S$ and agrees with $\,\tau$ on $M-S$.
\end{definition}

It may seem strange that we require $\,\sigma\,$ to vanish to infinite order along $\,S\,$ since vanishing to order $k\geq 0$ already makes the extension unique. Our definition is motivated by the applications we have in mind, which involve distributions vanishing to infinite order along a submanifold. Saying ``$\,\tau$ has a canonical extension, and the extension vanishes to infinite order" would sound too awkward! We shall be careful to distinguish between distributions vanishing to infinite order along $S$, which requires $\tau$ to be defined on all of $M$, and possessing a canonical extension across $S$, which applies when $\tau$ is defined on the complement of $S$.

Simple but prototypical examples of canonical extensions of distributions arise as follows. Let $\tau_0$ be a distribution on the space $\,\R^{1+r}\,$, with coordinates $(x,y) = (x,y_1,\dots y_r)$. We suppose that $\tau_0$ is periodic of period $1$ in all the variables, and that its Fourier series involves no terms independent of the variable $x$:
\begin{equation}
\label{periodic_dist1}
\tau_0(x,y) \ = \  {\sum}_{m\neq 0}\ {\sum}_{n\in \Z^r}\ c_{m,n}\, e(mx + ny)\ \ \ \ \ \ (\,e(u)=_{\text{def}}e^{2\pi i u}\ )\,;
\end{equation}
here $ny$ is shorthand for $\,\sum_j n_j y_j\,$, of course. The distribution
\begin{equation}
\label{periodic_dist2}
\tau(\,x\,,\,y\,) \ \ = \ \ \tau_0(\,1/x\,,\,y\,)
\end{equation}
is well defined on the complement of the hypersurface $S = \{x=0\}\subset \R^{1+r}$.

\begin{prop}\label{oscillation}
The distribution
\[
\tau(\,x\,,\,y\,) \ \ = \ \ {\sum}_{m\neq 0}\ {\sum}_{n\in \Z^r}\ c_{m,n}\, e(m/x + ny)
\]
has a canonical extension across $S$. In particular, each of the summands $c_{m,n} e(m/x + ny)$ extends canonically across $S$. The sum of the canonical extensions of the summands converges in the strong distribution topology and agrees with the canonical extension of $\tau$.
\end{prop}

To put the proposition into perspective, we should remark that vanishing to order $k\leq \infty$ along a submanifold $S\subset M$ does not define a closed subspace of $C^{-\infty}(M)$ in the strong distribution topology, or even in the weak dual topology: the distribution $\,e(t/x)$ vanishes to infinite order at $x=0$ when $t\neq 0$, but converges to $\,1$ in the weak distribution topology, as $\,t\to 0$.

\begin{proof} The Fourier coefficients $c_{m,n}$ of the distribution $\tau_0$ grow at most polynomially with the indices. Thus, for $\,k\,$ sufficiently large,
\begin{equation}
\label{periodic_dist3}
F_k(x,y) \ = \ (2\pi i)^{-3k}\, {\sum}_{m\neq 0}\ {\sum}_{n\in \Z^r}\ \frac{c_{m,n}\, e(mx + ny)}{m^k(\|n\|^2 + 1)^k}
\end{equation}
is a continuous periodic function, and
\begin{equation}
\label{periodic_dist4}
\tau_0(x,y) \ \ = \ \ \textstyle\frac{\partial^k\ }{\partial x^k} \bigl(\textstyle\sum_j \frac{\partial^2\ }{\partial y_j^2} - 4\pi^2\bigr)^k F_k(x,y)\,.
\end{equation}
An application of the chain rule gives the equation
\begin{equation}
\label{periodic_dist5}
\begin{aligned}
\tau(x,y)\, = \, x^k\,D_k \, G_k(x,y)\,,\ \ \ \text{with}\ \ G_k(x,y) \ = \ F_k(1/x,y)
\\
\text{and}\ \ \ D_k  \ = \ x^{-k}(-x^2\textstyle\frac{\partial\ }{\partial x})^k \bigl(\textstyle\sum_j \frac{\partial^2\ }{\partial y_j^2} - 4\pi^2\bigr)^k.
\end{aligned}
\end{equation}
Since $\,F_k(x,y)$ is bounded and continuous, the function $G_k(x,y)$ is defined, continuous, and bounded on $\,\R^{1+r}\!-S$, hence globally defined as $L^\infty$ function on $\,\R^{1+r}$. Applying $D_k$ to various powers of $x$ one finds that this differential operator is smooth even along $S$ and in fact tangential to $S$. This shows: for $k$ sufficiently large, $\sigma = x^kD_k G_k(x,y)$ is an extension of $\tau$ which vanishes to order $k$ along $S$. Because of \lemref{support}, $\sigma$ does not depend on the choice of $k$, and therefore vanishes to infinite order along $S$.

If the integer $k$ in (\rangeref{periodic_dist2}{periodic_dist3}) is chosen large enough, the Fourier series for $F_k(x,y)$ converges uniformly. Consequently the series
\begin{equation}
\label{periodic_dist6}
G_k(x,y) \ = \ (2\pi i)^{-3k}\, {\sum}_{m\neq 0}\ {\sum}_{n\in \Z^r}\ \frac{c_{m,n}\, e(m/x + ny)}{m^k(\|n\|^2 + 1)^k}
\end{equation}
converges in $L^1_{\text{loc}}(\R^{1+r})$, and that in turn implies convergence of the series
\begin{equation}
\label{periodic_dist7}
\sigma(x,y) \ = \ {\sum}_{m\neq 0}\ {\sum}_{n\in\Z^r}\,\ c_{m,n}\biggl( x^k\ D_k\ \frac{e(m/x + ny)}{m^k(\|n\|^2 + 1)^k(2\pi i)^{3k}\, } \biggr)
\end{equation}
in the strong distribution topology. The expression in parentheses represents the canonical extension of $e(m/x+ny)$, so the final assertion of the proposition follows.
\end{proof}

For the remainder of this section, $f \in C^\infty(M)$ will denote a real-valued function which has no critical points on its zero set. Then
\begin{equation}
\label{multi_mellon1}
S \ \ = \ \ \{\, p \in M\, \mid \, f(p)\,=\,0\,\}
\end{equation}
is a closed submanifold, of codimension one, and $f$ vanishes on $S$ to exactly first order. For $\alpha, \beta \in \C$ and $\delta \in \Z/2\Z$, the function $(\sg f)^\delta |f|^\alpha (\log |f|)^\beta$ is smooth on the complement of $S$. Thus, for any $\sigma \in C^{-\infty}(M)$, we may regard $(\sg f)^\delta|f|^\alpha(\log |f|)^\beta\sigma$ as a well defined distribution on $M-S$.

\begin{prop}\label{xalpha} If $\sigma \in C^{-\infty}(M)$ vanishes to order $k$ along $S$, and if $\,\operatorname{Re}\alpha > -k-1$, $\,(\sg f)^\delta|f|^\alpha(\log |f|)^\beta\sigma\in C^{-\infty}(M-S)$ has an extension $\,\tau\in C^{-\infty}(M)$, such that:\newline
\noindent {\rm a)}\, If $\,0\leq \ell < \operatorname{Re}\alpha + k$ for some integer $\ell$, or if $\beta=0$ and $\,0\leq \ell \leq \operatorname{Re}\alpha + k$, $\,\tau$ vanishes to order $\ell$ along $S$.\newline
\noindent {\rm b)}\, The extension $\tau=\tau(\alpha,\beta)$ depends holomorphically on $\alpha$ and $\beta$,\newline
\noindent in the sense that the integral of $\,\tau(\alpha,\beta)$ against any compactly supported smooth measure is holomorphic in the region $\{(\alpha,\beta)\in \C^2 \mid \operatorname{Re}\alpha > -k-1\}$. The conditions {\rm \,a)\,} and {\rm \,b)\,} determine the extension uniquely. In particular, if $\,\sigma$ vanishes to infinite order along $\,S$, $\,(\sg f)^\delta|f|^\alpha(\log |f|)^\beta\sigma$ has a canonical extension, which depends holomorphically on $\,(\alpha,\beta)\in\C^2$.
\end{prop}

\begin{proof}
This is a local problem, which needs to be verified only near points of $S$. Recall that $\mathcal I_S$ denotes the ideal sheaf of $S$. Because of our hypotheses,
\begin{equation}
\label{multi_mellon2}
\text{for each}\ \ n\in\N\,,\ \ f^n\ \ \text{generates}\ \ \mathcal I_S^n\,.
\end{equation}
Given $p\in S$, we choose an open neighborhood $U_p$ of $p$ as in \defref{inf_order}. Since $p$ is not a critical point of $f$, we can shrink $U_p$, if necessary, and introduce local coordinates $(x_1,\dots\,,x_r)$ on $U_p$, such that $x_1=f$. According to (\ref{multi_mellon2}), the functions $f_j$ in the statement of \lemref{inf_order} are divisible by $f^k=\,x_1^k\,$,\, hence
\begin{equation}
\label{multi_mellon3}
\sigma \ = \ x_1^k\ {\sum}_{1 \leq j \leq N} \ g_j \,D_j\, h_j\,,\ \ \text{with}\ \ g_j\in C^\infty(U_p)\,,\ \ h_j\in L^\infty_{\text{loc}}(U_p)\,,
\end{equation}
and $D_j$ tangential to $S$. We temporarily relax the hypothesis on $\ell$, requiring only that
\begin{equation}
\label{temp_hypo}
0 \leq \ell < \operatorname{Re}\alpha + k + 1\,.
\end{equation}
At least one such integer $\ell$ exists since $\operatorname{Re}\alpha > - k - 1$. To simplify various formulas, we set $\,\tilde\alpha = \alpha + k - \ell$, $\,\tilde\delta = \delta + k - \ell$. The hypotheses of the proposition require $\operatorname{Re}\tilde\alpha \geq 0$ or $\operatorname{Re}\tilde\alpha > 0$ depending on whether $\beta=0$ or not, but (\ref{temp_hypo}) allows $\,\operatorname{Re}\tilde\alpha > -1$. In either case, on the complement of $S$,
\begin{equation}
\label{multi_mellon4}
\begin{aligned}
&(\sg x_1)^\delta \,|x_1|^\alpha \,(\log |x_1|)^\beta \,\sigma \ =
\\
&\qquad\qquad = \ x_1^\ell\ {\sum}_{1 \leq j \leq N} \  g_j \,D_j\bigl((\sg x_1)^{\tilde\delta}\, |x_1|^{\tilde\alpha}\, (\log |x_1|)^\beta\, h_j\bigr)
\\
&\qquad\qquad\qquad - \ x_1^\ell\ {\sum}_{1 \leq j \leq N} \  g_j \,[\,D_j\,,(\sg x_1)^{\tilde\delta}|x_1|^{\tilde\alpha} (\log |x_1|)^\beta\,]\, h_j\,;
\end{aligned}
\end{equation}
here $[\,D_j\,,(\sg x_1)^{\tilde\delta}|x_1|^{\tilde\alpha} (\log |x_1|)^\beta\,]$ denotes the commutator of the differential operator $D_j$ with $(\sg x_1)^{\tilde\delta}|x_1|^{\tilde\alpha} (\log |x_1|)^\beta$, viewed as $0$-th order operator. Since
\begin{equation}
\label{multi_mellon5}
(\sg x_1)^{\tilde\delta}\, |x_1|^{\tilde\alpha}\, (\log |x_1|)^\beta\, h_j\ \in \ \begin{cases}\, L^\infty_{\text{loc}}(U_p) &\ \text{if $\,\operatorname{Re}\tilde\alpha> 0$} \\ \, L^\infty_{\text{loc}}(U_p) &\ \text{if $\beta=\operatorname{Re}\tilde\alpha=0$}\\ \, L^1_{\text{loc}}(U_p) &\ \text{if $\,\operatorname{Re}\tilde\alpha > -1$}\,,\end{cases}
\end{equation}
the first term on the right of (\ref{multi_mellon4}) is a distribution on all of $U_p$, which vanishes to order $\ell$ along $S\cap U_p$ when $\,\operatorname{Re}\tilde\alpha> 0$, or when $\beta=0$ and $\,\operatorname{Re}\tilde\alpha\geq 0$.

We now examine the second term on the right of (\ref{multi_mellon4}). In terms of the coordinates, the tangentiality of $D_j$ to $S$ means
\begin{equation}
\label{multi_mellon6}
D_j \ =\ {\sum}_{0 \leq s_1,\dots,s_r \leq d} \,\ a_{j;\,s_1,\dots,s_r} \ x_1^{s_1}\,\textstyle\frac{\partial^{s_1}\ }{\partial x_1^{s_1}}\, \frac{\partial^{s_2}\ }{\partial x_2^{s_2}}\, \dots \, \frac{\partial^{s_r}\ }{\partial x_r^{s_r}}\,,
\end{equation}
with coefficients $a_{j;\,s_1,\dots,s_r}\in C^\infty(U_p)$.  For computing the commutator, we note that $(\sg x_1)^{\tilde\delta}\,|x_1|^{\tilde\alpha} \,(\log |x_1|)^\beta$ commutes with the coefficients $a_{j;\,s_1,\dots,s_r}$ and with the derivatives $\frac{\partial\ }{\partial x_i}$, $i>1$. On the other hand,
\begin{equation}
\label{multi_mellon7}
\begin{aligned}
&[ \, x_1^s\,\textstyle\frac{\partial^s\ }{\partial x_1^s} \,, \,(\sg x_1)^{\tilde\delta}\,|x_1|^{\tilde\alpha}\, (\log |x_1|)^\beta\,]\, h_j\ =
\\
&\qquad = \ {\sum}_{s'\!,\,\,\beta'\!,\,\,\eta}\, P_{s'\!,\,\beta'\!,\,\eta}(\tilde\alpha,\beta)\, x_1^{s'}\,\textstyle\frac{\partial^{s'}\ }{\partial x_1^{s'}} \,\bigl( (\sg x_1)^\eta \, |x_1|^{\tilde\alpha}\, (\log |x_1|)^{\beta'}\, h_j\bigr)\,,
\end{aligned}
\end{equation}
as can be checked by induction on $s$; here $s'$ runs from $0$ to $s-1$, $\beta'$ ranges over the set $\{\beta - j \mid 0\leq j \leq s\}$, $\eta$ ranges over $\Z/2\Z$, and $P_{s'\!,\,\beta'\!,\,\eta}(\tilde\alpha,\beta)$ denotes a polynomial function of $\tilde\alpha$ and $\beta$. Most crucially, both sides of the equation involve the same complex power $|x_1|^{\tilde\alpha}$. We conclude that the second term on the right of (\ref{multi_mellon4}) has the same appearance as the first. It, too, represents a distribution on $U_p$. The sum of the two terms defines an extension of $(\sg x_1)^\delta |x_1|^\alpha(\log |x_1|)^\beta\sigma$ across $S\cap U_p$; when $\,0 \leq \ell < \operatorname{Re}\alpha + k\,$ or when $\beta=0$ and $\,0 \leq \ell \leq \operatorname{Re}\alpha + k$, this extension vanishes to order $\ell\geq 0$ along $S\cap U_p$. In that case, \lemref{support} guarantees the uniqueness of the local extensions, which then define a global extension of $(\sg f)^\delta |f|^\alpha(\log |f|)^\beta \,\sigma$ across $S$. The holomorphic dependence of (\ref{multi_mellon4}) on $\alpha$, which we are about to establish, implies the uniqueness of the local extensions even without the relaxed hypothesis (\ref{temp_hypo}). Alternatively the uniqueness of the local extensions can be deduced from remark \ref{improved_support}. In any case, we have extended $(\sg f)^\delta|f|^\alpha (\log |f|)^\beta \sigma\,$ across $S$ with the required order of vanishing.

The holomorphic dependence of the extension is a again a local problem, which needs to be verified only near points $p\in S$. We choose a coordinate neighborhood $U_p$ of $p$ as in the preceding argument. We must show that the integral of the local extension against any smooth measure $\psi\, dx_1\dots\, dx_r$, with compact support in $U_p$, depends holomorphically on $\alpha$ and $\beta$, provided $\operatorname{Re}\alpha > -k-1$, of course. The formula (\ref{multi_mellon4}), with $\,\ell=0$, expresses the local extension as a sum of two terms which, as we have argued, are really of the same type. It therefore suffices to show that
\begin{equation}
\label{multi_mellon8}
\begin{aligned}
(\tilde\alpha,\beta)\ &\mapsto \ \int_{U_p}\!\! \psi \,  g_j \,D_j\bigl((\sg x_1)^{\tilde\delta}\, |x_1|^{\tilde\alpha}\, (\log |x_1|)^\beta\, h_j\bigr)\,dx_1\,dx_2 \dots\, dx_r
\\
&= \ \int_{U_p}\!(\sg x_1)^{\tilde\delta}\, |x_1|^{\tilde\alpha}\, (\log |x_1|)^\beta\, h_j \,D_j^*\bigl(\psi \, g_j\bigr)\,dx_1\,dx_2 \dots\, dx_r
\end{aligned}
\end{equation}
describes a holomorphic function on the region $\{(\tilde\alpha,\beta)\in \C^2 \mid \operatorname{Re}\tilde\alpha > -1 \}$. Holomorphic dependence is clear if one integrates only over $U_p \cap \{ |x_1|\geq \delta \}$, for any small $\delta>0$. As $\delta$ tends to zero, these truncated integrals converge to the complete integral, locally uniformly in $\tilde\alpha$ and $\beta$. This implies the holomorphic nature of the integral (\ref{multi_mellon8}).
\end{proof}

\begin{remark}\label{negative_order}
In the case of a closed submanifold $S\subset M$ of codimension one, one can formally define the vanishing of a distribution $\,\sigma$ along $S$ to a negative power{\,\rm :} we say that $\,\sigma$ vanishes along $S$ to order $-k$, $k>0$, if locally near any point $p\in S$, $f^k\sigma$ vanishes to order $\,0\,$ for some, or equivalently any $C^\infty$ function $f$ which vanishes on $S$ exactly to first order. With this more general definition, the statement and proof of proposition \ref{xalpha} remain valid.
\end{remark}

We continue with the hypotheses of \propref{xalpha}. Since $f$ vanishes exactly to first order on $S$, $M$ is the disjoint union of $S$ and the two open subsets $\{\,f>0\,\}$ and $\{\,f<0\,\}$, both of which must be non-empty unless $S=\emptyset$. Equivalently, $M$ is the union of the two closed subsets $\{\,f\geq 0\,\}$, $\{\,f\leq 0\,\}$, which intersect exactly in $S$.

\begin{lem}\label{truncation}
If $\,\sigma\in C^{-\infty}(M)$ vanishes to order $k\geq 0$ along $S$, there exist distributions $\,\sigma_{f\geq 0}$ and $\,\sigma_{f\leq 0}$, both also vanishing to order $k$ along $S$, such that
\[
\sigma\ =\ \sigma_{f\geq 0} + \sigma_{f\leq 0}\,,\, \ \ \operatorname{supp}(\sigma_{f\geq 0}) \subset \{\,f\geq 0\,\}\,,\, \ \ \operatorname{supp}(\sigma_{f\leq 0}) \subset \{\,f\leq 0\,\}\,.
\]
These conditions determine $\,\sigma_{f\geq 0}$ and $\,\sigma_{f\leq 0}$ uniquely.
\end{lem}

\begin{proof}
Lemma \ref{support} implies the uniqueness. We may therefore argue locally, on some open neighborhood $U_p$ of $p\in S$, as in \defref{inf_order}, on which we represent $\,\sigma$ as
\begin{equation}
\label{inf_order7}
\sigma \ = \  {\sum}_{1 \leq j \leq N} \ f_j \,D_j\, h_j\,,\ \ \ \text{with}\ \ h_j \in L^\infty_{\text{loc}}(U_p)\,;
\end{equation}
the $D_j$ are tangential to $S$ and the $f_j\in C^\infty(U_p)$ vanish on $S$ to order $k$. Let $\,\chi_{f> 0}$ denote the characteristic function of the set $\{\,f>0\,\}$. Since $\chi_{f> 0}\,h_j$ is locally bounded,
\begin{equation}
\label{truncation1}
\sigma_{f\geq 0}\ \ =_{\text{def}}\ \ {\sum}_{1 \leq j \leq N} \ f_j \,D_j\bigl( \chi_{f> 0}\,h_j\bigr)\ \in \ C^{-\infty}(U_p)
\end{equation}
vanishes to order $k$ along $S\cap U_p$, has support in $S\cap\{\,f\geq 0\,\}$, and agrees with $\chi_{f> 0}\,\sigma$ on $U_p - S$. We define $\,\sigma_{f\leq 0}$ analogously. Then $\,\sigma_{f\geq 0} + \sigma_{f\leq 0} = \sigma$ at least on the complement of $S$, hence on all of $U_p$.
\end{proof}

As one consequence of the lemma, a distribution $\,\sigma$ on $\,M-S\,$ has an extension $\tau\in C^\infty(M)$ which vanishes to order $0\leq k \leq \infty$ along $S$ when it can be extended in this way ``from both sides of $S$". More precisely:

\begin{cor}\label{both_sides}
A distribution $\,\sigma\in C^{-\infty}(M-S)$ can be extended to a distribution $\,\tau$ on $M$ which vanishes along $S$ to order $0\leq k \leq \infty$ if and only if each $p\in S$ has an open neighborhood $U_p$ on which there exist distributions $\,\tau_{p,+}$, $\,\tau_{p,-}$, both vanishing to order $k$ along $\,U_p \cap S$, such that $\,\sigma = \tau_{p,+}$ on $\,U_p \cap \{\,f>0\,\}\,$ and $\,\sigma = \tau_{p,-}$ on $\,U_p \cap\{\,f<0\,\}$.
\end{cor}

The necessity of this condition is clear. To see the sufficiency, we note that $\,\tau = (\tau_{p,+})_{f\geq 0} + (\tau_{p,+})_{f\leq 0}\in C^{-\infty}(U_p)$ agrees with $\,\sigma$ on $U_p- S$ and vanishes to order $k$ along $U_p\cap S$. Lemma \ref{inf_order} then implies $\,\tau=\sigma$ on all of $U_p$ and ensures that the local extensions define a global distribution $\,\tau$ with the required properties.

\section{The case of one variable}\label{onevariablesec}

Many of our applications involve distributions on the real line or the compactified real line $\,\R\mathbb P^1 = \R \cup \{\infty\}$. For such distributions we shall make the definition of vanishing to order $k\geq 0$ more concrete. To simplify the discussion, we consider vanishing at $0$ or, occasionally, at $\infty$. By translation, that covers other points, too. Recall that a distribution vanishes to infinite order at a point if it vanishes to order $k$, for every $k\geq 0$.

\begin{lem}\label{one_var_def} Let $I\subset \R$ be an open interval containing the origin, and $k$ a non-negative integer. The following conditions on a distribution $\sigma \in C^{-\infty}(I)$ are equivalent:\newline
\noindent{\rm a)}\,\ $\sigma$ vanishes to order $k$ at the origin.\newline
\noindent{\rm b)}\,\ There exists an open interval $J$, with $\,0\in J\subset I$, an integer $N\geq 0$, and functions $h_j \in L^\infty_{\text{\rm{loc}}}(J)$, $0\leq j\leq N$, such that on $J$, $\sigma = \sum_{j=0}^N\,x^{k+j}\frac{d^j\ }{d x^j} h_j$.\newline
\noindent{\rm c)}\,\ There exists an open interval $J$, with $\,0\in J\subset I$, an integer $N\geq 0$, and functions $h_j \in L^\infty_{\text{\rm{loc}}}(J)$, $0\leq j\leq N$, such that on $J$, $\sigma = \sum_{j=0}^N\,\frac{d^j\ }{d x^j}(x^{k+j}h_j)$.\newline
\noindent{\rm d)}\,\ There exists an open interval $J$, with $\,0\in J\subset I$, an integer $N\geq 0$, and a function $h \in L^\infty_{\text{\rm{loc}}}(J)$, such that on $J$, $\sigma = x^k(\frac{d\ }{d x}\circ x)^N h$\,.\newline
\noindent If $\,\sigma \in C^{-\infty}(I)$ vanishes to order $k+1$ at the origin, it satisfies the three conditions {\rm \,b) -- d)\,} even with $h_j \in C(J)$, respectively $h \in C(J)$. If $k\geq 1$, and if $\,\sigma$ satisfies any of the conditions {\rm \,b) -- d)\,}, but with $h_j \in L^1_{\text{\rm{loc}}}(J)$, respectively $\,h \in L^1_{\text{\rm{loc}}}(J)$, then $\,\sigma$ vanishes to order $k-1$ at the origin.
\end{lem}

\begin{proof} A differential operator $D$ on the interval $J$ is tangential to the co\-dimension one submanifold $\{0\}\subset \R$ if and only if it can be expressed as a sum $D=\sum_{j=0}^N \,g_j x^j\frac{d^j\ }{d x^j}$, with $C^\infty$ coefficients $g_j$. If $f_j\in C^\infty(J)$ vanishes to order $k$ at $0$, the quotient $x^{-k}f_j$ is smooth, so
\begin{equation}
\label{one_var_1}
{\sum}_{j=0}^N\, f_j\,g_j\,x^j\,\textstyle\frac{d^j\ }{d x^j}\, = \, x^k\,\displaystyle{\sum}_{j=0}^N \, \tilde g_j\,x^j\textstyle\frac{d^j\ }{d x^j}\,,\,\ \text{with}\, \ \tilde g_j = x^{-k} f_jg_j\in C^\infty(J)\,.
\end{equation}
The $\tilde g_j$ can be moved across the derivatives, introducing new terms of order less than $N$, but with one or more ``excessive" power of $x$. Those can be moved across the derivatives, too, until eventually one obtains an expression of the type
\begin{equation}
\label{one_var_2}
{\sum}_{j=0}^N\ f_j\,g_j\,x^j\,\textstyle\frac{d^j\ }{d x^j}\ =\  x^k\,\displaystyle{\sum}_{j=0}^N \,x^j\textstyle\frac{d^j\ }{d x^j}\circ  \breve g_j\,\,,\,\ \text{with}\, \ \breve g_j \in C^\infty(J)\,.
\end{equation}
Thus a) implies b). By induction on $N$, one can show that the linear span of the differential operators $x^{k+j}\textstyle\frac{d^j\ }{d x^j}$, for $0\leq j\leq N$, coincides with the linear span of $\textstyle\frac{d^j\ }{d x^j}\circ x^{k+j}$, $0\leq j\leq N$, and also with the linear span of $x^k(\textstyle\frac{d\ }{d x}\circ x)^j$, $0\leq j\leq N$. In particular, b) is equivalent to c), and c) implies the existence of an expression
\begin{equation}
\label{one_var_3}
\sigma \ \ = \ \ {\sum}_{j=0}^N\,\ x^k\,\textstyle(\frac{d\ }{d x}\circ x)^j\, h_j\,,\ \ \ \text{with}\ \ h_j \in L^\infty_{\text{\rm{loc}}}(J)\,.
\end{equation}
The antiderivative $H_j$ of a function $h_j\in L^\infty_{\text{\rm{loc}}}(J)$, normalized by the condition $H_j(0)=0$, is H\"older continuous of index $1$, which makes  $x^{-1}H_j(x)$ locally bounded on $J$ and continuous on $J-\{0\}$, hence locally $L^\infty$ on $J$. Repeating this process $N-j$ times, one can solve the equation $\,(\frac{d\ }{d x}\circ x)^{N-j}\tilde h_j = h_j\,$ for $\tilde h_j\in L^\infty_{\text{\rm{loc}}}(J)$. This turns (\ref{one_var_3}) into the equation asserted by d), with $h = \sum_j \tilde h_j$, so c) implies d). The differential operator $(\frac{d\ }{d x}\circ x)^N$ is tangential to $\{0\}\subset\R$, so d) certainly implies a). At this point, we have established the equivalence of a) -- d).

If $\,\sigma$ vanishes to order $k+1$ at $0$, the condition b)\, with $k+1$ in place of $k$ gives the expression $\,\sigma = \sum_j x^{k+1+j}\frac{d^j\ }{d x^j} h_j = \sum_j x^{k+j+1}\frac{d^{j+1}\ }{d x^{j+1}} H_j$, where $H_j\in C(J)$ again denotes an antiderivative of $h_j$. In other words, we can replace the $h_j$ in b)\, by continuous functions if we lower the integer $k$ by $1$; this improves on remark \ref{cont_inf_order} in the special case of a one dimensional manifold. The equivalence of b) -- d) in the setting of continuous functions $h_j$ and $h$ follows from the same arguments as in the $L^\infty_{\text{\rm{loc}}}$ setting.

Every locally $L^1$ function $h_j$ has a continuous antiderivative $H_j$. Thus, if $\,\sigma$ satisfies b) or c) with $h_j \in L^1_{\text{\rm{loc}}}(J)$, we can write $\,\sigma = x^{k-1}\sum_j D_j H_j$, with $D_j$ tangential to $\,\{0\,\}$ and $H_j$ continuous, hence locally $L^\infty$. Except for the notation, condition d) with $h \in L^1_{\text{\rm{loc}}}(J)$ can be treated the same way.
\end{proof}

If a distribution $\,\sigma$ vanishes to order $k\geq 0$ at $0$, the statement of \lemref{one_var_def} allows the open intervals $J$ to depend on the particular choice of $k$. It is not difficult to show that for $J$, one can take any bounded interval $J$ around the origin whose closure is contained in the domain of definition $I$. More importantly, if $\,\sigma$ has compact support, one can express it in terms of compactly supported functions $h_j$\,:

\begin{lem}\label{compactsupport} If $\,\sigma\in C_c^{-\infty}(I)$ vanishes to order $k\geq 0$ at $0\in I$, there exist presentations of $\,\sigma$ as in the statements {\,\rm b)\,} and {\,\rm c)\,} in \lemref{one_var_def}, but with $J=I$ and functions $h_j \in L^\infty(I)$ which vanish outside some compact subinterval of $I$. If $\,\sigma\in C_c^{-\infty}(I)$ vanishes to order $k+1$ at $0$, there exist presentation as in {\,\rm b)\,} and {\,\rm c)\,}, with $J=I$ and $h_j \in C_c(I)$.
\end{lem}

\begin{proof}
We use the notation of the statement and proof of \lemref{one_var_def}. Let us observe first of all that the equivalence b) and c) works independently of the degree of regularity of the functions $h_j$ and does not affect the size of their supports. We may therefore concentrate on the condition b). If $\,\sigma$ vanishes to order $k$ at $0$, there exists an open subinterval $J\subset I$ containing $0$ on which $\,\sigma$ can be expressed as in b). We choose $\,\phi\in C_c^\infty(J)$ such that $\,\phi(x)\equiv 1$ near $x=0$. Then $\,\sigma = \phi\,\sigma + (1-\phi)\sigma$, both summands have compact support, and $(1-\phi)\sigma$ vanishes near the origin. Since $\,\phi$ has compact support in $J$,
\begin{equation}
\label{compactsupport1}
\begin{aligned}
\phi\,\sigma\ &= \ {\sum}_{j=0}^N\ \phi\,x^{k+j}\,\textstyle\frac{d^j\ }{dx^j}\,h_j\
\\
&= \ \displaystyle{\sum}_{j=0}^N\ x^{k+j}\,\textstyle\frac{d^j\ }{dx^j}\,(\phi\,h_j) \ - \ \displaystyle{\sum}_{j=0}^N\ x^{k+j}\,\textstyle[\,\frac{d^j\ }{dx^j},\,\phi\,]\,h_j\,,
\end{aligned}
\end{equation}
as an identity of compactly supported distributions on $I$, or even $\R$. The functions $\phi h_j$ lie in $L^\infty(I)$ and vanish outside the support of $\phi$. The commutator $[\,\frac{d^j\ }{dx^j}\,,\,\phi\,]$ can be expressed as a sum $\sum_{i=0}^{j-1}\frac{d^i\ }{dx^i}\circ \psi_{j,i}$ with coefficients $\,\psi_{j,i}$ which are linear combinations of derivatives of $\,\phi$. We can move the ``excessive" $(j-i)$-th power of $x$ to the right, as explained in the proof of \lemref{one_var_def}. The upshot is an expression for $\,\phi\,\sigma$ of the same type as b), in terms of $L^\infty$ functions $h_j$ whose support is contained in the support of $\phi$. If $\,\sigma$ vanishes to order $k+1$ at the origin, \lemref{one_var_def} allows us to require the functions $h_j$ in (\ref{compactsupport1}) to lie in $C(J)$. After the manipulation which was just described, the redefined functions $h_j$ are continuous, and their support still lies in the support of $\,\phi$. Thus, in either case, $\,\phi\sigma$ has been expressed in the form in which $\,\sigma$ needs to be expressed.

As a distribution of compact support, $(1-\phi)\sigma = F^{(N)}$ has a continuous $N$-th anti-derivative $F$, for every sufficiently large integer $N$; $F$ need not have compact support, of course. Since the origin does not lie in the support of $(1-\phi)\sigma$, which is compact and contained in $I$, there exists $\,\psi\in C_c^\infty(I)$ such that $\,\psi(x)\equiv 0\,$ near $\,x=0\,$ and $\,\psi(x)\equiv 1$ on an open neighborhood of the support of $(1-\phi)\sigma$. Hence
\begin{equation}
\label{compactsupport2}
(1-\phi)\sigma\, = \, \psi\,F^{(N)}\, = \, x^{k+N}\textstyle\frac{d^N\ }{dx^N} \bigl(x^{-k-N}\psi F\bigr)  -\,  x^{k+N}\, [\,\textstyle\frac{d^N\ }{dx^N},\,x^{-k-N}\psi\,]F\,.
\end{equation}
Since $\psi$ has compact support and vanishes near $x=0$, $\,x^{-k-N}\psi F$ lies in $C_c(I)$. For the same reason we can transform $\,x^{k+N}[\,\textstyle\frac{d^N\ }{dx^N},\,x^{-k-N}\psi\,]F$ into a sum $\,\sum_{j=0}^{N-1}x^{k+j}\textstyle\frac{d^j\ }{dx^j}F_j\,$ with $\,F_j\in C_c(I)$. Hence also $(1-\phi)\sigma$ has an expression of the required type.
\end{proof}

Because of \lemref{support}, a distribution $\,\sigma\in C^{-\infty}(I)$, defined on an open interval $I$, and vanishing to order $k\geq 0$ at $0\in I$, is completely determined by its restriction to $I-\{0\}$. To make this precise, we choose a cutoff function $\,\phi\in C^\infty_c(I)$ such that $\phi \equiv 1$ near the origin, and we define $\,\phi_t(x)=\phi(x/t)$, for $0< t \leq 1$. Then $\,\phi_t$ also has compact support in $I$, and $(1-\phi_t)\sigma$ vanishes near the origin -- in particular, $(1-\phi_t)\sigma$ depends only on the restriction of $\sigma$ to $I-\{0\}$.

\begin{lem}\label{approximation} If $\,\sigma\in C^{-\infty}(I)$ vanishes to order $k\geq 0$ at $0\in I$, $(1-\phi_t)\sigma \to \sigma$ as $t\to 0$, in the strong distribution topology.
\end{lem}

Convergence in the strong distribution topology implies convergence in the weak dual topology. Hence, for any $\psi\in C^{\infty}(I)$,
\begin{equation}
\label{one_var_4}
\int_I \sigma(x)\, \psi(x)\, dx \ \ = \ \ {\lim}_{t\to 0} \int_I \sigma(x)\, (1-\phi_t)(x)\,\psi(x)\, dx\,,
\end{equation}
provided $\,\sigma\in C^{-\infty}(I)$ vanishes to order $k\geq 0$ at $0\in I$.

\begin{proof}
We apply the criterion b) in \lemref{one_var_def}, with $k=0$. For $t$ small enough, the support of $\,\sigma - (1-\phi_t)\sigma = \phi_t \sigma\,$ is contained in $J$, so
\begin{equation}
\label{one_var_5}
\begin{aligned}
\phi_t \, \sigma \ \ &= \ \ {\sum}_{j=0}^N\ x^{j}\,\phi_t\,\textstyle\frac{d^j\ }{d x^j} h_j
\\
&=\ \ {\sum}_{j=0}^N\ \textstyle\frac{d^j\ }{d x^j}(\,x^{j}\,\phi_t h_j)\ \ - \ \ \displaystyle{\sum}_{j=0}^N\ [\,\textstyle\frac{d^j\ }{d x^j}\,,\,x^{j}\,\phi_t\,]\, h_j\,.
\end{aligned}
\end{equation}
As $t$ tends to $0$, the support of $\phi_t$ shrinks down to zero linearly in $t$. Thus $x^{j}\phi_t h_j \to 0$ in $L^1$ norm, hence in the strong distribution topology. Differentiation is continuous with respect to the strong dual topo\-logy. This allows us to conclude that the first term on the right of (\ref{one_var_5}) tends to $0$ as $t\to 0$. The commutator $[\,\textstyle\frac{d^j\ }{d x^j}\,,\,x^{j}\phi_t\,]$ can be expressed as a linear combination
\begin{equation}
\label{one_var_6}
[\,\textstyle\frac{d^j\ }{d x^j}\,,\,x^{j}\,\phi_t\,] \ \ = \ \ \displaystyle{\sum}_{\stackrel{\scriptstyle{i\geq 0,\,\ell\geq 0}}{1\leq i+\ell\leq j}} \,\ c_{i,\ell}\,t^{-\ell}\textstyle\frac{d^{j-i-\ell}\ }{d x^{j-i-\ell}}\circ \bigl(x^{j-i}\,\phi^{(\ell)}(x/t)\bigr) \,.
\end{equation}
In the summation $j-i\geq \ell$, and the diameter of the support of $\phi^{(\ell)}(x/t)$ is $O(t)$, so $t^{-\ell}\,x^{j-i}\,\phi^{(\ell)}(x/t)\to 0$ in $L^1$ norm. Arguing as in the first case, we see that also the second term on the right of (\ref{one_var_5}) tends to $0$ as $t\to 0$.
\end{proof}

Lemma \ref{approximation} implies an analogous statement about distributions defined near $\infty$ in $\,\R\mathbb P^1=\R\cup\{\infty\}$. Let $I\subset\R\mathbb P^1$ be a connected open neighborhood of $\infty$. We choose a function $\phi\in C^\infty_c(\R)$ as before, i.e., with $\phi\equiv 1$ near $0$, and we again define $\phi_t(x) = \phi(x/t)$, but this time for $t\geq 1$. Then $\phi_t$ vanishes near $\infty$, but $\phi_t(x) \to 1$ as $t\to \infty$, for any $x\in\R$.

\begin{cor}\label{approximationcor} If $\,\sigma\in C^{-\infty}(I)$ vanishes to order $k\geq 0$ at $\infty$, $\phi_t\sigma \to \sigma$ as $t\to \infty$, in the strong distribution topology, and hence also in the weak dual topology.
\end{cor}

This follows from \lemref{approximation} via the change of coordinates $x\rightsquigarrow 1/x$ and the substitution of $\phi(x)$ for $(1-\phi)(1/x)$.

Just as in the case of \lemref{support}, the proof of \lemref{approximation} establishes more than is claimed by its statement. If the functions $h_j$ in (\ref{one_var_4}) lie only in $L^1_{\text{loc}}(J)$, we can still conclude that both $x^{j}\phi_t h_j$ and $t^{-\ell}x^{j-i}\,\phi^{(\ell)}(x/t) h_j$ tend to $0$ in $L^1$ norm as $t\to 0$. The lemma therefore remains valid in this greater degree of generality:

\begin{remark}\label{approx_rem}
If $\,\sigma\in C^{-\infty}(I)$ can be expressed as $\,\sigma=\sum_{j=0}^N x^j\frac{d^j\ }{d x^j}h_j$, with $h_j\in L^1_{\text{\rm{loc}}}(J)$, $(1-\phi_t)\sigma$ converges to $\,\sigma$ as $t\to 0$, in the strong distribution topology.
\end{remark}

Recall the notion of a tempered distribution: a distribution $\,\sigma\in C^{-\infty}(\R)$ such that the integration pairing $C^\infty_c(\R)\ni\psi\mapsto \int_\R \sigma\,\psi\, dx$ extends continuously from $C^\infty_c(\R)$ to the space of Schwartz functions $\Sch(\R)$. Equivalently, tempered distributions on the real line can be characterized as those which arise as the $\ell$-th derivative $h^{(\ell)}(x)$, for some $\ell\in\N$, of a continuous function $h(x)$ growing at most polynomially as $|x|\to\infty$. Any distribution which can be extended from $\,\R$ to a distribution on $\,\R\mathbb P^1$ necessarily has this property.

The integration pairing exhibits the space of tempered distributions $\Sch'(\R)$ as the continuous dual of the Schwartz space $\Sch(\R)$. The Fourier transform
\begin{equation}
\label{fourier1}
\Sch(\R)\ \ni\ f(x)\ \mapsto\ \widehat f(y) \ \ = \ \ \int_\R f(x)\,e(-xy)\,dx\,,\ \ \ \bigl(\,e(u) = e^{2\pi i u}\,\bigr)
\end{equation}
sends Schwartz functions to Schwartz functions. Since
\begin{equation}
\label{fourier2}
f(x) \ \ = \widehat{\widehat f\,\,}\!(-x)\ \ = \ \ \int_\R \widehat f(y)\,e(xy)\,dy\ \ \ \ \ \ \ \bigl(\,f \in \Sch(\R)\,\bigr)\,,
\end{equation}
by Fourier inversion, the Fourier transform establishes an automorphism of $\Sch(\R)$. This makes it possible to define the Fourier transform $\,\widehat\sigma$ of a tempered distribution $\,\sigma$ by the formula
\begin{equation}
\label{fourier3}
\int_\R \widehat \sigma(x)\, f(x)\,dx \ \ = \ \ \int_\R \sigma(y)\,\widehat f(y)\,dy\ \ \ \ \ \ \ \bigl(\,f\in\Sch(\R)\,\bigr)\,,
\end{equation}
which reduces to Parseval's identity when the tempered distribution $\,\sigma$ happens to be a Schwartz function.

Our next result generalizes the one-variable version of \propref{oscillation}. We consider a periodic distribution without constant term,
\begin{equation}
\label{fourier4}
\tau(x) \ \ = \ \ {\sum}_{n\neq 0}\ a_n\, e(nx)\,.
\end{equation}
Like any periodic distribution, $\tau$ is tempered. Its inverse Fourier transform
\begin{equation}
\label{fourier5}
\widehat \tau(x) \ \ = \ \ {\sum}_{n\neq 0}\ a_n\,\delta_n(x)\ \ \ \ \ \ \ \bigl(\,\delta_n\,=\,\text{delta function at $n$}\,\bigr)\,,
\end{equation}
vanishes identically near $x=0$, hence vanishes to infinite order at $0$. According to \propref{oscillation}, $\tau$ extends canonically across $\infty$. This is one instance of a general phenomenon:

\begin{thm}\label{fourier}
If a tempered distribution $\,\sigma$ on the real line vanishes to order $k\geq 0$ at the point $\,0\,$, then its Fourier transform can be extended to a distribution on $\R\cup\{\infty\}$ which vanishes to order $k$ at $\,\infty\,$. In particular, if $\,\sigma$ vanishes to infinite order at $0$, the Fourier transform $\,\widehat\sigma$ has a canonical extension across $\infty$.
\end{thm}

The theorem does not have a converse: for example, if $\widehat\sigma\in\Sch(\R)$ has compact support, $\sigma$ is a smooth function whose Taylor series at the origin need not vanish to any order. Theorem \ref{4nmelthm} below identifies the obstruction to the converse of the theorem.

\begin{proof} We choose a cutoff function $\,\phi\in C^\infty_c(\R)$ such that $\phi \equiv 1$ near $0$. Then $\,\sigma=\phi\,\sigma + (1-\phi)\,\sigma\,$ is the sum of two tempered distributions, one of which has compact support and vanishes to infinite order at $0$, whereas the other vanishes identically on a neighborhood of zero. It therefore suffices to deal separately with these two cases.

Let us suppose first that the tempered distribution $\,\sigma\,$ has support away from the origin. Like any tempered distribution, we can express $\,\sigma\,$ as the $\ell$-th derivative of a continuous function $F$, for some $\ell\geq 0$, with $F$ growing at most polynomially. Since $\,\sigma = F^{(\ell)}$ vanishes on some open interval $J$ around $\,0\,$, the restriction of $F$ to $J$ must be a polynomial of degree at most $\ell-1$. We can subtract the polynomial from $F$, which allows us to assume that $F \equiv 0\,$ on $J$. Then, for $k$ sufficiently large, $F_k(x) =_{\text{def}} x^{-k}F(x)$ is not only continuous but decays like, say, $\,|x|^{-2}$ as $|x| \to \infty$. That makes the Fourier transform $\widehat F_k$ bounded and continuous. Since $\,\sigma(x) = \frac{d^\ell\ }{dx^\ell}(x^k F_k(x))$,
\begin{equation}
\label{fourier6}
\begin{aligned}
\widehat\sigma(1/x)\ &= \ (-1)^k\,(2\pi i)^{\ell-k}\,x^{-\ell}\,\left(\textstyle\frac{d^k\ }{dx^k}\widehat{F_k}\right)(1/x)
\\
&= \ (-1)^k\,(2\pi i)^{\ell-k}\,x^{-\ell}\,\left(-x^2\textstyle\frac{d\ }{dx}\right)^k\left(\widehat{F_k}(1/x)\right)\,,
\end{aligned}
\end{equation}
for all sufficiently large $k$ but with $\ell$ fixed. As a bounded function which is continuous away from the origin, $\widehat{F_k}(1/x)$ certainly lies in $L^\infty_{\text{loc}}(\R)$. Also, $(-x^2\textstyle\frac{d\ }{dx})^k$ can be expressed as $x^k D_k$, in terms of a differential ope\-rator $D_k$ which is tangential to $\{0\}\subset \R$. Thus (\ref{fourier6}) defines an extension of $\,\widehat\sigma$ across $\infty$ which vanishes there to order $k-\ell$ for all large $k$, hence to infinite order.

For the remaining case, we suppose that $\,\sigma$ has compact support and vanishes to order $k\geq 0$ at $0$. We use lemmas \ref{one_var_def} and \ref{compactsupport} to write
\begin{equation}
\label{fourier7}
\begin{aligned}
\sigma(x) \ = \ {\sum}_{0 \leq j \leq N} \, \textstyle\frac{d^{j}\ }{dx^{j}}\left( x^{k+j} h_j(x)\right)\,, \ \ \ &\text{with}\ \ h_j \in L^\infty(\R)
\\
&\text{and}\ \ h_j(x)\equiv 0 \ \ \text{for} \ \ |x|\gg 1\,.
\end{aligned}
\end{equation}
Then $\widehat h_j \in C^\infty(\R)$ and $\,\widehat\sigma(x) =  (-1)^{k+j}(2\pi i)^{-k}\,{\sum}_{j=0}^N\, x^j\left(\frac{d^{k+j}}{dx^{k+j}}\widehat{h_j}\right)(x)$, hence
\begin{equation}
\label{fourier8}
\widehat\sigma(1/x) \ =\  (-1)^{k+j}\,(2\pi i)^{-k}\,{\sum}_{0 \leq j \leq N}\ x^{-j}\left(-x^2\textstyle\frac{d\ }{dx}\right)^{k+j}\left(\widehat{h_j}(1/x)\right)\,.
\end{equation}
The function $\,x\mapsto\widehat{h_j}(1/x)$ is bounded, and is smooth except at $0$, hence locally $L^\infty$. The differential operator $\,x^{-j}(-x^2\textstyle\frac{d\ }{dx})^{k+j}$ can be expressed as a linear combination of $\,x^{k+i}\frac{d^i\ }{dx^i}$, $0\leq i \leq k+j$. We have therefore extended $\,\widehat\sigma$ across $\infty$, where the extension vanishes to order $k$.
\end{proof}

According to our convention, distributions are dual to compactly supported smooth measures. A distribution $\,\sigma$ defined on some neighborhood $U$ of $\infty$ in the compactified real line $\,\RP^1 = \R \cup \{\infty\}$ can be integrated against a smooth measure supported in $U$, or equivalently, against $g(1/x)\,dx$, where $g \in C^\infty_c(1/U)$ must vanish to second order at $0$, to balance the second order pole of $\,dx$ at $\infty$. If $\sigma$ vanishes to second order at $\infty$, it can absorb the second order pole. Thus, in this situation, the change of variables formula
\begin{equation}
\label{changeofvariables}
\int_U \sigma(x)\,g(1/x)\,dx \ = \ \int_{1/U} x^{-2}\,\sigma(1/x)\,g(x) \,dx
\end{equation}
can be legitimately applied to any $g \in C^\infty_c(1/U)$, without requiring $g$ to vanish at $0$.

We now suppose that $\,\sigma \in \Sch'(\R)$ vanishes to order $k_0\geq 0$ at $0$ and has an extension across $\infty$ which vanishes there to order $k_\infty\geq 0$. We also suppose that $k_0+k_\infty \geq 1$. According to proposition \ref{xalpha}, for $\delta\in\Z/2\Z$ and $\operatorname{Re}s > -k_0$, the distribution $(\sg x)^\delta |x|^{s-1} \sigma(x)$ has an extension across $0$ which depends holomorphically on $s$. Similarly, there exists an extension across $\infty$, with holomorphic dependence on $s$, for $\operatorname{Re}s < k_\infty+2$. Equivalently, $(\sg x)^\delta |x|^{s-1} \sigma(x)dx$ is well defined and holomorphic near $x=\infty$ when $\operatorname{Re}s < k_\infty$, since $dx$ has a second order pole at $\infty$. Putting the two statements together, we see that $(\sg x)^\delta |x|^{s-1} \sigma(x)dx$, with $-k_0< \operatorname{Re}s < k_\infty$, can be regarded as global ``measure with distribution coefficients" on the compact manifold $\R\cup\{\infty\}$. As such, it can be integrated against the constant function $1$. This allows us to define the signed Mellin transform of $\,\sigma$,
\begin{equation}
\label{mellin1}
M_\delta\sigma\,(s) \ = \ \int_{\R} (\sg x)^\delta \,|x|^{s-1}\, \sigma(x)\,dx \qquad \bigl(\, -k_0 < \operatorname{Re}s < k_\infty \,\bigr)\,,
\end{equation}
as a holomorphic function of $s$ in the indicated region. The notation takes some license: we are really integrating over the compactified real line, and the integrand needs to be extended across $\infty$ in the described manner. It is sometimes convenient to split up the integral (\ref{mellin1}) into two integrals over bounded intervals. For that purpose, we choose a cutoff function $\,\phi\in C^\infty_c(\R)$ which is identically equal to $1$ near $x=0$. Then $\,\sigma = \phi\,\sigma + (1-\phi)\,\sigma$, and
\begin{equation}
\label{mellin2}
\begin{aligned}
&M_\delta\sigma\,(s) \ = \ \int_{\R} (\sg x)^\delta \,|x|^{s-1}\,\phi(x)\, \sigma(x)\,dx\ +
\\
&\qquad\qquad\qquad\qquad+  \ \int_{\R} (\sg x)^\delta \,|x|^{-s-1}\, (1 - \phi)(1/x)\,\sigma(1/x)\,dx\,.
\end{aligned}
\end{equation}
Both integrals have compactly supported integrands, and thus can be regarded as integrals over bounded intervals. The change of variables formula (\ref{changeofvariables}), which justifies the passage from (\ref{mellin1}) to (\ref{mellin2}), can also be used to switch the roles of $0$ and $\infty$\,:
\begin{equation}
\label{mellin3}
M_\delta\sigma\,(s) \, = \, M_\delta\widetilde\sigma\,(-s) \ \ \ \bigl( -k_0 < \operatorname{Re}s < k_\infty\bigr)\,,\ \ \ \text{with}\ \ \widetilde\sigma(x)\, = \, \sigma(1/x)\,.
\end{equation}
The relation between $\,\sigma$ and $\,\widetilde\sigma$ is to be understood as an equality of distri\-butions on $\R\cup\{\infty\}$\,, of course.

Recall the generalization of the notion of vanishing to order $k$ introduced by remark \ref{negative_order}: $\,\sigma$ vanishes at $0$ to order $-k$, $k\geq 1$, if $x^k \sigma$ vanishes to order zero. As in \defref{inf_order}, vanishing to order $k$ implies vanishing to order $\ell$, for any integer $\ell\leq k$. With the extended definition, the discussion of the Mellin transform still applies. More precisely, we can define the signed Mellin transform $M_\delta \sigma$ of a tempered distribution $\,\sigma$ which vanishes to order $k_0\in \Z$ at $0$ and has an extension across $\infty$ vanishing there to order $k_\infty\in\Z$, provided ${k_0+k_\infty\geq 1}$. In this situation, $M_\delta \sigma$ is defined as holomorphic function on the region $\{-k_0 < \operatorname{Re}s < k_\infty\}$, and the identities (\rangeref{mellin2}{mellin3}) remain valid.  A word of caution: if $k_0 + k_\infty \leq 0$ even with maximal choices for $k_0$ and $k_\infty$, it is still possible for the two integrals in (\ref{mellin2}) to have analytic continuations whose domains intersect on some vertical strip. In that case, it is not legitimate to think of the sum of the two integrals, on their common domain of definition, as the Mellin transform of $\,\sigma$.

Any distribution $\,\sigma$, defined on some neighborhood of $0$, can be expressed locally as a $k$-th derivative of some continuous function, in which case $\,\sigma$ vanishes to order $-k$ according to the generalized definition. If $\,\sigma$ also has a canonical extension to $\R\cup\{\infty\}$, $M_\delta \sigma$ is then defined on the right half plane $\{k < \operatorname{Re}s\}$. Similarly, if $\,\sigma\in\Sch'(\R)$ vanishes to infinite order at $0$, the signed Mellin transform $M_\delta \sigma$ is defined on some left half plane $\{ \operatorname{Re}s < -k\}$. Lastly, if $\,\sigma$ both vanishes to infinite order at $0$ and has a canonical extension across $\infty$, $M_\delta \sigma$ is defined as an entire function.

For our next statement, we consider a tempered distribution $\,\sigma\in\Sch'(\R)$ which vanishes to order $k_0$ at $0$ and has an extension across $\infty$ which vanishes there to order $k_\infty$, with $k_0 + k_\infty \geq 1$. The product $x^{-1}\sigma(x)$ vanishes to order at least $k_0-1$ at $0$ and has an extension to $\R\cup\{\infty\}$ which vanishes at $\infty$ to order at least $k_\infty +1$. In this situation, $M_\delta\sigma\,(s)$ and $M_{\delta+1}(x^{-1}\sigma)(s+1)$ are both defined on the region $\{-k_0 < \operatorname{Re}s < k_\infty\}$. Differentiation also has this effect on $k_0$ and $k_\infty$, so $M_{\delta+1}\sigma'\,(s+1)$ is defined on the same region.

\begin{prop}\label{mellin_properties}
Under the hypotheses that were just mentioned,
\begin{equation*}
M_{\delta+1}\sigma'\,(s+1) \ = \ -s\,M_\delta\sigma\,(s) \ = \ -s\,M_{\delta+1}(x^{-1}\sigma)\,(s+1)\,,
\end{equation*}
as equalities of holomorphic functions in the region $\{-k_0 < \operatorname{Re}s < k_\infty\}$. In particular, $M_\delta\sigma'\,(1)=0$ if $-k_0<0<k_\infty$.
\end{prop}

\begin{proof}
We can argue separately for the two summands in (\ref{mellin2}). In other words, we may suppose that $\,\sigma$ has compact support and vanishes at $0$ to order $k_0$. The second identity follows formally from the definition, and we may use it freely in the proof of the first identity. In particular, we may suppose $k_0=1$. Then, according to \lemref{one_var_def} and \lemref{compactsupport},
\begin{equation}
\label{mellin_properties1}
\sigma \ = \ {\sum}_{0\leq j\leq N}\ x^j\,\textstyle\frac{d^j\ }{dx^j}\,h_j\,,\ \ \ \text{with}\ \ h_j \in C_c(\R)\,.
\end{equation}
Taking one summand at a time and using the identity we already know reduces the problem to the case of a distribution $\,\sigma=h^{(j)}$ which arises as the $j$-th derivative of a function $h\in C_c(\R)$. According to the definition of the signed Mellin transform,
\begin{equation}
\label{mellin_properties2}
M_\delta h^{(j)}\,(s) \ = \ \biggl(\, {\prod}_{1\leq \ell \leq j}\ (\ell-s) \biggr) \int_\R (\sg x)^{\delta+j}\,|x|^{s-j-1}\,h(x)\,dx
\end{equation}
in the region $\{\operatorname{Re}s > j\}$. When we substitute $j+1$ for $j$, the identity we want follows for $\,\sigma=h^{(j)}$ and hence any compactly supported $\,\sigma$, at least when $\operatorname{Re}s \gg 0$. If $\,\sigma$ has compact support and vanishes to order $k_0$ at $0$, as we had assumed, both sides of the equation are known to be holomorphic to the right of the line $\operatorname{Re}s=-k_0$. The equation to be proved extends to this region by analytic continuation.
\end{proof}

Some examples may be instructive at this point. Dirac's function $\delta_0$ based at the origin vanishes there to order $-1$, since it can be written, locally near $0$, as the first derivative of a bounded measurable function. That is the lowest possible value for the order of vanishing, since \lemref{support} rules out order of vanishing zero. The signed Mellin transform of $\delta_0$ has meaning as holomorphic function on the region $\{\operatorname{Re}s > 1\}$. Parity considerations show that the odd Mellin transform $M_1\delta_0$ vanishes identically. But $x\delta_0 = 0$ and $M_0\delta_0\,(s)=M_1(x\delta_0)\,(s-1)$ by \propref{mellin_properties}, so the even Mellin transform $M_0\delta_0$ also vanishes. One can argue similarly for the derivatives of $\delta_0$, which together with $\delta_0$ span the space of distributions supported at the origin, hence:

\begin{cor}\label{mellin_of_dirac}
If $\,\sigma\in\Sch'(\R)$ has support at the origin and vanishes there to order $-k$, $M_\delta\sigma\,(s)=0$, as an identity on the region $\{\operatorname{Re}s > k\}$.
\end{cor}

The signed Mellin transform $M_\delta f$ of a Schwartz function $f$, of parity $\delta$, equals twice the usual Mellin transform of $f$. Thus well known results about the usual Mellin transform can be restated as follows:

\begin{lem}\label{mellin_of_schwartz} The signed Mellin transform $M_\delta f(s)$ of a Schwartz function $f\in \Sch_\delta(\R)$ extends meromorphically to the complex plane, with poles only at integral points $\,s=-n$, $\,n\geq 0$, $\,n\equiv \delta$ modulo $2$, all of first order, with residue $2\,(n!)^{-1}\,f^{(n)}(0)$ at $s=-n$. The Mellin transforms of $f$ and $f'$ are related by the identity $M_{\delta+1}f'\,(s+1) \ = \ -s\,M_\delta f\,(s)$
\end{lem}

If a function $f\in C^\infty(\R)$ vanishes to order $k\geq 0$ at the origin, it also vanishes to order $k$ in the sense of definition \ref{inf_order} when $f$ is regarded as a distribution. The converse of this statement is not equally obvious, but can be deduced from properties of the Mellin transform. To see this, we choose a cutoff function $\,\phi\in C_c^\infty(\R)$ which is identically equal to $1$ on some neighborhood of the origin. If $f$ vanishes to order $k$ in the sense of definition \ref{inf_order}, so does $\,\phi f$, which vanishes near $\infty$, hence extends canonically across $\infty$. That makes $M_\delta (\phi f)$ well defined and holomorphic for $\operatorname{Re}s > -k$, for both choices of $\delta$; in particular, $M_\delta (\phi f)$ has no poles at $s=-n$, $0\leq n<k$. In the region $\{\operatorname{Re}s > 0\}$, the definition of the Mellin transform of $\,\phi f$, viewed as distribution, agrees with the definition of $M_\delta (\phi f)$ when $\,\phi f$ is regarded as Schwartz function. At this point, \lemref{mellin_of_schwartz} implies:

\begin{cor}\label{equiv_of_def}
For a function $f\in C^\infty(\R)$, vanishing to order $k$ at $0$ in the usual sense is equivalent to vanishing to order $k$ in the sense of distributions, according to \defref{inf_order}.
\end{cor}

The corollary applies in particular to the constant function $1$, which has order of vanishing $k_0=k_\infty=0$ both at the origin and at $\infty$. Our definition of the signed Mellin transform does not apply directly since $k_0+k_\infty=0$, but there is a weaker notion, as we shall explain next. We choose a cutoff function $\,\phi\in C_c(\R)$ such that $\,\phi(x)\equiv 1$ near $x=0$. Then $1 = \phi + (1-\phi)$, and $M_\delta\phi\,(s)$, $M_\delta(1-\phi)\,(s)$ are defined on the right and left half plane, respectively.

\begin{lem}\label{mellin_of_one}
If $\,\phi\in C_c(\R)$, $\,\phi(x)\equiv 1$ near the origin, the signed Mellin transform $M_\delta\phi\,(s)$ extends meromorphically from $\{\operatorname{Re}s > 0\}$ to the entire complex plane, with at most a single pole at $s=0$, of order one and residue $2$ or $0$, depending whether $\delta=0$ or $\delta=1$. Similarly $M_\delta(1-\phi)\,(s)$ extends meromorphically from $\{\operatorname{Re}s < 0\}$ to $\C$. The sum of these two meromorphic continuations vanishes identically.
\end{lem}

\begin{proof}
The description of the poles and residues follows from \lemref{mellin_of_schwartz}, which also provides an explicit meromorphic continuation via the equation
\begin{equation}
\label{mellin_of_schwartz1}
M_{\delta+1}\phi'\,(s+1) \ = \ -s\,M_\delta\phi\,(s)\,.
\end{equation}
What matters here is the entirety of $M_{\delta+1}\phi'\,(s+1)$, which follows from the fact that $\,\phi'$ has compact support and vanishes identically near $x=0$. According to (\ref{mellin3}),
\begin{equation}
\label{mellin_of_schwartz2}
M_\delta(1-\phi)\,(s) \ = \ M_\delta\psi\,(-s)\ \ \ \bigl(\,\operatorname{Re}s < 0\,\bigr)\,,\ \ \text{with}\ \ \psi(x)=1-\phi(1/x)\,.
\end{equation}
Like $\,\phi$, the function $\,\psi$ has compact support and vanishes near $x=0$, hence
\begin{equation}
\label{mellin_of_schwartz3}
M_{\delta+1}\psi'\,(1-s) \ = \ s\,M_\delta\psi\,(-s)\,,
\end{equation}
is an entire function, in analogy to (\ref{mellin_of_schwartz1}). To complete the proof, it suffices to establish the equality of entire functions
\begin{equation}
\label{mellin_of_schwartz4}
M_{\delta+1}\phi'\,(s+1) \ = \ M_{\delta+1}\psi'\,(1-s)\,,
\end{equation}
which follows from (\ref{mellin3}) and the identity $\,\psi'(1/x)=x^2\phi'(x)$.
\end{proof}

As the final example of this section, we consider a periodic distribution without constant term $\,\tau(x)=\sum_{n\neq 0} \, a_n\, e(nx)$. Appealing either to \propref{oscillation} or \thmref{fourier}, we can conclude that $\,\tau$ has a canonical extension across infinity. Its signed Mellin transform is then well defined on some half plane $\{\,\operatorname{Re}s > k\,\}$. The Fourier coefficients $a_n$ grow at most polynomially with $n$, so the Dirichlet series $\sum_{n>0}\,a_n\,n^{-s}$ and $\sum_{n>0}\,a_{-n}\,n^{-s}$ converge absolutely for $\,\operatorname{Re}s \gg 0$.

\begin{lem}\label{mellin_and_dirichlet}
The signed Mellin transform of $\,\tau = \sum_{n\neq 0} a_n\,e(nx)$ is given by the formula
\[
M_\delta\tau\,(s) \ \ = \ \ G_\delta(s)\ {\sum}_{n\neq 0}\ (\sg n)^\delta \, a_n \,|n|^{-s}\ \ \ \text{for}\ \ \operatorname{Re}s \gg 0\,,
\]
with $\,G_0(s)=2(2\pi)^{-s}\Gamma(s)\cos(\pi s/2)$ and $\,G_1(s)=2i(2\pi)^{-s}\Gamma(s)\sin(\pi s/2)$.
\end{lem}

\begin{proof} To shorten the various formulas we only discuss the case $\delta=0$; the other case can be treated exactly the same way. We choose $k_0\in\N\,$ large enough to ensure $\sum_{n\neq 0} |a_n|\,|n|^{-k_0}<\infty$. Then, for $k\geq k_0$,
\begin{equation}
\label{mellin_and_dirichlet1}
F_k(x) \ = \ \textstyle\sum_{n \neq 0}\, a_n \, (2\pi i n)^{-k}\, e(nx)
\end{equation}
converges absolutely and uniformly to a continuous function $\,F_k(x)$, such that $\,F_k^{(k)}(x) = \tau(x)$. As before, we pick a cutoff function $\phi \in C^\infty_c(\R)$, with $\phi(x)\equiv 1$ for $x$ near $0$. Then, in analogy to (\ref{mellin2}),
\begin{equation}
\label{mellin_and_dirichlet2}
\begin{aligned}
&M_0\tau\,(s) \ \ = \ \ \int_{\R} |x|^{s-1}\, \phi(x)\,\tau(x)\,dx
\\
&\qquad\qquad\qquad\qquad + \ \int_{\R} |x|^{-s-1}\, (1 - \phi)(1/x)\,\tau(1/x)\,dx\,.
\end{aligned}
\end{equation}
The first integral on the right is to be interpreted as
\begin{equation}
\label{mellin_and_dirichlet3}
\begin{aligned}
\int_{\R} |x|^{s-1}\, \phi(x)\,\tau(x)\,dx\ &= \ \int_{\R} |x|^{s-1}\, \phi(x)\,F_{k_0}^{(k_0)}(x)\,dx
\\
&= \ (-1)^{k_0}\int_{\R} \textstyle\frac{d^{k_0}\ }{dx^{k_0}}\bigl(|x|^{s-1}\, \phi(x)\bigr)\,F_{k_0}(x)\,dx\,,
\end{aligned}
\end{equation}
which has meaning for $\,\operatorname{Re}s > k_0 + \epsilon$, for any $\epsilon > 0$, since then
\begin{equation}
\label{mellin_and_dirichlet4}
\textstyle\frac{d^{k_0}\ }{dx^{k_0}}\,\bigl(\,|x|^{s-1}\, \phi(x)\,\bigr)\ \ = \ \ O(\,|x|^{\epsilon-1})\ \ \ \ \bigl(\, |x|\ll 1  \,\bigr)\,.
\end{equation}
We substitute the uniformly convergent series (\ref{mellin_and_dirichlet1}) in this identity, interchange the order of summation and integration, then reverse the integration by parts, and conclude
\begin{equation}
\label{mellin_and_dirichlet5}
\int_{\R}  |x|^{s-1}\, \phi(x)\, \tau(x)\,dx \ \ =\ \ {\sum}_{n \neq 0}\,\ a_n  \int_{\R} |x|^{s-1}\, \phi(x)\, e(nx)\,dx\,,
\end{equation}
still for $\,\operatorname{Re}s > k_0 + \epsilon$.

We now choose $k_1\in\N$, $k_1 > k_0$, and suppose that $\,\operatorname{Re}s < k_1 - \epsilon$, for some $\epsilon > 0$. To paraphrase the proof of \propref{xalpha} in light of \lemref{one_var_def}, one extends the distribution $\,|x|^{-s-1}\tau(1/x)$ across $x=0$ by means of the formula
\begin{equation}
\label{mellin_and_dirichlet6}
\begin{aligned}
&|x|^{-s-1}\, \tau(1/x) \ \ =\ \ |x|^{-s-1}\, \bigl( -x^2 \textstyle\frac{d\ }{dx} \bigr)^{k_1}\bigl(F_{k_1}(1/x)\bigr)
\\
&\qquad\qquad=\ {\sum}_{0\leq j\leq k_1}\ q_j(s)\,\textstyle\frac{d^j\ }{dx^j}\bigl( (\sg x)^{k_1+j}\,|x|^{-s+k_1+j-1}\,F_{k_1}(1/x) \bigr)\,,
\end{aligned}
\end{equation}
with suitable polynomials $q_j(s)$; this depends on the fact that
\begin{equation}
\label{mellin_and_dirichlet7}
(\sg x)^{k_1+j}\,|x|^{-s+k_1+j-1}\,F_{k_1}(1/x)\ \ = \ \ O(|x|^{\epsilon-1})\ \ \ \ \bigl(\, |x|\ll 1  \,\bigr)\,,
\end{equation}
because $F_{k_1}(1/x)$ is bounded and $\,\operatorname{Re}s < k_1 - \epsilon$. We substitute the series (\ref{mellin_and_dirichlet1}) for $F_{k_1}$ in (\ref{mellin_and_dirichlet6}) and interchange the order of differentiation and summation; that is legitimate because the series converges absolutely. When we integrate the resulting formula against the smooth -- even at $\,x=0\,$! -- compactly supported function $(1-\phi)(1/x)$ and work backwards, we obtain the formula
\begin{equation}
\label{mellin_and_dirichlet8}
\begin{aligned}
&\int_{\R} |x|^{-s-1}\, (1 - \phi)(1/x)\,\tau(1/x)\,dx\ \ =
\\
&\qquad\qquad=\ \ {\sum}_{n\neq 0}\ a_n\, \int_{\R} |x|^{-s-1}\, (1 - \phi)(1/x)\,e(n/x)\,dx\,.
\end{aligned}
\end{equation}
The integral on the right is not an ordinary integral, but rather denotes the integral of the canonical extension of the distribution $\,|x|^{-s-1}e(n/x)$ against the smooth, compactly supported measure $(1 - \phi)(1/x)\,dx$.

We now make the change of variables $\,x\rightsquigarrow x/n$ in (\ref{mellin_and_dirichlet5}), the change of variables $\,x\rightsquigarrow n/x$ in (\ref{mellin_and_dirichlet8}), and combine the resulting formulas with (\ref{mellin_and_dirichlet2}):
\begin{equation}
\label{mellin_and_dirichlet9}
M_0\tau\,(s) \ \ = \ \ {\sum}_{n\neq 0}\ |n|^{-s}\,a_n\int_{\R} |x|^{s-1}\, e(x)\,dx\,.
\end{equation}
In deriving this identity, we have assumed that $\,k_0 + \epsilon < \operatorname{Re}s< k_1 - \epsilon$, which in particular implies absolute convergence of the series. The integral is an ordinary, convergent integral near $x=0$, but near $x=\infty$, it must be interpreted as the integral of the canonical extension of the distribution $|x|^{s+1} e(x)$ against the smooth measure $x^{-2} dx$. Taken in this sense, the integral represents a holomorphic function at least on the region $\{\,\operatorname{Re}s > 0\,\}$. But this function is well known: on the subregion $\{\,0<\operatorname{Re}s <1\,\}$, the integral converges conditionally, to the limit $G_0(s)$. The lemma follows.
\end{proof}

\section{Fourier And Mellin Transforms}\label{4nmelsec}

In this section we use the signed Mellin transform to characterize tempered distributions which vanish to infinite order at $0$ and extend canonically across $\infty$, and we relate the Fourier and Mellin transforms of such distributions.

We start out by introducing the two main results. For $\delta\in\Z/2\Z$, we let $\Sch_\delta(\R)$ denote the space of Schwartz functions of parity $\delta$,
\begin{equation}
\label{parity1}
\Sch_\delta(\R) \ \ = \ \ \{\,f\in \Sch(\R) \, \mid \, f(-x)\,=\,(-1)^\delta f(x) \}\,,
\end{equation}
and similarly $\Sch'_\delta(\R)$ the space of tempered distributions of parity $\delta$. Then
\begin{equation}
\label{parity2}
\Sch(\R) \ \ = \ \ \Sch_0(\R)\oplus\Sch_1(\R)\,,\ \ \ \ \Sch'(\R) \ \ = \ \ \Sch'_0(\R)\oplus\Sch'_1(\R)\,.
\end{equation}
When it is defined, the signed Mellin of a tempered distribution $\,\sigma$ satisfies the relation
\begin{equation}
\label{parity3}
M_\delta\widetilde{\sigma} \ = \ (-1)^\delta \,M_\delta\sigma\ \ \ \ \ \text{if}\ \ \ \ \widetilde\sigma(x)\,=\,\sigma(-x)\,,
\end{equation}
for entirely formal reasons. In particular, the even Mellin transform of an odd distribution vanishes identically, as does the odd Mellin transform of an even distribution.

We shall say that a holomorphic function $H(s)$, defined on the region $\{\,a<\operatorname{Re}s<b\,\}$, $\,-\infty\leq a<b\leq\infty\,$, has ``locally uniform polynomial growth" on vertical lines if
\begin{equation}
\label{moderategrowth1}
\begin{aligned}
&|H(s)| \ = \ O(|s|^N)\ \ \ \text{as}\ \ |\operatorname{Im}s|\to\infty\,,
\\
&\qquad\qquad\text{for some $N\in\N$, locally uniformly in}\,\ \operatorname{Re}s\,.
\end{aligned}
\end{equation}
The integer $N$ may depend on $\operatorname{Re}s$, but must do so in a locally uniform manner. Equivalently, for each choice of $a'$, $b'$, with $a<a'<b'<b$, there must exist positive constants $T=T(a',b')$, $C=C(a',b')$, and $N=N(a',b')\in\N$, such that
\begin{equation}
\label{moderategrowth2}
a'\leq \operatorname{Re}s \leq b'\,,\ \ |\operatorname{Im}s|\geq T\ \ \ \Longrightarrow\ \ \ |H(s)|\leq C\,|\operatorname{Im}s|^N \,.
\end{equation}
For this reason, we shall also use the synonymous terminology ``$H(s)$ has moderate growth on vertical strips". We refer to the dual notion,
\begin{equation}
\label{rapiddecay1}
\begin{aligned}
&|H(s)| \ = \ O(|s|^{-N})\ \ \ \text{as}\ \ |\operatorname{Im}s|\to\infty\,,
\\
&\qquad\qquad\text{for every $N\in\N$, locally uniformly in}\,\ \operatorname{Re}s\,,
\end{aligned}
\end{equation}
by saying that ``$H(s)$ has locally uniform rapid decay" along vertical lines, or synonymously, that ``$H(s)$ decays rapidly on vertical strips". This condition holds if and only if, for all $a',\,b'\in\R$, with $a<a'<b'<b$, and every $N\in\N$, there exist positive constants $T=T(a',b',N)$, $C=C(a',b',N)$, such that
\begin{equation}
\label{rapiddecay2}
a'\leq \operatorname{Re}s \leq b'\,,\ \ |\operatorname{Im}s|\geq T\ \ \ \Longrightarrow\ \ \ |H(s)|\leq C\,|\operatorname{Im}s|^{-N} \,.
\end{equation}
The conditions (\rangeref{moderategrowth1}{rapiddecay2}) make sense even when $H(s)$ is meromorphic, provided the real parts of the poles of $H(s)$ have no accumulation points in the open interval $(a,b)$ -- in particular, when all the poles lie on the real line. We shall use the same terminology in that situation.

\begin{thm}\label{melthm}
For $\delta\in\Z/2\Z$, the Mellin transform $M_\delta$ establishes an isomorphism between
\[
\{\sigma\!\in\!\Sch'_\delta(\R) \mid \sigma\ \text{vanishes to infinite order at $0$, extends canonically across $\infty$}\}
\]
and
\[
\{\,H : \C \to \C\, \mid \, \text{$H$ is entire, of moderate growth on vertical strips}\,\}\,.
\]
\end{thm}

Recall that the signed Mellin transform $M_\delta\sigma(s)$ is well defined and regular for $\operatorname{Re} s \ll 0$ provided $\,\sigma\in\Sch'(\R)$ vanishes to infinite order at $x=0$, whereas $M_\delta\sigma(s)$ is well defined and regular for $\operatorname{Re} s \gg 0$ when $\sigma$ has a canonical extension across infinity. According to theorem \ref{fourier}, if $\,\sigma\in\Sch'(\R)$ vanishes to infinite order at $x=0$, the Fourier transform $\widehat\sigma$ extends canonically across infinity. In that situation, the domains of definition of $M_\delta\sigma(1-s)$ and $M_\delta\widehat\sigma(s)$ intersect in some half plane $\operatorname{Re} s \gg 0$. Our next statement also involves the function
\begin{equation}
\label{whyGdelta}
G_\delta(s)\ \ = \ \ \int_\R e(x)\,(\sg x)^\delta\,|x|^{s-1}\,dx\ \ \ \ \ (\,0 < \operatorname{Re} s < 1\,)\,.
\end{equation}
It extends meromorphically to the entire complex plane by virtue of the formula
\begin{equation}
\label{whyGdelta1}
G_\delta(s)\ \ = \ \ 2\,i^\delta\,(2\pi)^{-s}\,\Gamma(s)\,\cos(\pi(s-\delta)/2)\,,
\end{equation}
which follows from standard identities for the Mellin transform. Note that the current definition is consistent with the earlier use of the notation $G_\delta(s)$ in lemma \ref{mellin_and_dirichlet}. We shall also use the functional equation
\begin{equation}
\label{whyGdelta2}
G_\delta(s)\,G_\delta(1-s) \ \ = \ \ (-1)^\delta\,,
\end{equation}
which is equivalent to the Gamma identity $\Gamma(s)\Gamma(1-s)=\pi \csc(\pi s)$.

\begin{thm}\label{4nmelthm}
{\rm a)\,} If $\,\sigma\in\Sch'(\R)$ vanishes to infinite order at the origin, $M_\delta\widehat\sigma(s)\, = \, (-1)^\delta\, G_\delta(s)\,M_\delta\sigma\,(1-s)\ \ \text{for}\,\ \operatorname{Re} s \gg 0\,$.\newline
\noindent {\rm b)\,} Suppose that $\,\sigma\in\Sch_\delta'(\R)$ extends canonically across infinity and vanishes to order $k_0\geq 1$ at the origin. Then $\,\widehat\sigma$ vanishes to infinite order at the origin if and only if $\,G_\delta(s)\,M_\delta\sigma\,(1-s)\,$ is regular for $\,\operatorname{Re} s < 1$. When that is the case, $M_\delta\widehat\sigma(s)\, = \, (-1)^\delta\,G_\delta(s)\, M_\delta\sigma(1-s)\,$ for $\, \operatorname{Re} s < 1$.
\end{thm}

A periodic distribution without constant term $\,\tau(x)=\sum_{n\neq 0}\,a_ne(nx)\,$ has Fourier transform $\,\widehat\tau(x)=\sum_{n\neq 0}\,a_n \,\delta_n(x)$, which vanishes identically near the origin. Since $M_\delta\widehat\tau(s)=\sum_{n\neq 0}(\sg n)^\delta a_n |n|^{1-s}$, and since $\tau(-x)$ is the double Fourier transform of $\,\tau$, the theorem contains lemma \ref{mellin_and_dirichlet} as a special case.

The remainder of this section contains the proofs of theorems \ref{melthm} and \ref{4nmelthm}, which depend on similar arguments.

Recall \lemref{mellin_of_schwartz}, which relates the Mellin transform $M_\delta f$ of a Schwartz function $f$ to that of its derivative, and which asserts that $M_\delta f$ is regular on the complex plane, except for first order poles at non-positive integers of parity $\delta$.

\begin{lem}\label{mellem1} For every choice of $\epsilon,\,R\in\R$, with $0<\epsilon<R$, there exists a continuous seminorm $\nu_{\epsilon,R}:\Sch(\R)\to \R_{\geq 0}$ such that
\[
f\in\Sch(\R)\,,\ \ \epsilon\leq \operatorname{Re}s\leq R\ \ \ \Longrightarrow\ \ \ |M_\delta f(s)|\leq \nu_{\epsilon,R}(f)\,.
\]
\end{lem}

\begin{proof}
The family of seminorms $\mu_{n,k}(f)={\sup}_{x\in\R}\bigl((1+x^2)^n|f^{(k)}(x)|\bigr)$, indexed by integers $k,n\geq 0$, defines the topology of $\Sch(\R)$. We choose $n$ so that $2n> R$. Then, for $\epsilon \leq \operatorname{Re}s \leq R$,
\begin{equation}
\label{melbound1}
\begin{aligned}
&|M_\delta f(s)| \ =\ \bigl|\,\int_\R (\sg x)^\delta\,|x|^{s-1}\,f(x)\,dx \bigr|\ \leq \ \int_\R |x|^{\operatorname{Re}s-1}\,|f(x)|\,dx
\\
&\ \ \leq \ \mu_{n,0}(f) \int_\R (1+x^2)^{-n}\,|x|^{\operatorname{Re}s-1}\,dx\ \leq \ \biggl(\frac2\epsilon + \frac2{(2n-R)}\biggr)\,\mu_{n,0}(f)\,.
\end{aligned}
\end{equation}
Thus $\,\nu_{\epsilon,R}=2\bigl(\epsilon^{-1}+(2n-R)^{-1}\bigr)\mu_{n,0}$ will do.
\end{proof}

In particular, $M_\delta f(s)$ decays rapidly on vertical strips to the right of $\operatorname{Re}s=0$. Repeated application of the relation $-sM_\delta f\,(s)=M_{\delta+1} f'\,(s+1)$ implies
\begin{equation}
\label{melbound2}
M_\delta f\,(s)\ = \ (-1)^m\bigl(\, {\prod}_{0\leq \ell \leq m-1}\ (s+\ell)^{-1}\,\bigr)\, M_{\delta+m}f^{(m)}(s+m)\,,
\end{equation}
so the boundedness of the Mellin transform of $f^{(m)}\in\Sch(\R)$ to the right of the line $\operatorname{Re}s=0$ forces the bound $|M_\delta f\,(s)|=O(|\operatorname{Im}s|^{-m})$ to the right of $\operatorname{Re}s=-m$, for any $m\in\N$.

\begin{cor}\label{melcor1}
The signed Mellin transform $M_\delta f(s)$ of a Schwartz function $f$ decays rapidly on vertical strips.
\end{cor}

This is equivalent, of course, to the same, and well known, property of the ordinary Mellin transform. For $k_0,\,k_\infty\in\Z\cup\{+\infty\}$, we define
\begin{equation}
\label{k0kinf}
\begin{aligned}
\Sch'(\R)_{k_0,\,k_\infty}\ = \ \bigl\{\,&\sigma\in\Sch'(\R)\,\mid\,\sigma\,\ \text{vanishes to order $k_0$ at $0$, has}
\\
&\text{extension across $\infty$ vanishing to order $k_\infty$ at $\infty$}\bigr\}\,.
\end{aligned}
\end{equation}
We use the natural convention $\,k_0+k_\infty=+\infty\,$ when at least one of the {summands} has the value $+\infty$. Recall that the signed Mellin transform $M_\delta\sigma(s)$ of a distribution $\,\sigma\in\Sch'(\R)_{k_0,\,k_\infty}$ is well defined and lies in the function space
\begin{equation}
\label{Ok0kinf}
\begin{aligned}
&\O(\{-k_0<\operatorname{Re}s<k_\infty\})\ =
\\
&\qquad\qquad= \ \bigl\{\,H:\{-k_0<\operatorname{Re}s<k_\infty\}\to\C\,\mid \, \text{$H$ is holomorphic}\bigr\}\,,
\end{aligned}
\end{equation}
provided $k_0+k_\infty\geq 1$.

\begin{lem}\label{moderategrowth} If $\,\sigma\in\Sch'(\R)_{k_0,\,k_\infty}$\,, with $k_0+k_\infty\geq 1$, $\,M_\delta\sigma(s)$ has locally uniform polynomial growth on vertical lines.
\end{lem}

\begin{proof}
We express $M_\delta\sigma$ as the sum of two terms, as in (\ref{mellin2}). They are of the same type, so it suffices to show that the first term has locally uniform polynomial growth on vertical lines to the right of $\operatorname{Re}s=-k_0$. Changing notation from $\,\phi\,\sigma$ to $\,\sigma$, we may and shall suppose that $\,\sigma$ has compact support and vanishes to order $k_0$ at $0$. We appeal to lemmas \ref{one_var_def} and \ref{compactsupport}, to write
\begin{equation}
\label{moderategrowth4}
\begin{aligned}
\sigma(x)\ =\ {\sum}_{0\leq j\leq N}\ x^{k_0+j}\,\textstyle\frac{d^j\ }{dx^j}\,h_j(x)\ \ \ \text{with}\ \ \ h_j \in L^\infty(\R)
\\
\text{and}\ \ \ h_j(x)\equiv 0 \ \ \ \text{for}\ \ |x|\gg 1\,,
\end{aligned}
\end{equation}
which then implies
\begin{equation}
\label{moderategrowth5}
\begin{aligned}
&\int_\R (\sg x)^\delta\,|x|^{s-1}\,\sigma(x)\,dx\ \ =
\\
&\ \ = \ {\sum}_{0\leq j\leq N}\ (-1)^j\int_\R h_j(x)\,\textstyle\frac{d^j\ }{dx^j}\bigl((\sg x)^{\delta+k_0+j}\,|x|^{s+k_0+j-1}\bigr)\,dx\\
&\ \ = \ {\sum}_{0\leq j\leq N}\, {\prod}_{0\leq\ell< j}(-s - k_0 - \ell)\!\int_\R \!\! h_j(x)\textstyle(\sg x)^{\delta+k_0}|x|^{s+k_0-1}dx\,.
\end{aligned}
\end{equation}
For $0<\epsilon \ll 1\,$, $\,R\gg 0\,$, and $\,-k_0+\epsilon\leq\operatorname{Re}s\leq R$, the integral on the right can be bounded in terms of $\,\epsilon$, $R$, the supremum of the $|h_j|$ and the support of the $h_j$, entirely without reference to $\,\operatorname{Im}s$.
\end{proof}

To paraphrase \lemref{moderategrowth}, the signed Mellin transform $M_\delta$ defines a linear map from $\Sch'(\R)_{k_0,\,k_\infty}$ to the subspace
\begin{equation}
\label{Ok0kinfmg}
\begin{aligned}
\!\!\O_{\text{pg}}(\{-k_0<\operatorname{Re}s<k_\infty\}) = \text{space of all}\ H\in \O(\{-k_0<\operatorname{Re}s<k_\infty\})
\\
\text{which have locally uniform polynomial growth on vertical lines}\,.
\end{aligned}
\end{equation}
For our next statement we consider a function $H\in\O_{\text{pg}}(\{-k<\operatorname{Re}s<1\})$, $k\geq 0$. Because of \lemref{mellin_properties} and \corref{melcor1}, for any $f\in \Sch_\delta(\R)$  and any $s_0$ in the interval $(0,k+1)$, the function $s\mapsto M_\delta f\,(s)\,H(1-s)$ is smooth and decays rapidly on the vertical line $\operatorname{Re}s = s_0$. It is therefore integrable over that line.

\begin{lem}\label{mellem3} For $H\in\O_{\text{\rm{pg}}}(\{-k<\operatorname{Re}s<1\})$, the linear function
\[
\Sch_\delta(\R)\,\ni\, f\ \ \mapsto \ \ \frac 1{4\pi i}\int_{\operatorname{Re}s = s_0}M_\delta f\,(s)\, H(1-s)\,ds
\]
is continuous with respect to the topology of $\,\Sch(\R)$. It does not depend on the particular choice of $s_0$, $0< s_0< k+1$.
\end{lem}

\begin{proof} The independence of $s_0$ follows from the Cauchy integral theorem and a limiting argument. Lemma \ref{mellem1} and (\ref{melbound2}) bound $(s_0^2+(\operatorname{Im}s)^2)^{m/2}M_\delta f\,(s)$, for $\operatorname{Re}s=s_0$, in terms of a continuous seminorm, applied to $f^{(m)}$. The lemma follows, since $f\mapsto f'$ is continuous with respect to the topology of $\,\Sch(\R)$.
\end{proof}

The integration pairing exhibits $\,\Sch'_\delta(\R)$ as the continuous dual of $\,\Sch_\delta(\R)$. Hence \lemref{mellem3} implicitly defines a linear map
\begin{equation}
\label{moderategrowth6}
\begin{aligned}
&\Phi_\delta\, :\, \O_{\text{pg}}(\{-k<\operatorname{Re}s<1\})\ \ \longrightarrow\ \ \Sch'_\delta(\R)\,,\ \ \ \text{such that}
\\
&\qquad\qquad\int_\R f(x)\,\Phi_\delta H\,(x)\,dx\, =\,  \frac 1{4\pi i}\int_{\operatorname{Re}s = s_0}  M_\delta f\,(s)\, H(1-s)\,ds
\end{aligned}
\end{equation}
for all $f\in\Sch_\delta(\R)$; the particular choice of $s_0$, $0< s_0 < k+1$, does not matter. The identity remains correct for $f\in\Sch(\R)$ since both sides vanish when $f\in\Sch_{\delta+1}(\R)$.

\begin{lem}\label{mellem4} If $k\geq 0$, the identity
\[
\int_\R f(x)\,\sigma(x)\,dx \ = \ \frac 1{4\pi i}\int_{\operatorname{Re}s = s_0} M_\delta f\,(s)\, M_\delta\sigma\,(1-s)\,ds\,,
\]
holds for all $\,\sigma\in\Sch'(\R)_{k,\,1}$, $f\in \Sch_\delta(\R)$, and $\,0 < s_0 < k+1$.
\end{lem}

\begin{proof} We already know that $s_0$ may be chosen anywhere in the interval $(0,k+1)$. Both sides of the identity vanish if $\,\sigma$ has the opposite parity to $\delta$. We may therefore suppose $\,\sigma\in\Sch'_\delta(\R)$. Recall the definition of $\,\sigma_{x\geq 0}$, $\,\sigma_{x\leq 0}$ in \lemref{truncation}. Since $\,\sigma_{x\leq 0}(-x)=(-1)^\delta\sigma_{x\geq 0}(x)$ and $\,f\in\Sch_\delta(\R)$,
\begin{equation}
\label{moderategrowth9}
\int_\R  f(x)\,\sigma(x)\,dx \ \ = \ \ 2 \int_\R f(x)\,\sigma_{x\geq 0}(x)\,dx\,.
\end{equation}
We now impose the temporary hypothesis
\begin{equation}
\label{moderategrowth10}
\operatorname{supp}\,\sigma\ \ \text{is compact and does not contain the origin}\,.
\end{equation}
Then $\sigma_{x\geq 0}(x)$ has compact support in $\R_{>0}$, which justifies the change of variables $x\rightsquigarrow e^x$ in the integral on the right in (\ref{moderategrowth9}):
\begin{equation}
\label{moderategrowth11}
\int_\R  f(x)\,\sigma(x)\,dx\ \ =  \ \ 2 \int_\R (e^{x/2} f(e^x))(e^{x/2} \sigma_{x\geq 0}(e^x))\,dx\,.
\end{equation}
The function $x\mapsto e^{x/2}f(e^x)$ decays rapidly, along with all its derivatives, both as $x\to -\infty$ and $x\to +\infty$. Thus $e^{x/2}f(e^x)$ is a Schwartz function, with Fourier transform
\begin{equation}
\label{moderategrowth12}
\begin{aligned}
\mathcal F \bigl( e^{x/2}f(e^x) \bigr)(y)\, =  \int_\R e^{x/2}\,f(e^x)\,e(-xy)\,dx \, &=  \int_{\R_{>0}} \!\!\!\! f(x)\,x^{-1/2-2\pi i y}\,dx
\\
= \, \frac 12\int_\R (\sg x)^\delta f(x)\,|x|^{-1/2-2\pi i y}\,dx\, &= \, \frac 12\,M_\delta f\,(1/2-2\pi i y)\,.
\end{aligned}
\end{equation}
Like any distribution with compact support, $x \mapsto e^{x/2}\sigma_{x\geq 0}(e^x)$ has a smooth Fourier transform, which can be computed using the change of variables $x\rightsquigarrow \log x$, along the same lines as (\ref{moderategrowth12}):
\begin{equation}
\label{moderategrowth13}
\begin{aligned}
&\mathcal F \bigl( e^{x/2}\sigma_{x\geq 0}(e^x) \bigr)(y) \ =\ \int_{\R_{>0}} \sigma_{x\geq 0}(x)\,x^{-1/2-2\pi i y}\,dx
\\
&\qquad\ \ = \ \frac 12\int_\R (\sg x)^\delta \sigma(x)\,|x|^{-1/2-2\pi i y}\,dx\, = \ \frac 12\,M_\delta \sigma\,(1/2-2\pi i y)\,.
\end{aligned}
\end{equation}
This is a tempered distribution. In view of (\rangeref{moderategrowth11}{moderategrowth13}), we find
\begin{equation}
\label{moderategrowth14}
\begin{aligned}
\int_\R  f(x)\,\sigma(x)\,dx \ &= \ 2\!\int_\R  \bigl( e^{x/2}f(e^x)\bigr)\bigl( e^{x/2}\sigma(e^x)\bigr)\, dx
\\
&= \ 2\!\int_\R \! \mathcal F\bigl( e^{x/2}f(e^x)\bigr)(-y)\ \mathcal F\bigl( e^{x/2}\sigma(e^x)\bigr)(y)\, dy
\\
&= \ \frac 12\int_\R  M_\delta f\,(1/2+2\pi i y)\,M_\delta \sigma\,(1/2-2\pi i y)\, dy
\\
&= \ \frac 1{4\pi i}\int_{\operatorname{Re}s = 1/2} M_\delta f\,(s)\, M_\delta\sigma\,(1-s)\,ds\,,
\end{aligned}
\end{equation}
still under the simplifying hypothesis (\ref{moderategrowth10}).

To deal with the general case, we choose a cutoff function $\,\phi\in C_c^\infty(\R)$ such that $\,\phi(x)\equiv 1$ near $x=0$ and $\,\phi(-x)\equiv \phi(x)$. Then, for $t>0$,
\begin{equation}
\label{moderategrowth15}
\psi_t(x)\ =_{\text{def}} \ \bigl(1-\phi(x/t)\bigr)\phi(tx)\ \in \ C_c^\infty(\R)
\end{equation}
is an even function, which vanishes near $x=0$. In particular, $\,\psi_t\sigma \in \Sch'_\delta(\R)$ satisfies (\ref{moderategrowth10}). According to \lemref{approximation} and \corref{approximationcor}, $\,\psi_t\sigma \to \sigma$ in the strong distribution topology, as $t\to 0$. This is a statement not just about convergence in $C^{-\infty}(\R)$. Recall that $\,\sigma$ extends across $\infty$ and vanishes there to order $1$;  convergence takes place in $C^{-\infty}(\R\cup\{\infty\})$ when $\,\sigma$ is replaced by this extension. Any Schwartz function extends naturally to a $C^\infty$ function on $\R\cup\{\infty\}$, and the integral of a Schwartz function against $\,\sigma$ can be re-interpreted as an integral over $\R\cup\{\infty\}$ -- this, too, follows from \corref{approximationcor}. We conclude:
\begin{equation}
\label{moderategrowth16}
\int_\R  f(x)\,\psi_t(x)\,\sigma(x)\,dx \ \ \longrightarrow \ \ \int_\R  f(x)\,\sigma(x)\,dx\,, \ \ \ \ \text{as}\ \ t \to 0\,.
\end{equation}
To complete the proof, it suffices to show: there exists $m >0$ such that
\begin{equation}
\label{moderategrowth17}
\begin{aligned}
\bigl((\operatorname{Im}s)^2 + 1\bigr)^{\!-m}\,\bigl| M_\delta(\psi_t\sigma)\,(s)\, - \, M_\delta\sigma\,(s)\bigr|\ \to \ 0\ \ \ \text{as}\ \ t \to 0\,,
\\
\text{locally uniformly in}\ \ -k < \operatorname{Re}s < 1\,.
\end{aligned}
\end{equation}
Since $f$ is a Schwartz function, $\bigl((\operatorname{Im}s)^2 + 1\bigr)^m M_\delta f\,(s)$ decays rapidly along vertical strips. Thus (\rangeref{moderategrowth16}{moderategrowth17}) and the identity (\ref{moderategrowth14}), with $\,\psi_t\sigma$ in place of $\,\sigma$, imply the identity asserted by the lemma.

The verification of (\ref{moderategrowth17}) splits into two local problems, one at $0$, the other at $\infty$. The coordinate change $x\rightsquigarrow 1/x$ relates the two, so we only need to treat the former. In other words, the problem can be solved by showing
\begin{equation}
\label{moderategrowth18}
\begin{aligned}
\bigl((\operatorname{Im}s)^2 + 1\bigr)^{\!-m}\,\bigl|
M_\delta\bigl((1-\phi(x/t))\,\sigma\bigr)(s)\, - \, M_\delta\sigma\,(s)\bigr|\ \to \ 0
\\
\text{as}\ \ t \to 0\,,\ \text{locally uniformly in}\ \ \operatorname{Re}s\,,
\end{aligned}
\end{equation}
for some $m>0$, and for $\operatorname{Re}s>-k$, provided $\,\sigma\in C^{-\infty}(\R)$ has compact support and vanishes to order $k\geq 0$ at the origin. Pointwise convergence follows directly from \lemref{approximation}, applied to the distribution $(\sg x)^\delta|x|^{s-1}\sigma$. To establish locally uniform convergence, we express $\,\sigma$ as in (\ref{fourier7}). Taking one term at a time, we find
\begin{equation}
\label{moderategrowth19}
\begin{aligned}
&\int_\R (\sg x)^\delta\,|x|^{s-1}\,\phi(x/t)\,\textstyle\frac{d^j\ }{dx^j}\bigl(x^{k+j}h_j(x)\bigr)dx\ \ =
\\
&\ \ \ \ = \ (-1)^j\int_\R x^{k+j}h_j(x)\,\textstyle\frac{d^j\ }{dx^j}\bigl((\sg x)^{\delta}\,|x|^{s-1}\,\phi(x/t)\bigr)\,dx\\
&\ \ \ \ = \ {\sum}_{0\leq \ell\leq j}\ c_{j,\,\ell}(s) \int_\R (x/t)^\ell h_j(x)\textstyle(\sg x)^{\delta+k}\,|x|^{s+k-1}\,\phi^{(\ell)}(x/t)\,dx\,,
\end{aligned}
\end{equation}
with suitably chosen constants $c_{j,\,\ell}(s)$ which depend polynomially on $s$. The support of the integrands shrinks down to $\{0\}$ linearly in $t$, and on the support, the integrands are bounded by a multiple of $|x|^{\operatorname{Re}s+k-1}$. We conclude that the integrals tend to $\,0\,$ for $\,\operatorname{Re}s>-k$, locally uniformly in $\,\operatorname{Re}s$. The factor $\bigl((\operatorname{Im}s)^2 + 1\bigr)^{\!-m}$ compensates for the $c_{j,\,\ell}(s)$, so (\ref{moderategrowth18}) follows.
\end{proof}

As one consequence of lemma \ref{mellem4}, the linear map $\Phi_\delta$ defined in (\ref{moderategrowth6}) constitutes a left inverse of the signed Mellin transform:
\begin{equation}
\label{moderategrowth20}
\Phi_\delta\bigl( M_\delta\sigma\bigr)\  = \  \sigma \ \ \ \ \text{if}\ \ \sigma\in\Sch'(\R)_{k,\,1}\,,\ \ k\geq 1\,.
\end{equation}

\begin{lem}\label{mellem5}
The linear map $\Phi_\delta : \O_{pg}(\{0< \operatorname{Re}s<1\})\to\Sch'_\delta(\R)$ is injective.
\end{lem}

\begin{proof}
If $H \in \O_{pg}(\{0< \operatorname{Re}s<1\})$ and $\Phi_\delta H=0$,
\begin{equation}
\label{moderategrowth21}
\int_{\operatorname{Re}s=1/2} M_\delta f\,(s)\,H(1-s)\,ds\  = \  0 \ \ \ \ \text{for all}\ \ f\in\Sch_\delta(\R)\,.
\end{equation}
This will be the case in particular if $f(x)=(\sg x)^\delta |x|^{-1/2}\phi(\log|x|)$ for some $\phi\in C_c(\R)$. In this situation, by (\ref{moderategrowth12}),
\begin{equation}
\label{moderategrowth22}
M_\delta f\,(1/2+2\pi i y)\  = \ 2\, (\mathcal F \phi)(-y)\,.
\end{equation}
We may regard $y\mapsto H(1/2-2\pi i y)$ as a tempered distribution, since it is function of moderate growth. According to (\rangeref{moderategrowth21}{moderategrowth22}), the Fourier transform of this distribution annihilates $\,\phi$, which is an arbitrary smooth function of compact support. But then $H(1/2-2\pi i y)\equiv 0$, which forces $H=0$.
\end{proof}

\begin{lem}\label{mellem6}
If $k_0\geq 0\,$, $k_\infty\geq 1\,$, and $\epsilon>0\,$,\newline
\noindent a)\, $H \in \O_{pg}(\{-k_0-\epsilon< \operatorname{Re}s<k_\infty\})$ implies that $\Phi_\delta H$ vanishes to order $k_0$ at $0$, and\newline
\noindent b)\, $H \in \O_{pg}(\{-k_0< \operatorname{Re}s<k_\infty+\epsilon\})$ implies that $\Phi_\delta H$ has an extension across $\infty$ vanishing there to order $k_\infty$.
\end{lem}

\begin{proof}
We begin with the proof of a). Because of \lemref{mellin_of_schwartz}, $M_\delta(xf')\,(s)=-sM_\delta f\,(s)$ for all $f\in \Sch_\delta(\R)$. Via the defining relation (\ref{moderategrowth6}), this translates into the equation $\frac{d\ }{dx}(x\,\Phi_\delta H)=\Phi_\delta((1-s)H(s))$, which is equivalent to
\begin{equation}
\label{moderategrowth24}
x\,\textstyle\frac{d\ }{dx}\Phi_\delta H\ = \ -\Phi_\delta(s H(s))\,.
\end{equation}
Similarly the identity $M_\delta f\,(s+1)=M_{\delta+1}(xf(x))\,(s)$ for $f\in \Sch_\delta(\R)$ translates into
\begin{equation}
\label{moderategrowth24a}
x\,(\Phi_{\delta} H)(x) \ = \ \Phi_{\delta+1}(s\mapsto H(s+1))\,(x)\,.
\end{equation}
The image under $\Phi_\delta$ of the constant function $1$ is the distribution whose integral against a test function $f\in\Sch_\delta(\R)$ equals $(4\pi i)^{-1}\int_{\operatorname{Re}s=1/2}M_\delta f\,(s)ds$. In view of (\ref{moderategrowth12}), this is the integral of the inverse Fourier transform of the constant function $1$ against $x\mapsto e^{x/2}f(e^x)$; in other words,
\begin{equation}
\label{moderategrowth25}
\Phi_\delta 1\,(x)\ = \ \frac 12\,\bigl(\delta_1(x) + (-1)^\delta\delta_{-1}(x)\bigr)\,,
\end{equation}
where $\delta_n(x)$ denotes Dirac's delta function based at $x=n$. Taken together, (\ref{moderategrowth24}) and (\ref{moderategrowth25}) show that $\Phi_\delta$ maps the space of polynomials $\C[s]$ to the space of linear combinations of delta functions and their derivatives at $x=1$ and $x=-1$. All of these vanish to infinite order at $0$ and extend canonically across $\infty$, so all polynomial functions $H(s)$ satisfy the lemma.

The hypothesis of locally uniform polynomial growth allows us to choose $N\geq 0$ and $C>0$ so that
\begin{equation}
\label{moderategrowth26}
|\operatorname{Re}s + k_0|\leq  \epsilon/2\,,\ \ |\operatorname{Im}s|\geq 1\ \ \Longrightarrow\ \ |H(s)|\ \leq \ C\,|\operatorname{Im}s|^{N-2}\,.
\end{equation}
This inequality remains valid, with a possibly different $C$, if we subtract any polynomial of degree $N-2$. Since the lemma holds for polynomials, we are free to assume that $H(s)$ has a zero of order $N-1$ at $s=-k_0$. Then
\begin{equation}
\label{moderategrowth27}
\widetilde H(s)\ =_{\text{def}}\ s^{-N}\,H(s-k_0)\ \ \text{has at most a first order pole at $s=0$}\,,
\end{equation}
is otherwise regular for $-\epsilon< \operatorname{Re}s<k_0+k_\infty$, and has locally uniform polynomial growth on vertical lines. In particular, we can apply $\Phi_\delta$ to $\widetilde H$. We use the defining relation (\ref{moderategrowth6}), with $s_0=1-\epsilon/2$, and then calculate as in (\rangeref{moderategrowth12}{moderategrowth14}): for $f\in\Sch_\delta(\R)$,
\begin{equation}
\label{moderategrowth28}
\begin{aligned}
\int_\R f(x)\,\Phi_\delta &\widetilde H\,(x)\,dx\ = \ \frac{1}{4\pi i}\int_{\operatorname{Re}s=1-\epsilon/2} M_\delta f(s)\, \widetilde H(1-s)\,ds
\\
&=\ \frac{1}{2}\int_\R M_\delta f(1-\epsilon/2+2\pi i y)\, \widetilde H(\epsilon/2-2\pi i y)\,dy
\\
&= \ \int_\R \mathcal F \bigl( e^{(1-\epsilon/2)x} f(e^x) \bigr)(-y) \, \widetilde H (\epsilon/2 - 2\pi i y)\,dy
\\
&= \ \int_\R e^{(1-\epsilon/2)x} f(e^x)\,\, \mathcal F \bigl( \widetilde H (\epsilon/2 + 2\pi i y)\bigr)(x)\,dx\,.
\end{aligned}
\end{equation}
The integral on the right represents the integration pairing between the Schwartz function $\,e^{(1-\epsilon/2)x}f(e^x)$ against the Fourier transform of the tempered distribution $\,\widetilde H(\epsilon/2 + 2\pi i y)$. Because of (\ref{moderategrowth26}), this tempered distribution is actually a function in $L^1(\R)\cap L^2(\R)$, so both the integral on the right of (\ref{moderategrowth28}) and the Fourier transform itself can be calculated as ordinary, absolutely convergent integrals:
\begin{equation}
\label{moderategrowth29}
\begin{aligned}
&\int_\R e^{(1-\epsilon/2)x} f(e^x)\,\, \mathcal F \bigl( \widetilde H (\epsilon/2 + 2\pi i y) \bigr)(x)\,dx\ =
\\
&\qquad\qquad=\ \int_{\R_{>0}}  f(x)\,\,x^{-\epsilon/2}\, \mathcal F \bigl( \widetilde H (\epsilon/2 + 2\pi i y) \bigr)(\log x)\,dx
\\
&\qquad\qquad=\ \int_{\R_{>0}} \int_\R f(x)\, \widetilde H (\epsilon/2 + 2\pi i y)\,x^{-\epsilon/2-2\pi i y}\,dy\,dx
\\
&\qquad\qquad=\ \int_\R f(x)\left( \frac{(\sg x)^\delta}{4\pi i}\int_{\operatorname{Re}s=\epsilon/2} \widetilde H (s)\,|x|^{-s}\,ds\right)dx\,.
\end{aligned}
\end{equation}
Putting together (\ref{moderategrowth28}) and (\ref{moderategrowth29}) and appealing once more to (\ref{moderategrowth26}), we find that
\begin{equation}
\label{moderategrowth30}
\Phi_\delta \widetilde H\,(x)\  = \ \frac{(\sg x)^\delta}{4\pi i}\int_{\operatorname{Re}s=\epsilon/2} \widetilde H(s)\,|x|^{-s}\,ds
\end{equation}
is continuous for $x\neq 0$. When we shift the line of integration across the origin, we pick up a residue from the first order pole of $\widetilde H(s)$ at $s=0$, and otherwise get the same integral, now over the line $\operatorname{Re}s=-\epsilon/2$. Hence $\Phi_\delta \widetilde H\,(x)$ is bounded near the origin. In view of (\rangeref{moderategrowth24}{moderategrowth24a}), the relationship (\ref{moderategrowth27}) between $H$ and $\widetilde H$ implies
\begin{equation}
\label{moderategrowth31}
\Phi_{\delta+N} H\,(x)\ = \ x^{k_0}\,(-x\,\textstyle\frac{d\ }{dx})^N\, \Phi_\delta\widetilde H\,(x)\,,\ \ \ \ \text{with}\ \ \Phi_\delta\widetilde H\in L^\infty(\R)\,.
\end{equation}
But $\delta\in\Z/2\Z$ can take either value, so $\Phi_\delta H$ vanishes to order $k_0$ at the origin, as asserted in statement a).

The restriction of $\Phi_\delta H$ to $\R-\{0\}$ is completely determined by the relation (\ref{moderategrowth6}) corresponding only to test functions $f\in C_c(\R)$ which vanish near the origin. If $f$ has this property, then so does $x^{-1}f(1/x)$, and both $M_\delta f(s)$ and $M_\delta (x^{-1}f(1/x))\,(s) = M_{\delta+1} (f(x))\,(1-s)$ are entire. We can then shift the line of integration in (\ref{moderategrowth6}) across the origin if necessary, and conclude
\begin{equation}
\label{moderategrowth23}
\f 1x\,\bigl(\Phi_\delta H\bigr)(1/x)  =  \bigl(\Phi_{\delta + 1} H^-\bigr)(x)\,\ \text{on}\,\ \R-\{0\}\,,\ \text{with}\,\ H^-(s) = H(1-s)\,.
\end{equation}
Since the change of variables $x\rightsquigarrow 1/x$ interchanges $0$ and $\infty$, and since the passage from $\,H\,$ to $\,H^-$ has the effect of replacing the hypotheses of b) with those of a), the preceding argument now also implies b).
\end{proof}

At this point we have assembled all the pieces for the proof of theorem \ref{melthm}. Lemma \ref{moderategrowth} tells us that $M_\delta$ induces a linear map $M_\delta : \Sch'_\delta(\R)_{\infty,\infty} \to \mathcal O_{\text{pg}}(\C)$, and lemmas \ref{mellem3} and \ref{mellem4} produce a left inverse $\Phi_\delta : \mathcal O_{\text{pg}}(\C) \to \Sch'_\delta(\R)$. The left inverse takes values in $\Sch'_\delta(\R)_{\infty,\infty}$ by \lemref{mellem6}, and is injective by \lemref{mellem5}, hence defines a two-sided inverse.

We begin the proof of \thmref{4nmelthm} with an observation about the signed Mellin kernel $\,(\sg x)^\delta |x|^{s-1}$, which visibly defines a tempered distribution if $\,\operatorname{Re}s>0$. Integration by parts can be used to extend the definition to all $s\in\C$, $s\notin (2\Z+\delta)\cap\Z_{\leq 0}$, and the resulting tempered distribution depends meromorphically on $\,s$.

\begin{lem}\label{fourier_of_mellin_kernel}
$\ \mathcal F\bigl( x\mapsto (\sg x)^\delta |x|^{s-1}\bigr)\, = \, (-1)^\delta  G_\delta(s)(\sg x)^\delta|x|^{-s}\,.$
\end{lem}

\begin{proof} We choose a cutoff function $\,\phi\in C_c^\infty(\R)$, such that $\phi(x)\equiv 1$ near $x=0$. Then
\begin{equation}
\label{fourier_of_mellin_kernel1}
(\sg x)^\delta |x|^{s-1}\ \ = \ \ \phi(x)\,(\sg x)^\delta |x|^{s-1}\ + \ (1-\phi(x))\,(\sg x)^\delta |x|^{s-1}\,,
\end{equation}
and both summands on the right have meromorphic continuations. We compute the Fourier transform separately for each summand. The first summand is an $L^1$ function when $\operatorname{Re}s>0$, and lies in $L^1(\R)\cap L^2(\R)$ when $\operatorname{Re}s>1/2$. In the latter case, at least, the Fourier transform can be calculated as an ordinary integral. The integral converges, of course, for $\operatorname{Re}s>0$; by analytic continuation,
\begin{equation}
\label{fourier_of_mellin_kernel2}
\mathcal F\bigl( x\mapsto \phi(x)\,(\sg x)^\delta |x|^{s-1}\bigr)(y)\  = \  \int_\R\phi(x)\,(\sg x)^\delta |x|^{s-1}\,e(-xy)\,dx
\end{equation}
for $\operatorname{Re}s>0$, as an integral in the ordinary sense. The second summand lies in $L^1(\R)\cap L^2(\R)$ when $\operatorname{Re}s<0$, in which case the Fourier transform is given by an ordinary integral. That integral exists as a conditionally convergent integral even when $\operatorname{Re}s<1$, provided $y\neq 0$. Arguing by analytic conti\-nuation, one sees that the integral (\ref{fourier_of_mellin_kernel2}), with $(1-\phi)$ in place of $\phi$, represents the Fourier transform, restricted to $\,\R-\{0\}$, in the wider range $\,\operatorname{Re}s<1$. We conclude:
\begin{equation}
\label{fourier_of_mellin_kernel3}
\begin{aligned}
&\mathcal F\bigl( x\mapsto (\sg x)^\delta |x|^{s-1}\bigr)(y)\bigl|_{\{y\neq 0\}}\ = \ \int_\R (\sg x)^\delta |x|^{s-1}\,e(-xy)\,dx\
\\
&\ \ =\ (-\sg y)^\delta |y|^{-s}\!\! \int_\R (\sg x)^\delta |x|^{s-1} e(x)\,dx \, = \, (-\sg y)^\delta  G_\delta(s)\,  |y|^{-s}
\end{aligned}
\end{equation}
in the range $0<\operatorname{Re}s<1$. In other words, the two sides of the identity asserted by the lemma differ by a distribution supported at the origin. But the region $\{0<\operatorname{Re}s<1\}$ is invariant under $\,s\mapsto 1-s$,\, so $\,(\sg x)^\delta |x|^{s-1}$ and $\,(\sg x)^\delta |x|^{-s}$ play essentially symmetric roles. Taking the Fourier transform, we find that the two sides of the identity also differ by the Fourier transform of a distri\-bution supported at the origin, i.e., by a polynomial. That is a contradiction unless (\ref{fourier_of_mellin_kernel3}) remains correct even around the origin. The lemma follows by meromorphic continuation.
\end{proof}

The Fourier transform preserves the parity of Schwartz functions and of tempered distributions. Thus, in proving part a) of \thmref{4nmelthm}, we may as well suppose that $\,\sigma\in\Sch'_\delta(\R)$, in which case $\,\widehat\sigma$ also has parity $\delta$. The identity we need to prove is equivalent to the corresponding identity with $x\,\sigma(x)$ in place of $\sigma$ and $\delta+1$ in place of $\,\delta\,$ -- this follows from \propref{mellin_properties} and the identity $sG_\d(s)=-2\pi iG_{\d+1}(s+1)$, which is equivalent to the Gamma identity $s\G(s)=\G(s+1)$. Division by $x$ does not affect vanishing to infinite order at the origin. We may therefore suppose, without loss of generality, that $\sigma$ is the restriction to $\R$ of a distribution on $\R\cup\{\infty\}$ which vanishes to order $k_\infty\geq 1$ at infinity. In that case $M_\delta\sigma(s)$ is holomorphic on $\{\,\operatorname{Re}s<1\,\}$, of moderate growth on vertical strips. Stirling's formula implies that $G_\delta(s)$ has moderate growth on vertical strips, and from the definition one can see that $G_\delta(s)$ has no poles to the right of $\operatorname{Re}s=0$. In particular, $G_\delta(s)M_\delta\sigma\,(1-s)$ lies in the space $\mathcal O_{\text{pg}}(\{0<\operatorname{Re}s<\infty\})$, on which $\Phi_\delta$ is injective. Since $\Phi_\delta(M_\delta\sigma)=\sigma$, and in view of \lemref{mellem4}, the assertion of part a) of the theorem comes down to the equality
\begin{equation}
\label{moderategrowth32}
\begin{aligned}
\int_\R f(x)\,\widehat\sigma(x)\,dx\ &= \ \frac{(-1)^\delta}{4\pi i}\int_{\operatorname{Re}s=1/2}\!\! M_\delta f\,(s)\,G_\delta(1-s)\,M_\delta\sigma\,(s)\,ds
\\
&=\ \frac{(-1)^\delta}{4\pi i}\int_{\operatorname{Re}s=1/2}\!\!\!\! M_\delta f\,(1-s)\,G_\delta(s)\,M_\delta\sigma\,(1-s)\,ds\,,
\end{aligned}
\end{equation}
for all $f\in\Sch_\delta(\R)$. We use Parseval's identity $\int_\R f(x)\widehat \sigma(x)dx=\int_\R\widehat f(x)\sigma(x)dx$ and \lemref{mellem4} to write (\ref{moderategrowth32}) in the equivalent form
\begin{equation}
\label{moderategrowth33}
\begin{aligned}
&\int_{\operatorname{Re}s=1/2} M_\delta \widehat f\,(s)\,M_\delta\sigma\,(1-s)\,ds\ \ =
\\
&\qquad\qquad=\ \ (-1)^\delta\int_{\operatorname{Re}s=1/2} M_\delta f\,(1-s)\,G_\delta(s)\,M_\delta\sigma\,(1-s)\,ds\,,
\end{aligned}
\end{equation}
again for all $f\in\Sch_\delta(\R)$. This reduces part a) of the theorem to the identity
\begin{equation}
\label{moderategrowth34}
M_\delta \widehat f\,(s)\ = \ (-1)^\delta  G_\delta(s)\,M_\delta f\,(1-s)\ \ \ \ \text{for all}\ \ \ f\in\Sch(\R)\,,
\end{equation}
which is a direct consequence of \lemref{fourier_of_mellin_kernel}.

If $\sigma\in\Sch_\delta'(\R)$ satisfies the hypotheses of part b) -- i.e., vanishing to order at least zero at the origin and extending canonically across $\infty$ -- the Mellin transform $M_\delta\sigma(s)$ is holomorphic on $\{\,\operatorname{Re}s < 0\,\}$, of moderate growth on vertical strips. For the ``if" statement, we suppose that $G_\delta(s)M_\delta\sigma(1-s)$ has no poles for $\operatorname{Re}s < 1$. As was mentioned earlier, $G_\delta(s)$ has moderate growth on vertical strips, so the product $G_\delta(s)M_\delta\sigma(1-s)$ has that property as well. Thus (\ref{moderategrowth6}) and \lemref{mellem6} guarantee the existence of some $\tilde\sigma \in\Sch_\delta'(\R)$ such that $\tilde\sigma$ vanishes to infinite order at $x=0$ and
\begin{equation}
\label{moderategrowth35}
\int_\R f(x)\,\tilde\sigma(x)\,dx\  = \ \frac{(-1)^\delta}{4\pi i}\int_{\operatorname{Re}s=1/2} M_\delta f\,(s)\,G_\delta(1-s)M_\delta\sigma(s)\,ds \,,
\end{equation}
for all $f\in\Sch(\R)$. The change of variables $s\rightsquigarrow 1-s$, the identity (\ref{moderategrowth34}), and \lemref{mellem4} transform this into the equation $\int_\R f(x)\tilde\sigma(x)dx = \int_\R \widehat f(x)\sigma(x)dx$, so $\widehat\sigma=\tilde\sigma$ vanishes to infinite order at $x=0$, as asserted. To establish the ``only if" statement, we suppose that $\widehat\sigma$ vanishes to infinite order at $x=0$ and apply the first part of the theorem to $\widetilde\sigma$\,: $M_\delta\sigma(s)= G_\delta(s)M_\delta\widehat\sigma(1-s)$ for $\operatorname{Re}s \gg 0$; the factor $(-1)^\delta$ has disappeared since the inverse Fourier transform and the Fourier transform of any $\sigma\in\Sch_\delta'(\R)$ are related by this sign factor. Appealing to the functional equation (\ref{whyGdelta2}) for $G_\delta(s)$, we find
\begin{equation}
\label{moderategrowth36}
M_\delta\widehat\sigma\,(s)\  = \  (-1)^\delta\,G_\delta(s)\,M_\delta\sigma\,(1-s)\,.
\end{equation}
Since $\sigma$ vanishes to order at least one at $x=0$, $\widehat\sigma\in \Sch'(\R)_{\infty,1}$ by \thmref{fourier}, so $M_\delta\widehat\sigma(s)$ is regular for $\operatorname{Re}s < 1$, as remained to be shown.

\section{Examples: Applications to $L$-functions}
\label{examplesec}

The methods developed in the previous two sections can be used to prove the analytic continuation and functional equations of various $L$-functions. We shall show how this works by giving the proofs for the Riemann zeta function, Dirichlet $L$-functions, and $L$-functions for automorphic forms on $GL(2,\R)$. These are exceedingly well-known results, of course -- the aim is to illustrate our technique, not to explain the results. Our paper \cite{MS2} contains more substantial applications.

We begin with the Riemann zeta function $\,\zeta(s)=\sum_{n\geq 1} n^{-s}$. The point of departure is the tempered distribution $\,\delta_\Z(x)=\sum_{n\in\Z}\delta_n(x)$, i.e., the sum of the Dirac delta functions based at all the integers. The Poisson summation formula for $\Z$ can be paraphrased by the identity $\,\delta_\Z = \widehat\delta_\Z$. Since $\,\widehat\delta_n(x)=e(-nx)$, we can write this in the equivalent form
\begin{equation}
\label{poisson1}
{\sum}_{n\neq 0}\ \delta_n(x)\ - \ (1 - \phi(x))\ \ = \ \ {\sum}_{n\neq 0}\ e(nx)\ - \ \delta_0(x)\ + \ \phi(x)\,;
\end{equation}
here $\,\phi\in C_c(\R)$ denotes a cutoff function such that $\,\phi(x)\equiv 1$ near $x=0$. The distribution on the left of the equality sign vanishes identically near $x=0$, and therefore vanishes to infinite order at $0$. The distribution on the right differs from a compactly supported distribution by one that is periodic, without constant term. Thus, according to \propref{oscillation}, the right hand side has a canonical extension across $\infty$. Since the two distributions are equal, the discussion leading up to \propref{mellin_properties} allow us to conclude that
\begin{equation}
\label{zeta1}
\begin{aligned}
&M_0\bigl( \textstyle{\sum}_{n\neq 0}\ e(nx) - \delta_0(x) + \phi(x)\bigr)(s) \ \ =
\\
&\ \ \ \ \ \ =\ \ M_0\bigl( \textstyle{\sum}_{n\neq 0}\ \delta_n(x) - 1 + \phi(x)\bigr)(s)\ \ \ \text{is an entire function.}
\end{aligned}
\end{equation}
These Mellin transforms are globally defined. That is not the case for the summands on the left and right of the identity (\ref{poisson1}): both $\sum_{n\neq 0}\,\delta_n(x)$ and $\,(1-\phi(x))$ have Mellin transforms, in the sense of our definition, only for $\operatorname{Re}s<0$, hence
\begin{equation}
\label{zeta2}
\begin{aligned}
&M_0\bigl( \textstyle{\sum}_{n\neq 0}\ \delta_n(x) - 1 + \phi(x)\bigr)(s)\ \ =
\\
&\qquad\ \ \ = \ \ M_0\bigl( \textstyle{\sum}_{n\neq 0}\ \delta_n(x)\bigr)(s)\ +\ M_0(\phi(x)-1)(s)\ \ \ \text{for $\,\operatorname{Re}s<0$}\,.
\end{aligned}
\end{equation}
Quite similarly,
\begin{equation}
\label{zeta3}
\begin{aligned}
&M_0\bigl( \textstyle{\sum}_{n\neq 0}\ e(nx) - \delta_0(x) + \phi(x)\bigr)(s)\ \ =
\\
&\ \ \ \ = \ \ M_0\bigl( \textstyle{\sum}_{n\neq 0}\ e(nx)\bigr)(s)\ -\ M_0\delta_0\,(s)\ +\ M_0\phi\,(s)\ \ \text{for $\,\operatorname{Re}s>1$}\,.
\end{aligned}
\end{equation}
One can appeal to corollary \ref{approximationcor} to justify the heuristically obvious equation
\begin{equation}
\label{zeta4}
M_0\bigl(\textstyle{\sum}_{n\neq 0}\ \delta_n(x)\bigr)(s)\ =  \ 2\,\zeta(1-s)\qquad (\,\operatorname{Re}s<0\,)\,.
\end{equation}
On the other hand, lemma \ref{mellin_and_dirichlet} implies
\begin{equation}
\label{zeta5}
M_0\bigl(\textstyle{\sum}_{n\neq 0}\, e(nx)\bigr)(s)\ =  \ 4\,(2\pi)^{-s}\,\Gamma(s)\,\cos(\pi s/2)\,\zeta(s)\ \ \ \ (\,\operatorname{Re}s>1\,)\,.
\end{equation}
The function $\,\phi$ satifies the hypotheses of \lemref{mellin_of_one}, hence
\begin{equation}
\label{zeta6}
\begin{aligned}
M_0\phi\,(s)\,\ \text{and}\,\ M_0(\phi-1)\,(s)\,\ \text{extend meromorphically to $\C$,}
\\
\text{the two extensions coincide, and have no singularities}
\\
\text{except for a simple pole at $s=0$, with residue $2$}\,.
\end{aligned}
\end{equation}
The analytic continuation of $\,\zeta(s)-1/(s-1)$ follows from (\rangeref{zeta1}{zeta2}), (\ref{zeta4}), and (\ref{zeta6}). The functional equation can be read off from (\rangeref{zeta1}{zeta6}) and the identity $\,M_0\delta_0\,(s)=0$, $\,\operatorname{Re}s>1$, which is a special case of corollary \ref{mellin_of_dirac}.

The case of Dirichlet $L$-functions is simpler from the analytic point of view, but requires some combinatorics. We recall the definition of a Dirichlet character modulo $q>1$: a multiplicative function $\chi: \Z\to\C$ obtained from a character of the multiplicative group $(\Z/q\Z)^*$, which is extended to $\Z/q\Z$ by setting it equal to zero on non-units, and then lifted from $\Z/q\Z$ to $\Z$. One calls the Dirichlet character $\,\chi$ primitive if the character of $(\Z/q\Z)^*$ which it encodes is not lifted from a quotient $(\Z/q'\Z)^*$, corresponding to a proper divisor $q'|q$. In the primitive case, one calls $q$ the conductor of $\,\chi$. Like $\zeta(s)$, the Dirichlet series
\begin{equation}
\label{dirichlet1}
L(s,\chi)\ =  \ {\sum}_{n=1}^\infty\ \chi(n)\,n^{-s}
\end{equation}
converges absolutely and uniformly for $\,\operatorname{Re}s>1$. This is the $L$-function of the primitive Dirichlet character $\,\chi$.

For the remainder of this dicussion we fix a particular primitive Dirichlet character $\,\chi$ and the conductor $q>1$. Then $\,\bar\chi$ corresponds to the reciprocal character of $(\Z/q\Z)^*$ and is also primitive. The tempered distribution
\begin{equation}
\label{dirichlet2}
\tau_\chi(x)\ =  \ {\sum}_{n\in\Z}\ \chi(n)\,\delta_n(x)\ = \ {\sum}_{a\in (\Z/q\Z)^*}\, \chi(a)\,{\sum}_{n\in\Z,\, n\equiv a\,(q)} \,\delta_n(x)
\end{equation}
vanishes near the origin, since $\,\chi(0)=0$. To calculate  $\widehat\tau_\chi$, we note that $\sum_{m\in\Z}\delta_{a+mq}(x)= q^{-1}\delta_\Z((x-a)/q)$ has $\sum_{n\in\Z}e(-na/q)\,\delta_n(qx)$ as Fourier transform, which implies
\begin{equation}
\label{dirichlet3}
\widehat\tau_\chi(x)\ = \ {\sum}_{a\in (\Z/q\Z)^*}\,{\sum}_{n\in\Z}\ \chi(a)\,e(-na/q)\,\delta_n(qx)\,.
\end{equation}
A basic identity for Dirichlet characters asserts that
\begin{equation}
\label{dirichlet4}
{\sum}_{a\in (\Z/q\Z)^*}\ \chi(a)\,e(na/q)\ =\ g_\chi\,\bar\chi(n)\,,
\end{equation}
with $g_\chi=\sum_{a\in (\Z/q\Z)^*}\chi(a)e(a/q)$ denoting the so-called Gauss sum. Hence
\begin{equation}
\label{dirichlet5}
\widehat\tau_\chi(x)\ = \ \chi(-1)\,g_\chi\,\tau_{\bar\chi}(qx)
\end{equation}
also vanishes near the origin. Appealing to \thmref{fourier} we see that both $\,\tau_\chi$ and $\,\widehat\tau_{\chi}$ have canonical extensions across $\infty$, which lets us conclude that
\begin{equation}
\label{dirichlet6}
M_\delta\tau_{\!\chi}\,(s)\ \ \text{and}\ \ \ M_\delta\widehat\tau_{\!\chi}\,(s)\ \ \text{are entire holomorphic functions}.
\end{equation}
Theorem \ref{4nmelthm} relates the two Mellin transforms:
\begin{equation}
\label{dirichlet7}
M_\delta\tau_{\!\chi}\,(s)\ = \ G_\delta(s)\,M_\delta\widehat\tau_{\!\chi}\,(1-s)\,,
\end{equation}
with $\,G_0(s)=2(2\pi)^{-s}\Gamma(s)\cos(\pi s/2)$ and $\,G_1(s)=2i(2\pi)^{-s}\Gamma(s)\sin(\pi s/2)$. We now fix $\delta\in\Z/2\Z$ so that
\begin{equation}
\label{dirichlet8}
\chi(-1)\ = \ \bar\chi(-1)\ = \ (-1)^\delta\,.
\end{equation}
Then, in analogy to (\ref{zeta4}),
\begin{equation}
\label{dirichlet9}
M_\delta\tau_{\!\chi}\,(s)\ = \ 2\,L(1-s,\chi)\,,\ \ \ M_\delta\tau_{\!\bar\chi}(qx)\,(s)\ = \ 2\,q^{-s}\,L(1-s,\bar\chi)\,,
\end{equation}
both in the range $\operatorname{Re}s<0$. The identities (\rangeref{dirichlet5}{dirichlet9}) give the analytic continuation of the two $L$-functions, as well as the functional equation
\begin{equation}
\label{dirichlet10}
L(s,\chi)\ = \ (-1)^\delta\,G_\delta(1-s)\,g_\chi\,L(s,\bar\chi)\,.
\end{equation}
This can be written in more symmetric form. For details, and for the history of the functional equation, we refer the reader to \cite{daven}.

To keep the discussion of $GL(2)$ reasonably succinct, we define the notion of a $GL(2,\Z)$-automorphic distribution without any further introduction. Motivation and a much more general notion of automorphic distribution can be found in \cites{S,MS2}. We fix parameters $\nu\in\C$, $\delta\in\Z/2\Z$, and define $V_{\nu,\delta}^{-\infty}$ as the vector space of pairs of distributions $(\tau,\widetilde\tau)$, with $\tau,\widetilde\tau\in C^{-\infty}(\R)$ related by the equation
\begin{equation}
\label{autdist1}
\widetilde\tau(x)\ \ = \ \ |x|^{2\nu-1}\,\tau(-1/x)\qquad\qquad (\,x\neq 0\,)\,.
\end{equation}
Then $\,\tau$ determines $\,\widetilde\tau$ except at $x=0$. We therefore may, and shall, think of vectors in $V_{\nu,\delta}^{-\infty}$ as a distribution $\,\tau$, together with the datum of a specific extension of $\,x\mapsto|x|^{2\nu-1}\tau(-1/x)$ across $x=0$. The group $G=GL(2,\R)$ acts on $V_{\nu,\delta}^{-\infty}$ by the rule
\begin{equation}
\label{autdist2}
\bigl(\pi_{\nu,\delta}(g)\tau\bigr)(x)\, = \, \frac{(\sg\det g)^\delta}{|cx+d|^{1-2\nu}}\ \tau\bigl(\,\frac{ax+b}{cx+d}\,\bigr)\ \ \text{if}\ \ g^{-1} = \, \ttwo{a}{b}{c}{d}\in G\,.
\end{equation}
At points where the denominator $cx+d$ vanishes this identity retains meaning when re-written in terms of $\,\widetilde\tau$. In any case, $\,\pi_{\nu,\delta}$ defines a representation of $G$ on $V_{\nu,\delta}^{-\infty}$. By definition, the invariants for $\,\Gamma=GL(2,\Z)\,$,
\begin{equation}
\label{autdist3}
(V_\nu^{-\infty})^\Gamma\ = \ \{\,\tau\in V_{\nu,\delta}^{-\infty}\, \mid \, \pi_{\nu,\delta}(\gamma)\tau=\tau\,\ \text{for all}\,\ \gamma\in\Gamma\}\,,
\end{equation}
constitute the space of $GL(2,\Z)$-automorphic distributions corresponding to $(\nu,\delta)$. For a $\Gamma$-invariant distribution $\,\tau$, the invariance condition $\,\pi_{\nu,\delta}(g)\tau=\tau$,\, with $\,a=d=0$, $\,b=-c=1\,$ in (\ref{autdist2}), implies $\,\widetilde\tau=\tau$, so we no longer need to specify $\,\widetilde\tau$ separately.

To see how automorphic distributions arise from $GL(2,\Z)$-automorphic forms in the usual sense, we first consider a modular form of weight $2k$ and parity $\delta$, i.e., a holomorphic function $F(z)$ on $\,\C-\R$ which grows at most polynomially in $\,y=\operatorname{Im}z\,$ as $\,|y|\to\infty$, and satifies the automorphy condition
\begin{equation}
\label{holom1}
F(z)\, = \, (cz+d)^{-2k}\,F\bigl(\,\frac{az+b}{cz+d}\,\bigr) \ \ \text{for all}\ \ \gamma = \ttwo{a}{b}{c}{d}\in SL(2,\Z)\,,
\end{equation}
as well as the parity condition $\,F(-z)=(-1)^\delta F(z)$. The limit
\begin{equation}
\label{holom2}
\tau(x)\ \ = \ \ {\lim}_{y\to 0^+} \ F(x+iy)
\end{equation}
converges in the strong distribution topology to an automorphic distribution $\,\tau\in V_{\nu,\delta}^{-\infty}$, with $\nu=1/2-k$ \cite{S}. Next we consider a Maass form of parity $\delta$, i.e., an $SL(2,\Z)$-invariant eigenfunction $F(x+iy)$ of the hyperbolic Laplace operator $\,\Delta$ on the upper half plane $\,\mathcal H$, of polynomial growth in $y$, which obeys the parity condition $F(-x+iy)=(-1)^\delta F(x+iy)$. We choose $\nu\in\C$ so that $\,\Delta\, F = (1/4-\nu^2)F$. In this situation, $\,F(x+iy)$ has an asymptotic expansion for $\,y\to 0^+$,
\begin{equation}
\label{maass1}
F(x+iy)\ \sim \ y^{1/2-\nu} \textstyle{\sum}_{n\geq 0}\, \tau_{\nu,n}(x)\,y^n\ + \ y^{1/2+\nu} \textstyle{\sum}_{n\geq 0}\, \tau_{-\nu,n}(x)\,y^n\,,
\end{equation}
and $\,\tau = \tau_{\nu,0}\,$ is the automorphic distribution $\,\tau\in V_{\nu,\delta}^{-\infty}$ corresponding to $\,F$ \cite{S}. Then $\nu$ is determined only up to sign, but the two automorphic distributions in $\,V_{\pm\nu,\delta}^{-\infty}\,$ are related by the so-called standard intertwining ope\-rator $V_{-\nu,\delta}^{-\infty}\to V_{\nu,\delta}^{-\infty}$. The passage from a modular form or Maass form $\,F$ to the automorphic distribution $\,\tau$ can be reversed: there is a simple, explicit formula for $\,F$ in terms of $\,\tau$ \cite{S}.

The automorphy condition (\rangeref{autdist2}{autdist3}), with $a=b=d=1$ and $c=0$, implies $\,\tau(x+1)=\tau(x)$, so $\,\tau\in V_{\nu,\delta}^{-\infty}$ can be developed as a Fourier series:
\begin{equation}
\label{autdist4}
\tau(x)\ \ =\ \  {\sum}_{n\in\Z}\ |n|^{-\nu}\,a_n\,e(nx)\,;
\end{equation}
the factor $|n|^{-\nu}$ has the effect of making the coefficients $a_n$ independent of the choice between $\nu$ and $-\nu$ in the Maass case, except for a normalizing factor. The parity condition implies
\begin{equation}
\label{autdist5}
a_n\ \ =\ \ (-1)^\delta\,a_{-n} \,,
\end{equation}
both in the holomorphic and the Maass case. We shall call $\,\tau$ cuspidal if
\begin{equation}
\label{autdist6}
a_0\ =\ 0\,,\ \ \text{and} \ \ \tau\,\ \text{vanishes to infinite order at $\,x=0$}\,.
\end{equation}
The first of these two conditions ensures that $\,\tau$ has a canonical extension across $\infty$, and the second can be paraphrased by saying that the automorphy condition $\,\tau(x)=|x|^{2\nu-1}\tau(-1/x)$ for $x\neq 0$ extends as an equality of canonical extensions to $\R\cup\{\infty\}$. One can show without great difficulty that our definition agrees with the usual notion of cuspidality for the modular form or Maass form from which $\,\tau$ was derived. Non-cuspidal automorphic forms should be thought of as attached to automorphic representations of $GL(1)$. As such, they are less interesting in the context of $GL(2)$. To simplify the arguments, we shall consider only cuspidal automorphic distributions $\,\tau\in V_\nu^{-\infty}$.

Because of (\rangeref{autdist4}{autdist5}), $\,\tau$ is completely determined by the Fourier coefficients $\,a_n$, $\,n>0$, or equivalently, by the Dirichlet series
\begin{equation}
\label{autdist7}
L(s,\tau)\ \ = \ \ {\sum}_{n\geq 1}\ a_n\, n^{-s}\,,
\end{equation}
which converges for $\,\operatorname{Re}s\gg 0$. This is the standard $L$-function of $\,\tau$, though in the holomorphic case, the definition (\ref{autdist7}) differs from the classical definition by an additive shift in the variable $s$, which makes the functional equation relate values at $s$ and $1-s$, rather than at $s$ and $2k-s$.

We can appeal either to \lemref{mellin_and_dirichlet} or to \thmref{4nmelthm} and corollary \ref{approximationcor}, to conclude
\begin{equation}
\label{autdist8}
M_\delta \tau\,(s)\ \ = \ \ 2\,G_\delta(s)\,L(s+\nu,\tau)\qquad\qquad (\,\operatorname{Re}s \gg 0\,)\,.
\end{equation}
The hypothesis of cuspidality, in conjunction with \lemref{moderategrowth}, implies
\begin{equation}
\label{autdist9}
M_\delta \tau\,(s)\ \ \text{is entire, of moderate growth on vertical strips}\,.
\end{equation}
Since $\,\tau(x)=|x|^{2\nu-1}\tau(-1/x)$, the definition of $M_\delta$ and the change of variables formula (\ref{changeofvariables}) lead to the identity
\begin{equation}
\label{autdist10}
M_\delta \tau\,(s)\ \ =\ \ (-1)^\delta\,M_\delta \tau\,(1-s-2\nu)\,.
\end{equation}
At this point, (\rangeref{autdist8}{autdist10}) provide a meromorphic continuation of $L(s,\tau)$ and the functional equation
\begin{equation}
\label{autdist11}
G_\delta(s-\nu)\,L(s,\tau)\ \ =\ \ (-1)^\delta\,G(1-s-\nu)\,L(1-s,\tau)\,,
\end{equation}
which can be re-written in various equivalent ways. To see that $L(s,\tau)$ is entire, we observe that
\begin{equation}
\label{autdist12}
\widehat\tau(x)\ \ =\ \ {\sum}_{n\in\Z}\ |n|^{-\nu}\,a_n\,\delta_n(x)
\end{equation}
vanishes near the origin and extends canonically across $\infty$ by (\ref{autdist6}) and \thmref{fourier}. Consequently
\begin{equation}
\label{autdist13}
L(\tau,s)\ \ =\ \ 1/2\ M_\delta\widehat\tau\,\,(1+\nu-s)\qquad\qquad(\,\operatorname{Re}s \gg 0\,)
\end{equation}
extends to an entire function, of moderate growth on vertical strips.

We can also use our methods to prove a ``converse theorem", which reconstructs a cuspidal automorphic distribution $\,\tau\in V_{\nu,\delta}^{-\infty}$ from its $L$-function and functional equation. We fix $(\nu,\delta)$ as before, and suppose that $a_n$, $n\geq 1$, is a sequence of complex numbers which grows at most polynomially with $n$. Then
\begin{equation}
\label{autdist14}
L(s)\ \ = \ \ {\sum}_{n\geq 1}\ a_n\, n^{-s}
\end{equation}
converges for $\,\operatorname{Re}s \gg 0$, and
\begin{equation}
\label{autdist15}
\tau(x)\ \ = \ \ {\sum}_{n\neq 0}\ |n|^{-\nu} \, a_n\, e(nx)\,,\qquad \text{with $a_{-n}=(-1)^\delta a_n$}\,,
\end{equation}
converges to a periodic distribution $\,\tau$ without constant term. In particular, $\,\tau$ has a canonical extension across $\infty$. We had just argued that if $\,\tau$ also vanishes to infinite order at $0$, then both $L(s)$ and $G_\delta(s-\nu)L(s)$ extend to entire functions, of moderate growth on vertical strips. We shall now reverse that argument: we make the holomorphic extension and growth behavior of $L(s)$ and $G_\delta(s-\nu)L(s)$ the hypothesis, and shall deduce that $\,\tau$ vanishes to infinite order at $0$.

Indeed, theorems \ref{melthm} and \ref{4nmelthm} guarantee the existence of $\,\sigma\in\Sch'_\delta(\R)$ such that both $\,\sigma$ and its Fourier transform $\,\widehat\sigma$ vanish to infinite order at the origin, $M_\delta\sigma\,(s)=2G_\delta (s) L(s+\nu)$, and $M_\delta\widehat\sigma\,(s)=2 L(1-s+\nu)$. These identities hold globally. The identities $M_\delta\tau\,(s)=2G_\delta (s) L(s+\nu)$ and $M_\delta\widehat\tau\,(s)=2 L(1-s+\nu)$, which can be derived just as in the proof of the functional equation, only hold for $\,\operatorname{Re}s \gg 0\,$ and $\,\operatorname{Re}s \ll 0\,$, respectively. We can push up the order of vanishing at the origin to at least one when we multiply $\,\tau$ by a high enough power $x^n$ of the variable $x$. Doing so has no effect on the existence of a canonical extension across $\infty$. On the other hand, multiplication by $x^n$ shifts both the argument $s$ of the Mellin transform and the parity $\delta$ by $n$. Hence, for $n$ large enough,
\begin{equation}
\label{autdist16}
M_{\delta+n}(x^n\sigma)\,(s)\ =\ M_{\delta+n}(x^n\tau)(s) \ \ \ \ \ \ \text{if}\ \ \ \operatorname{Re}s>-1\,,
\end{equation}
which according to (\ref{moderategrowth20}) implies $x^n\sigma=x^n\tau$, or equivalently $\tau=\sigma+P(\frac{d\ }{dx})\delta_0$, for some polynomial $P(X)\in\C[X]$. Taking the Fourier transform, we find that $\widehat\sigma(x)$ and $\widehat\tau(x)=\sum_{n\neq 0} |n|^\nu a_n \delta_n(x)$ -- both of which have canonical extensions across $\infty$ -- differ by $P(2\pi i x)$. That, in conjunction with corollary \ref{equiv_of_def}, forces $P(X)$ to vanish as polynomial, so $\,\tau=\sigma$ does vanish to infinite order at $x=0$. For future reference, we summarize:

\begin{prop}\label{L_and_inf_order}
Fix $(\nu,\delta)\in\C\times\Z/2\Z$ and suppose that the Dirichlet series $L(s)$ and the distribution $\,\tau\in\Sch'_\delta(\R)$ are related as in {\rm(\rangeref{autdist14}{autdist15})}. Then $\,\tau$ vanishes to infinite order at $x=0$ if and only if both $L(s)$ and $G_\delta(s-\nu)L(s)$ extend to entire functions, of moderate growth on vertical strips.
\end{prop}

Let us continue with the proof of the converse theorem. We not only suppose that $L(s)$ and $G_\delta(s-\nu)L(s)$ extend to entire functions, of moderate growth on vertical strips, but also impose the functional equation
\begin{equation}
\label{autdist17}
G_\delta(s-\nu)\,L(s)\ \ =\ \ (-1)^\delta\,G(1-s-\nu)\,L(1-s)\,.
\end{equation}
As we have seen, $\,\tau$ vanishes to infinite order at $0$ and has a canonical extension across $\infty$. Hence there exists $\,\widetilde\tau\in\Sch'_\delta(\R)$, also vanishing to infinite order at $0$ and having a canonical extension across $\infty$, such that $\,\widetilde\tau(x)=|x|^{2\nu - 1}\tau(-1/x)$ for $x\neq 0$. The functional equation is equivalent to the equality $M_\delta\widetilde\tau=M_\delta\tau$, hence by \thmref{melthm}, to the equality $\,\widetilde\tau=\tau$. The pair $\tau,\,\widetilde\tau=\tau$ defines a vector in $V_{\nu,\delta}^{-\infty}$ which is $\pi_{\nu,\delta}$-invariant under
\begin{equation}
\label{autdist18}
\gamma_n \ = \ \ttwo{1}{n}{0}{1},\ \ \ \gamma(\eta_1,\eta_2)\ = \ \ttwo{\eta_1}{0}{0}{\eta_2},\ \ \ w\ = \ \ttwo{0}{-1}{1}{0},
\end{equation}
with $n\in\Z$ and $\eta_1,\,\eta_2\in\{\pm1\}$; in the case of $\,\gamma_n$, the invariance follows from the periodicity of $\,\tau$ and its canonical extension across $\infty$, for $\gamma(\eta_1,\eta_2)$ from the parity condition on $\,\tau$, and for $w$ from the identity $\,\widetilde\tau=\tau$. The matrices (\ref{autdist18}) generate $\Gamma=GL(2,\Z)$, so $\,\tau$ does define a $\Gamma$-automorphic distribution. It satisfies the cuspidality condition (\ref{autdist6}) by construction. This completes the proof of the converse theorem.

\section{The operators $T_{\alpha,\eta}$}\label{talphasec}

In this section we introduce and study the operators $T_{\alpha,\eta}$, which play an important role in our proof of the Voronoi summation formula for $GL(3)$ \cite{MS2}.

For $\,\alpha\in\C$ and $\,\eta\in\Z/2\Z$, the integral that computes the Fourier transform in the expression
\begin{equation}
\label{talphadef}
T_{\alpha,\eta}(f)\ \ = \ \ \mathcal F\bigl(\, x \mapsto f(1/x)\,(\sg x)^\eta\,|x|^{-\alpha-1} \, \bigr)\qquad (\, f \in \Sch(\R)\,)
\end{equation}
converges absolutely when $\,\operatorname{Re}\alpha>0$. We can make sense of $T_{\alpha,\eta}(f)(y)$, at least as function on $\R-\{0\}$, for all $\,\alpha\in\C$ because the Fourier kernel $x\mapsto e(-xy)$, $\,y\neq 0$, vanishes to infinite order at $\infty$. With this extended definition the values $T_{\alpha,\eta}(f)(y)$, $\,y\neq 0$, depend holomorphically on $\,\alpha$.

\begin{lem}\label{talphalem}
The function $T_{\alpha,\eta}(f)(x)$, $\,x\neq 0$, is infinitely differentiable. It decays rapidly as $|x|\to\infty$, along with all its derivatives.
\end{lem}

\begin{proof}
For $\,\operatorname{Re}\alpha> 0$, the Fourier integral converges absolutely, so $T_{\alpha,\eta}(f)$ is bounded. Multiplying $T_{\alpha,\eta}(f)$ by $2\pi ix$ has the same effect as differentiating the argument of $\,\mathcal F$, which results in an expression of the same type, but with $\,\alpha$ raised by $1$. Repeating this reasoning gives the rapid decay of $\,T_{\alpha,\eta}(f)(x)$, for any $\,\alpha$. Differentiating the function $\,T_{\alpha,\eta}(f)$ also results in an expression of the same type, now with $\,\alpha$ lowered by $1$. Thus $\,T_{\alpha,\eta}(f)$ is differentiable, and the derivative decays rapidly, etc.
\end{proof}

We can apply the operator $T_{\alpha,\eta}$ also to distributions which extend cano\-nically across $\infty$. If $\,\sigma\in\Sch'(\R)$ has this property, $(\sg x)^\eta |x|^{-\alpha-1}\sigma(1/x)$ is well defined as distribution on $\,\R$, with holomorphic dependence on $\,\alpha$. Since $\,\sigma(1/x)$ extends across $\infty$, $x\mapsto \sigma(1/x)$ has the temperedness property at $\infty$, and that remains the case when we multiply this distribution with $(\sg x)^\eta |x|^{-\alpha-1}$. In short, $(\sg x)^\eta |x|^{-\alpha-1}\sigma(1/x)$ is a tempered distribution, to which we can apply the Fourier transform:
\begin{equation}
\label{talphadef'}
\begin{aligned}
T_{\alpha,\eta}(\sigma)\ \ = \ \ \mathcal F\bigl(\, x \mapsto \sigma(1/x)\,(\sg x)^\eta\,|x|^{-\alpha-1} \, \bigr)\ \in \ \Sch'(\R)
\\
\text{if $\,\sigma\in\Sch'(\R)$ has a canonical extension across $\infty$}.
\end{aligned}
\end{equation}
The two definitions (\rangeref{talphadef}{talphadef'}) are consistent: if we regard a Schwartz function $f$ as distribution with canonical extension across $\infty$, the definitions agree if $\,\operatorname{Re}\alpha> 0$ because the Fourier transform of a function in $L^2(\R)\cap L^1(\R)$ has unambiguous meaning. For other values of $\,\alpha$ we can argue by analytic continuation.

As things stand, we cannot compose two operators of the type (\ref{talphadef}). To remedy this deficiency, we shall first extend the domain of definition of $T_{\alpha,\eta}$, and then show that $T_{\alpha,\eta}$ maps this extended domain to itself.

\begin{definition}\label{sisdef}
A function $f\in C^\infty(\R-\{0\})$ has a singularity of type $(\alpha,\eta)\in \C\times\Z/2\Z$ at $x=0$ if there exist $C^\infty$ functions $f_0,\,f_1,\,\dots,\,f_n$, defined near $x=0$, such that
\[
f(x)\ \ = \ \ {\sum}_{0\leq j\leq n}\ (\sg x)^\eta\,|x|^\alpha\,(\log|x|)^j\,f_j(x)\ \ \ \ \text{for}\ \ \ 0 < |x| \ll 1\,.
\]
When it is chosen minimally, the integer $n$ will be called the index of the singularity. We say that $f$ has a simple singularity at $x=0$ if, locally near $x=0$, it can be expressed as a sum of functions $g_j$, $\,1\leq j \leq m$, each of which has a singularity of some type $(\alpha_j,\eta_j)\in \C\times\Z/2\Z$. We let $\,\Sch_{\text{\rm{sis}}}(\R)$ denote the space of functions $\,f\in C^\infty(\R-\{0\})$ which have a simple singularity at $x=0$ and decay rapidly, along with all of their derivatives, as $|x|\to\infty$.
\end{definition}

Multiplication with a $C^\infty$ function does not change the type of a singularity, and differentiation changes the type from $(\alpha,\eta)$ to $(\alpha-1,\eta+1)$. In particular, $\,\Sch_{\text{\rm{sis}}}(\R)$ is a module over the ring of linear differential operators with coefficients which are $C^\infty$, and which grow at most polynomially as $|x|\to\infty$, along with all their derivatives. Pointwise multiplication turns $\,\Sch_{\text{\rm{sis}}}(\R)$ into a ring.

The definition (\ref{talphadef}) of $T_{\alpha,\eta}(f)$ has meaning even for a function $f$ with a singularity of type $(\alpha,\eta)$, because \propref{xalpha} also allows for powers of $\log|x|$. We can therefore extend the definition to $\,\Sch_{\text{\rm{sis}}}(\R)$,
\begin{equation}
\label{talpha1}
T_{\alpha,\eta}\ :\ \Sch_{\text{\rm{sis}}}(\R)\ \ \longrightarrow\ \ C^\infty(\R-\{0\})\,.
\end{equation}
Lemma \ref{talphalem} remains valid in the current setting, though the proof needs to be adapted slightly.

\begin{thm}\label{talphathm1}
The operator $\,T_{\alpha,\eta}$ maps the space $\,\Sch_{\text{\rm{sis}}}(\R)$ to itself.
\end{thm}

Our proof will establish a quantitative version of the theorem, which pins down the potential singularities of $\,T_{\alpha,\eta}f$ in terms of those of $f$ and the parameter $(\alpha,\eta)$. We shall state the more refined version in the cases of interest to us, at the end of this section.

For the applications we need to know the effect of the adjoint of $\,T_{\alpha,\eta}$ on the level of distributions. If $\,\sigma\in\Sch'(\R)$ vanishes to infinite order at $x=0$, and if $f\in \Sch_{\text{\rm{sis}}}(\R)$ has a singularity of type $(\alpha,\eta)$, \propref{xalpha} identifies $\,f\sigma$ as the product of a smooth function and a distribution which vanishes to infinite order at $x=0$; away from the origin, $\,f\sigma$ may be regarded as the product of a Schwartz function and a tempered distribution. It therefore makes sense to integrate $\,f\sigma$ over $\,\R$. Taking linear combinations, we can define the integration pairing
\begin{equation}
\label{talpha2}
f\ \ \mapsto\ \ \int_\R f(x)\,\sigma(x)\,dx\qquad (\,f\in \Sch_{\text{\rm{sis}}}(\R)\,)
\end{equation}
on all of $\Sch_{\text{\rm{sis}}}(\R)$, against any tempered distribution $\,\sigma\in\Sch'(\R)$ which vanishes to infinite order at $x=0$. According to \thmref{fourier}, the Fourier transform $\,\widehat\sigma$ of any such $\,\sigma$ extends canonically across $\,\infty$. We may therefore regard $\,(\sg x)^\eta|x|^{\alpha-1}\widehat\sigma(1/x)\,$ as tempered distri\-bution -- see the discussion following the proof of \lemref{talphalem} -- which vanishes to infinite order at $x=0$. In short,
\begin{equation}
\label{talpha3}
T_{\alpha,\eta}^*(\sigma)(x)\ \ = \ \ (\sg x)^\eta\,|x|^{\alpha-1} \,\widehat\sigma(1/x)
\end{equation}
is a well defined distribution if $\,\sigma\in\Sch'(\R)$ vanishes to infinite order at $x=0$, and (\ref{talpha3}) defines a map $\,T_{\alpha,\eta}^*$ from the space of all such $\,\sigma$ to itself.

\begin{thm}\label{talphathm2}
The operator $\,T_{\alpha,\eta}^*$ is the adjoint of $\,T_{\alpha,\eta}$, in the sense that
\[
\int_\R T_{\alpha,\eta}(f)(x)\,\sigma(x)\,dx \ = \ \int_\R f(x)\,T_{\alpha,\eta}^*(\sigma)(x)\,dx
\]
if $\,\sigma\in\Sch'(\R)$ vanishes to infinite order at the origin and $\,f\in \Sch_{\text{\rm{sis}}}(\R)$
\end{thm}

The proofs of the two theorems occupy most of the remainder of this section. We begin with a decomposition of the space $\,\Sch_{\text{\rm{sis}}}(\R)$, which is formally similar to the classification of regular singularities in the theory of ordinary differential equations. Note that the rule
\begin{equation}
\label{sis1}
(\alpha_1,\eta_1)\ \preceq\ (\alpha_2,\eta_2)\ \ \Longleftrightarrow \ \ \alpha_2-\alpha_1\, \in \,(2\Z+\eta_1+\eta_2)\,\cap\,\Z_{\geq 0}
\end{equation}
defines an order relation on the set $\,\C\times\Z/2\Z$.

\begin{lem}\label{sisdecomp}
Each $\,f\in \Sch_{\text{\rm{sis}}}(\R)$ can be expressed as a sum
\[
f(x)\ \ = \ \ f_0(x)\ +\ {\sum}_{1\leq j\leq m}\ {\sum}_{0\leq \ell\leq n_j}\ (\sg x)^{\eta_j}\,(\log|x|)^\ell\,|x|^{\alpha_j}\,f_{j,\ell}(x)\,,
\]
with $f_0\in\Sch(\R)$ vanishing to infinite order at the origin, with $\,f_{j,\ell}\in\Sch(\R)$, $(\alpha_j,\eta_j)\in \C\times\Z/2\Z$, $m\geq 0$, $n_j\geq 0$, and $(\alpha_i,\eta_i) \npreceq(\alpha_j,\eta_j)$ unless $i=j$. If $f$ has parity $\delta\in\Z/2\Z$ -- i.e., if $f(-x)=(-1)^\delta f(x)$ -- one can choose $f_0\in\Sch_\delta(\R)$, $f_{j,\ell}\in\Sch_{\delta+\eta_j}(\R)$. The $\alpha_j$ and $n_j$ become uniquely determined when one requires that for each $j$, $f_{j,\ell}(0)\neq 0$ for some $\ell$, and that no $f_{j,n_j}$ vanishes to infinite order at the origin. The $f_{j,\ell}$ are unique up to addition of a Schwartz function which vanishes to infinite order at $x=0$.
\end{lem}

\begin{proof}
The definition of $\,\Sch_{\text{\rm{sis}}}(\R)$ provides a decomposition locally, near the origin, which can be made global by means of a suitable cutoff function. Since $(\alpha_j,\eta_j) \preceq(\alpha_i,\eta_i)$ implies
\begin{equation}
\label{sisdecomp1}
(\sg x)^{\eta_i}\,(\log|x|)^\ell\,|x|^{\alpha_i}\,\Sch(\R)\ \subset\ (\sg x)^{\eta_j}(\log|x|)^\ell\,|x|^{\alpha_j}\,\Sch(\R)\,,
\end{equation}
terms can be combined so as to satisfy the conditions on the $(\alpha_j,\eta_j)$. The function $f_0$ is needed in the decomposition only if $m=0$, i.e., if $f$ vanishes to infinite order at the origin. The uniqueness statements follow from the fact that $f(x)$ has an asymptotic expansion as $x\to 0$, which completely determines the Taylor series of the $f_{j,\ell}$ at the origin.
\end{proof}

We suspect that our next statement is known, though we have not been able to find it stated elsewhere.

\begin{lem}\label{charschwartz}
The signed Mellin transform $\,M_\delta$, $\,\delta\in\Z/2\Z$, establishes an isomorphism between $\Sch_\delta(\R)$ and the space of meromorphic functions $H(s)$ whose only singularities are first order poles at points in $(2\Z+\delta)\cap\Z_{\leq 0}$, and which decay rapidly along vertical lines, locally uniformly in $\,\operatorname{Re}s$.
\end{lem}

\begin{proof} Lemma \ref{mellin_of_schwartz} and corollary \ref{melcor1} tell us that $H(s)=M_\delta f\,(s)$, with $f\in\Sch_\delta$, has the properties asserted by the lemma. Inversion of the Fourier transform in (\ref{moderategrowth12}), the parity condition on $f$, and the change of variables $\,x\rightsquigarrow\log x\,$ make it possible to recover $f(x)$, for $x\neq 0$, from $H(s)=M_\delta f\,(s)$:
\begin{equation}\label{charschwartz1}
f(x)\ \ = \ \ \frac{(\sg x)^\delta}{4\pi i}\int_{\operatorname{Re}s=s_0} H(s)\,|x|^{-s}\,ds\qquad (\,x\neq 0\,,\ \ s_0 > 0\,)\,,
\end{equation}
at least if $s_0=1/2$, but then for other $s_0>0$ by a simple contour shift. If $H(s)$ has the required properties, we define its ``signed Mellin inverse" $\Phi_\delta H$ as the integral on the right of (\ref{charschwartz1}), initially as a function on $\,\R-\{0\}$. We should remark that the notation is consistent with our earlier definition (\ref{moderategrowth6}). Since $H(s)$ decays rapidly on the line of integration, we can differentiate under the integral sign, to conclude
\begin{equation}
\label{derivmdH}
\bigl(\textstyle\frac{d^k\ }{dx^k}\Phi_\delta H \bigr)(x)  \ \ = \ \ \displaystyle\frac{(\sg x)^{\delta+k}\,k!}{4\pi i}  \int_{\operatorname{Re}s=s_0}\,\binom{-s}{k}\, H(s)\,|x|^{-s-k} \,ds\,,
\end{equation}
again with $s_0>0$, but otherwise arbitrary. Shifting the contour to the right and using the rapid decay, we get the bound $(\frac{d^k\ }{dx^k}\Phi_\delta H)(x)=O(|x|^{-N})$ for all $k,\,N\geq 0$. That is the Schwartz condition at infinity. To show that $\Phi_\delta H$ has a smooth extension across $x=0$, we now shift the contour to the left. As we do so, we pick up residues when we move across points $s\in (2\Z+\delta)\cap \Z_{\leq 0}$, but only those for which $-s-k\geq 0$, since the other poles are canceled by the zeros of the binomial coefficient: for $N\in\N$, with $2N>k$,
\begin{equation}
\label{derivasy}
\begin{aligned}
\bigl(\textstyle\frac{d^k\ }{dx^k}\Phi_\delta H \bigr)(x) \ \  &=\ \ \frac {k!}2\ {\sum}_{\stackrel{\scriptstyle{k\,\leq\, n < \,2N\, }}{n\equiv \delta\,\operatorname{mod}2}} \ \binom{n}{k}\bigl(\operatorname{Res}_{s=-n}H(s)\bigr)\,x^{n-k}
\\
&+\ \ \displaystyle\frac{(\sg x)^{\delta+k}\,k!}{4\pi i}\,\binom{-s}{k}\,  \int_{\operatorname{Re}s=1/2-2N} H(s)\,|x|^{-s-k} \,ds\,.
\end{aligned}
\end{equation}
The sum provides an asymptotic expansion for $\textstyle\frac{d^k\ }{dx^k}\Phi_\delta H$, because the error term tends to zero faster than $|x|^{2N-k-1}$ as $x\to 0$. These asymptotic expansions are consistent with the identity $\frac{d\ }{dx}(\Phi_\delta H(s)) = \Phi_{\delta+1}((1-s) H(s-1))$, and therefore do define a smooth extension of $\Phi_\delta H$.
\end{proof}

\begin{cor}\label{charschwartzcor}
The signed Mellin transform $\,M_\delta$ establishes an isomorphism between $\,(\sg x)^\eta(\log|x|)^\ell |x|^\alpha \Sch_{\delta+\eta}(\R)$ and the space of meromorphic functions $H(s)$, whose only singularities are $(\ell+1)$-st order poles at points $s=\alpha-n$, with $n\in(2\Z+\delta+\eta)\cap\Z_{\geq 0}$, such that $P(s)H(s)$ has zero residues at all poles, for every polynomial $P(s)$ of degree $\ell-1$, and such that $H(s)$ decays rapidly along vertical lines, locally uniformly in $\,\operatorname{Re}s$.
\end{cor}

The condition on the poles of $H(s)$ can be paraphrased by saying that they are of ``pure order $\ell+1$", i.e., with principal part $a(s-s_0)^{-\ell-1}$ around any pole $s_0$.

\begin{proof}
Since $\,M_\delta((\sg x)^\eta(\log|x|)^\ell |x|^\alpha f)(s)=M_{\delta+\eta}((\log|x|)^\ell f)(s+\alpha)$, it suffices to deal with the case $\alpha=0$, $\eta=0$. We argue by induction on $\ell$, beginning with $\ell=0$ which reduces to \lemref{charschwartz}. For the induction step, we use the identity
\begin{equation}
\label{charschwartzcor1}
\textstyle\frac{d\ }{dt}M_\delta \bigl((\log|x|)^\ell\, f(x)\bigr)\,(s+t)\bigr|_{t=0}  \ \ = \ \ M_\delta \bigl((\log|x|)^{\ell+1}\, f(x)\bigr)\,(s)\,,
\end{equation}
coupled with the observation that differentiation maps the space of meromorphic functions $H(s)$ corresponding to $\ell\geq 0$ isomorphically onto the space corresponding to $\ell+1$. What matters here is the vanishing of the residues and the rapid decay, which excludes constants. Note that rapid decay is preserved by differentiation, as can be shown by means of the Cauchy integral formula.
\end{proof}

Because of (\ref{talphadef'}), we may regard $T_{\alpha,\eta}f$, for $f\in \Sch_{\text{sis}}(\R)$, as distribution with canonical extension across $\infty$. As such, its Mellin transform is defined for $\,\operatorname{Re}s\gg 0$, as is the Mellin transform of $f$. Recall the function $G_\delta(s)$, which was introduced in \lemref{mellin_and_dirichlet}.

\begin{lem}\label{talphamellem}
For $\,f\in \Sch_{\text{sis}}(\R)$, $\,M_\delta(T_{\alpha,\eta}f)\,(s)= (-1)^\delta G_\delta(s)\,M_{\delta+\eta}f\,(s+\alpha)\,$ on the common domain of definition.
\end{lem}

\begin{proof}
Substitution of either $|x|^\beta f(x)$ or $(\sg x)|x|^\beta f(x)$ for $f(x)$ has the same effect on both sides of the identity to be proved. Also, both sides depend holomorphically on $\alpha$. We are therefore free to suppose that $f(x)=(\log|x|)^\ell x^2 g(x)$ for some $\ell\geq 0$ and $g\in\Sch_{\eta+\delta}(\R)$, and that $\,\operatorname{Re}\alpha > 0$. In that case, when we consider $\,\mathcal F T_{\alpha,\eta}f$ and $\,T_{\alpha,\eta}f$ as distributions, the former vanishes to order $k_0=\infty$ at $x=0$ and has an extension across $\infty$ which vanishes there to order $k_\infty=1$. According to lemma \ref{talphalem}, $T_{\alpha,\eta}f$ has a canonical extension across $\infty$. As the Fourier transform of a function in $L^2(\R)\cap L^1(\R)$, $T_{\alpha,\eta}f$ is continuous, hence vanishes to order $k_0=0$ at $x=0$. In particular, $M_\delta$ maps both $\,\mathcal F T_{\alpha,\eta}f$ and $\,T_{\alpha,\eta}f$ into $\O_{pg}(\{0< \operatorname{Re}s<1\})$, in the notation of \S\ref{4nmelsec}. We can now argue exactly as in the proof of \thmref{4nmelthm} and conclude
\begin{equation}
\label{talphamellem1}
\begin{aligned}
M_\delta(T_{\alpha,\eta}f)\,(s)\ &=\ G_\delta(s)\,M_\delta(\mathcal F T_{\alpha,\eta}f)\,(1-s)
\\
&= \ G_\delta(s) \int_\R (\sg x)^\delta \,(-\sg x)^\eta\,|x|^{-1-\alpha-s}\,f(-1/x)\,dx
\\
&= \ (-1)^\delta G_\delta(s) \int_\R (\sg x)^{\delta + \eta}\,|x|^{s +\alpha - 1}\,f(x)\,dx\,,
\end{aligned}
\end{equation}
which is the assertion of the lemma.
\end{proof}

Let us summarize what we have shown so far. The signed Mellin transform $\,M_\delta$ maps $\,\Sch_{\text{sis}}(\R)$ to a space of meromorphic functions $\mathcal M_{\text{sis}}(\C)$ which we are about to define formally. Via $\,M_\delta$, the operator $\,T_{\alpha,\eta}$ corresponds to the map $\,H(s)\mapsto (-1)^\delta G_\delta(s) H(s+\alpha)$ from the space $\mathcal M_{\text{sis}}(\C)$ to itself.

\begin{definition}\label{merosisdef}
In the following, $\mathcal M_{\text{\rm{sis}}}(\C)$ shall denote the space of all meromorphic functions $H(s)$ on the complex plane, such that\newline
\noindent {\rm a)\,} the poles of $H(s)$ lie in a finite number of sets $\{\beta-2n \mid n\in \Z_{\geq 0}\}$, $\beta\in\C$;\newline
\noindent {\rm b)\,} the order of the poles of $H(s)$ is uniformly bounded; and\newline
\noindent {\rm a)\,} $H(s)$ decays rapidly along vertical lines, locally uniformly in $\,\operatorname{Re}s$.
\end{definition}

To complete the proof of \thmref{talphathm1}, we need a decomposition of the space $\mathcal M_{\text{\rm{sis}}}(\C)$ analogous to the decomposition of $\,\Sch_{\text{sis}}(\R)$ given in \lemref{sisdecomp}. This depends on certain Gamma identities. Recall that $\Gamma(s)$ is a meromorphic function which has a first order pole at every nonpositive integer, but has no other poles and no zeros. Stirling's formula provides an asymptotic expansion of $\,|\Gamma(s)|\,$ along vertical lines, whose first term describes the asymptotic behavior:
\begin{equation}
\label{stirling}
|\G(s)| \ \ \sim  \ \ \sqrt{2\pi}\,\,|\operatorname{Im}s|^{\operatorname{Re}s-1/2}\, e^{-\pi |\operatorname{Im}s|/2}\qquad\text{as}\ \ |\operatorname{Im}s|\rightarrow\infty\,.
\end{equation}
This, in conjunction with the Cauchy integral formula, implies bounds for the derivatives:
\begin{equation}
\label{stirdif}
\begin{aligned}
|\Gamma(s)|^{-1}\,|\Gamma^{(k)}(s)|\ \ \text{has polynomial growth along vertical}
\\
\text{lines, locally uniformly in $\operatorname{Re}s$}\,,
\end{aligned}
\end{equation}
for every $k\in\Z_{>0}$. Note that $\,\Gamma^{(k)}(s)$ has poles at the same points as $\,\Gamma(s)$, but the poles of $\,\Gamma^{(k)}(s)$ have ``pure order $k+1$", in the sense that the product $\,P(s)\,\Gamma^{(k)}(s)$ with any polynomial $P(s)$ of degree $k-1$ has zero residues.

\begin{lem}\label{gamlem}
Let $k_1$, $k_2$ be non-negative integers, and $\beta$ a complex number.\newline
\noindent {\rm a)\,} If $\beta\notin \Z$, there exist entire functions $F_{j,\ell}(s)$ such that $\,\Gamma(s)^{-1}F_{j,\ell}(s)$ has locally uniform polynomial growth on vertical lines, and
\[
\G^{(k_1)}(s)\,\G^{(k_2)}(s+\beta) \ =  \sum_{0\leq\ell\leq k_1}\! F_{1,\ell}(s+\beta)\,\G^{(\ell)}(s) \ + \sum_{0\leq\ell\leq k_2}\! F_{2,\ell}(s)\,\G^{(\ell)}(s+\beta)\,.
\]
\noindent {\rm b)} If $\beta=n\in\Z_{\ge 0}$, there exist entire functions $F_{j,\ell}(s)$ such that $\Gamma(s)^{-1}F_{j,\ell}(s)$ has locally uniform polynomial growth on vertical lines, and
\[
 \G^{(k_1)}(s)\,\G^{(k_2)}(s+n)\ =  \sum_{0\leq\ell\leq k_1} F_{1,\ell}(s)\G^{(\ell)}(s)\ + \sum_{k_1<\ell\leq k_1+k_2+1} F_{2,\ell}(s)\G^{(\ell)}(s+n)\,.
\]
\end{lem}

\begin{proof}
We begin with the two equivalent identities
\begin{equation}
\label{whoa}
\begin{aligned}
\G(s)\,\G(s+\beta)\, =\, \f{\pi}{ \sin(\pi\, \beta)}\( \f{e^{i \,\pi\,s}}{ \G(1-s-\beta)} \G(s)  -  \f{e^{i \,\pi(s+\beta)} }{\G(1-s)} \G(s+\beta) \)
\\
= \, \f{\pi}{ \sin(\pi\,\beta)}\( \f{e^{-i \,\pi\, s}}{\G(1-s - \beta)} \G(s) - \f{e^{-i \,\pi(s+\beta)}}{\G(1-s)}\G(s+\beta)\),
\end{aligned}
\end{equation}
which can be verified by multiplying both sides with $\G(1-s)\,\G(1-s-\beta)$. Since $\G(s)\G(1-s)=\pi/(\sin \pi s)$, the first identity then becomes equivalent to the trigonometric identity $\sin(\pi\,\beta) = e^{i \pi s} \sin(\pi( s+\beta)) - e^{i \pi (s+\beta)} \sin(\pi s)$. The second follows from the first because $\,\G(s)$ is real on the real axis. We re-write the first identity as
\begin{equation}
\label{newhoa}
\begin{aligned}
\G(s)\,\G(s+\beta) \ \ &= \ \ c(\beta)\,h(s+\beta)\,\G(s) \ + \ c(-\beta) \,h(s)\,\G(s+\beta)\,,
\\
\text{with}\ \ \ h(s) \ &=\ \f{e^{i\pi s}}{\G(1-s)}\ \ \ \text{and}\ \ \ c(\beta)\ =\  \f{\pi\,e^{-i \pi \beta}}{\sin(\pi \beta)}\,.
\end{aligned}
\end{equation}
Then $c(\beta)$ is periodic, meromorphic, and has only simple poles with residue 1 at the integers. The function $h(s)$ is entire. In the upper half plane only, $\,\Gamma(s)^{-1}h^{(k)}(s)$ has locally uniform polynomial growth along vertical lines, as follows from (\ref{stirdif}). In the lower half plane $\,\Gamma(s)^{-1}h^{(k)}(s)$ grows exponentially; more precisely, $\,\Gamma(s)^3 h^{(k)}(s)$ has locally uniform polynomial growth along vertical lines in the lower half plane. Had we used the second equation in (\ref{whoa}) instead of the first, we would have obtained the same type of expression, with $\bar h(\bar s)$ and $\bar c(\bar \beta)$ in place of $h(s)$ and $c(\beta)$, which would have resulted in exponential decay in the lower half plane and exponential growth in the upper half plane. We now suppose $\beta\notin\Z$. We temporarily treat $s_1=s$ and $s_2=s+\beta$ as independent variables, apply $(\f{d\ }{ds_1})^{k_1}(\f{d\ }{ds_2})^{k_2}$ to the identity (\ref{newhoa}), then substitute back $s$ and $\beta$. The result is an expression for $\G^{(k_1)}(s)\,\G^{(k_2)}(s+\beta)$ as a finite sum of products of derivatives of $\,c(\beta)$, $c(-\beta)$, $h(s)$, $h(s+\beta)$, $\G(s)$, and $\G(s+\beta)$. In view of the growth properties of the derivatives of $h(s)$ mentioned above, this exhibits $\G^{(k_1)}(s)\,\G^{(k_2)}(s+\beta)$ as a sum of the form asserted in the lemma, but with the required growth properties of the $F_{j,\ell}(s)$ satisfied only in the upper-half plane, and subject to the weaker condition of locally uniform polynomial growth of $\G(s)^3F_{j,\ell}(s)$ along vertical lines in the lower half plane.

As was just remarked, the upper and lower half planes play symmetric roles. We therefore get another expression of the same type, with coefficient functions with locally uniform polynomial growth in the lower half plane. To blend the two expressions, we use an ``analytic partition of unity'' created from the classical error function
\begin{equation}
\label{erfdef}
\operatorname{erf}(s) \ \ =  \ \ \f{2}{\sqrt{\pi}}\, \int_0^s \, e^{-z^2}\, dz\,.
\end{equation}
Note that $\,\operatorname{erf}(s)$ is an entire function, $\,\operatorname{erf}(-s) = -\operatorname{erf}(s)$, and $\,\operatorname{erf}(s)\to 1$ as $s$ approaches $\infty$ along the positive real axis. Simple estimates imply
\begin{equation}
\label{erfest}
\begin{aligned}
|1 - \operatorname{erf}(s)| \ \ =  \ \ O(e^{-(\operatorname{Re}s)^2})\ \ \ \text{as}\ \ \ \operatorname{Re}s\to +\infty\,,
\\
\text{locally uniformly in $\operatorname{Im}s$}\,;
\end{aligned}
\end{equation}
for details see \cite{Lebed}, for example. The related function $E(s)=\f 12(1+\operatorname{erf}(-i s))$ tends to $1$ as $s\to\infty$ along the positive imaginary axis, and to $0$ as $s\to\infty$ along the negative imaginary axis. In fact,
\begin{equation}
\label{Eest}
\begin{aligned}
|1 - E(s)| \ \ =  \ \ O(e^{-(\operatorname{Im}s)^2})\ \ \ \text{as}\ \ \ \operatorname{Im}s\to +\infty\,,
\\|E(s)| \ \ =  \ \ O(e^{-(\operatorname{Im}s)^2})\ \ \ \text{as}\ \ \ \operatorname{Im}s\to -\infty\,,
\\
\text{in both cases locally uniformly in $\operatorname{Re}s$}\,;
\end{aligned}
\end{equation}
this follows from (\ref{erfest}). We now take the expressions for $\G^{(k_1)}(s)\,\G^{(k_2)}(s+\beta)$ which we had derived, multiply the first -- i.e., the one which has the required growth behavior in the upper half plane -- with $E(s)$, and the second with $(1-E(s))$, then add the two. Because of (\ref{stirling}) and (\ref{Eest}), the resulting expression has the properties asserted in the first part of the lemma.

The proof of the second part is similar, though slightly more involved. Recall that $\beta\mapsto c(\beta)$ is periodic of period $1$ and has a first order pole at $\beta=0$, with residue $1$. We subtract off the pole, to make $\tilde c(\beta)=c(\beta)-1/\beta$ regular at the origin. Specializing (\ref{newhoa}) we now find
\begin{equation}
\label{gamlem1}
\begin{aligned}
&\G(s+s_1)\,\G(s+n+s_2)  \ =  \ \tilde c(s_2-s_1) \, h(s+n+s_2)\, \G(s+s_1) \ +
\\
&\ \ \ \ \ +\ \tilde c(s_1-s_2)\,h(s+s_1)\,\G(s+n+s_2)
\\
&\ \ \ \ \ +\ \f{1}{s_2-s_1} \bigl[ h(s+n+s_2)\, \G(s+s_1) \, - \,
\,h(s+s_1)\,\G(s+n+s_2)\bigr].
\end{aligned}
\end{equation}
We expand both sides of this identity as a Taylor series in powers of $s_1$ and $s_2$, then equate the coefficient of $s_1^{k_1}s_2^{k_2}$ and clear out the denominator $k_1!\,k_2!\,$. On the left hand side this gives us $\G^{(k_1)}(s)\G^{(k_2)}(s+n)$, which is the left hand side of the identity we want to prove. We shall show that this process, applied to the right hand side of (\ref{gamlem1}), gives us an expression of the type asserted in the lemma, but with coefficients $F_{j,\ell}(s)$ which satisfy the appropriate bound only in the upper half plane and the weaker condition of locally uniform polynomial growth of $\G(s)^3F_{j,\ell}(s)$ in the lower half plane. The two summands containing $\tilde c$ contribute the type of terms we expect, except that the summation in the second sum extends over $0\leq \ell \leq k_2$, instead of $k_1<  \ell \leq k_1 + k_2 + 1$ as claimed. Terms corresponding to $0 \leq \ell  \le k_1$, if any, can be absorbed by the first sum, thanks to the identity
\begin{equation}
\label{gamlem2}
\G^{(\ell)}(s+n) \ = \ {\sum}_{0\leq j\leq \ell}\ {\ell \choose j}\, \G^{(\ell-j)}(s)
\,\f{d^j\ }{ds^j}\bigl(s(s+1)\cdots(s+n-1)\bigr),
\end{equation}
which follows from the standard identity $\,\G(s+1)=s\G(s)$ by induction on $n$ and differentiation. Next we expand the last term on the right of (\ref{gamlem1}) as a Taylor series:
\begin{equation}
\label{gamlem3}
\begin{aligned}
\frac 1{s_2-s_1}\left[ \, \dots\,  - \, \dots \,  \right]  \, \ = \, \ \frac 1{s_2-s_1}\,{\sum}_{j_1,\,j_2\geq 0}\ b_{j_1,j_2}(s)\,s_1^{j_1}\,s_2^{j_2}\,,\ \ \ \text{with}
\\
b_{j_1,j_2}(s) \ \ = \ \ \f{1}{j_1!\, j_2!}\,\bigl[ h^{(j_2)}(s+n)\,\G^{(j_1)}(s) \,  - h^{(j_1)}(s)\,\G^{(j_2)}(s+n)\bigr].
\end{aligned}
\end{equation}
The series is formally divisible by $s_2-s_1$ because the other terms in (\ref{gamlem1}) have no singularity along the hyperplane $s_2=s_1$:
\begin{equation}
\label{gamlem4}
\frac 1{s_2-s_1}\,{\sum}_{j_1,\,j_2\geq 0}\ b_{j_1,j_2}(s)\,s_1^{j_1}\,s_2^{j_2}\ \ = \ \ {\sum}_{j_1,\,j_2\geq 0}\ a_{j_1,j_2}(s)\,s_1^{j_1}\,s_2^{j_2}\,.
\end{equation}
This implies $\,\sum_{0\leq j \leq k}\,b_{j,k-j}(s)=0\,$ for all $\,k$, and
\begin{equation}
\label{gamlem5}
a_{k_1,k_2}(s)\ \ = \ \ {\sum}_{0\leq \ell\leq k_1}\ b_{\ell,k_1+k_2+1-\ell}(s) \,.
\end{equation}
When we combine (\ref{gamlem5}) with (\ref{gamlem3}), we almost get the expression we want. Derivatives of $\,\G(s+n)$ of order $\ell\leq k_1$ constitute the only remaining obstacle, but they can be absorbed into the first sum, as before. To complete the proof, we repeat the ``partition of unity" argument used in proving part a) to construct an expression of the required type which has the appropriate growth behavior in both half planes.
\end{proof}

\begin{lem}\label{merodecomp}
Every $H\in\mathcal M_{\text{\rm{sis}}}(\C)$ can be expressed as a sum
\[
H(s) \ = \ H_0(s)\ + \ {\sum}_{1\leq j\leq m}\ {\sum}_{0\leq \ell\leq n_j} H_{j,\ell}(s)\,,
\]
in terms of an entire function $H_0(s)$ and meromorphic functions $H_{j,\ell}(s)$, satisfying the following properties: there exist $\beta_1,\,\beta_2,\,\dots,\,\beta_m\in\C$, such that $\beta_i-\beta_j\notin 2\Z$ for $i\neq j$ and\newline
\noindent {\rm a)\,} $H_{j,\ell}$ has poles only at the points in $\beta_j-2\Z_{\leq 0}$, all of order $\ell+1$;\newline
\noindent {\rm b)\,} $P(s)\,H_{j,\ell}(s)$ has zero residues for every polynomial $P(s)$ of degree $\ell-1$;\newline
\noindent {\rm c)\,} $H_0$ and the $H_{j,\ell}$ decay rapidly on vertical lines, locally uniformly in $\,\operatorname{Re}s$.\newline
\noindent The $\beta_j$ and $n_j$ become uniquely determined when one requires $H(s)$ to have an actual pole at each $\beta_j$, and when for each $j$, $H(s)$ has a pole of order exactly $n_j$ at $s=\beta_j - 2k$, for some $k\geq 0$. The $H_{j,\ell}(s)$ are unique up to addition of an entire function which decays rapidly along vertical lines.
\end{lem}

The condition on the poles of $H_{j,\ell}(s)$ means that they are of ``pure order $\ell + 1$", just as in the case of corollary \ref{charschwartzcor}.

\begin{proof}
The definition of $\mathcal M_{\text{\rm{sis}}}(\C)$ implies the existence of distinct complex numbers $\beta_j$, no two of which differ by an even integer, and $n_j\geq 0$, such that
\begin{equation}
\label{merodecomp1}
F(s)\, =_{\text{def}}\, \bigl( \ {\prod}_{1\leq j\leq m}\, \G(2s+\beta_j)^{n_j+1} \,\bigr)^{-1}H(s)\, \ \ \text{is an entire function.}
\end{equation}
Repeated application of the two identities in \lemref{gamlem} -- the second in particular with $n=0$ -- makes it possible to separate the poles of the product of Gamma functions: there exist entire functions $F_{j,\ell}(s)$ such that
\begin{equation}
\label{merodecomp2}
{\prod}_{1\leq j\leq m} \G(2s+\beta_j)^{n_j+1} \, = \, {\sum}_{1\leq j\leq m} {\sum}_{0\leq \ell\leq n_j} F_{j,\ell}(s)\,\G^{(\ell)}(2s+\beta_j)\,,
\end{equation}
with $\,\G(s)^{1-\sum_j (n_j+1)}F_{j,\ell}(s)\,$ having polynomial growth on vertical lines, locally uniformly in $\,\operatorname{Re}s$. The functions $\widetilde H_{j,\ell}(s)=F(s)F_{j,\ell}(s)\G^{(\ell)}(2s+\beta_j)$ add up to $H(s)$ and satisfy the condition a). According to (\rangeref{stirling}{stirdif}), all the $\,\G(2s+\beta_j)$ and their derivatives have the same type of growth behavior, up to polynomial growth. The $\widetilde H_{j,\ell}(s)$ therefore inherit the rapid decay from $H(s)$; i.e., they satisfy c) as well. We still need to modify the $\widetilde H_{j,\ell}(s)$ to establish the condition b). At this point, we may as well suppose that $m=1$ and $\beta_1=0$. We shall argue by induction on $n_1=n$. For $n=0$ the condition b) holds vacuously. For $n>0$, we set
\begin{equation}
\label{merodecomp3}
\begin{aligned}
H_{1,n}(s)\ &= \ 2^{-n}\,\textstyle\frac{d^n\ }{ds^n}\bigl( F(s)F_{1,n}(s)\G(2s) \bigr)
\\
&= \ \widetilde H_{1,n}(s)\ + \ \textstyle{\sum}_{1\leq i \leq n}\ 2^{-i} \binom{n}{i}\,\frac{d^i\ }{ds^i}\bigl( F(s)F_{1,n}(s)\bigr)\,\G^{(n-i)}(2s)\,.
\end{aligned}
\end{equation}
Then $H_{1,n}(s)$ has poles of ``pure order $n+1$", and $H(s)-H_{1,n}(s)$ has poles of order at most $n$. Differentiation does not change the order of growth or decay, so we have reduced $n$ by $1$ without affecting the other hypotheses. That completes the inductive argument.

The function $H_0$ is needed only when $m=0$, i.e., when $H(s)$ is entire. The conditions a), b) determine the principal part of each $H_{j,\ell}(s)$ at each of its poles, and that makes each $H_{j,\ell}(s)$ unique up to addition of an entire function. The uniqueness statement about the $\beta_j$ and $n_j$ is correct for purely formal reasons.
\end{proof}

\begin{cor}\label{melsis}
The signed Mellin transform $\,M_\delta$ induces an isomorphism
\[
M_\delta \, : \, \{\,f\in\Sch_{\text{\rm{sis}}}(\R) \mid f(-x)=(-1)^\delta f(x)\,\} \ \overset{\sim}{\longrightarrow} \ \mathcal M_{\text{\rm{sis}}}(\C)\,.
\]
\end{cor}

\begin{proof}
According to corollary \ref{charschwartzcor}, $\,M_\delta$ relates the decomposition of the space $\{\,f\in\Sch_{\text{\rm{sis}}}(\R) \mid f(-x)=(-1)^\delta f(x)\,\}$ in \lemref{sisdecomp} to the decomposition of $\,\mathcal M_{\text{\rm{sis}}}(\C)$ in \lemref{merodecomp}, and $\,M_\delta$ induces isomorphisms between components on the two sides that correspond to each other. Lemmas \ref{mellin_of_schwartz} and \ref{charschwartz} ensure that $\,M_\delta$ also relates the ambiguities in the two decompositions bijectively.
\end{proof}

The proof of \thmref{talphathm1} is essentially complete. For $f\in\Sch_{\text{\rm{sis}}}(\R)$, $\,T_{\alpha,\eta}(f)$ exists at least as distribution with canonical extension across $\infty$. Corollary \ref{melsis} and \lemref{talphamellem} guarantee the existence of some $\widetilde f\in \Sch_{\text{\rm{sis}}}(\R)$ such that $\,M_\delta T_{\alpha,\eta}(f)= M_\delta\widetilde f\,$ for both choices of $\,\delta$. But a distribution $\,\sigma\in\Sch'_\delta(\R)$, with canonical extension across $\infty$, is determined by $\,M_\delta\sigma$ up to a distribution supported at the origin -- see the discussion around (\ref{autdist16}). We conclude that $\,T_{\alpha,\delta}(f)$ agrees with $\,\widetilde f$ as function on $\,\R-\{0\}$, and consequently as function in $\,\Sch_{\text{\rm{sis}}}(\R)$. That, in effect, is the assertion of \thmref{talphathm1}.

We begin the proof of \thmref{talphathm2} with another lemma. The signed Mellin transform $\,M_\delta f\,(s)$ of any $\,f\in \Sch_{\text{\rm{sis}}}(\R)$ is regular for $\,\operatorname{Re}(s)\gg 0$ and decays rapidly along vertical lines. In section \S \ref{4nmelsec}, we saw that $\,M_\delta \sigma\,(s)$ is regular for $\,\operatorname{Re}(s)\ll 0$, provided $\,\sigma\in\Sch'_\delta (\R)$ vanishes to infinite order at the origin; moreover, $\,M_\delta \sigma\,(s)$ grows polynomially along vertical lines.

\begin{lem}\label{adjointlem}
If $\,\sigma\in\Sch'_\delta (\R)$ vanishes to infinite order at the origin,
\[
\int_\R f(x)\,\sigma(x)\, dx\ = \ \frac{1}{4\pi i}\int_{\operatorname{Re}(s)=s_0} M_\delta f\,(s)\,\,M_\delta \sigma\,(1-s)\,ds \qquad (\, s_0 \gg 0\,)\,,
\]
for any $\,f\in \Sch_{\text{\rm{sis}}}(\R)$. The abscissa of integration $\,s_0$ must be chosen so that the integrand is regular on, and to the right of, the line of integration, but $\,s_0$ is otherwise arbitrary.
\end{lem}

\begin{proof}
According to the comments before the statement of the lemma, the integrand is indeed regular on some right half plane and decays rapidly along vertical lines, as always locally uniformly in $\,\operatorname{Re}(s)$. It follows that the integral on the right converges and does not depend on the particular choice of $\,s_0$. Both sides of the equation vanish when $f$ has the parity opposite to $\,\delta$. We therefore may and shall assume that $\,f(-x)=(-1)^\delta f(x)$. Lemma \ref{sisdecomp} allows us to also suppose
\begin{equation}
\label{adjeq1}
f(x)\ =\ (\sg x)^\eta \, |x|^{\alpha}\, (\log|x|)^{n}\,g(x)\ \ \ \text{with}\ \ g\, \in \, \Sch_{\delta+\eta}(\R)\,,
\end{equation}
for some $(\alpha,\eta)\in \C \times \Z/2\Z$ and $n\geq 0$. In this situation,
\begin{equation}
\label{adjeq2}
\int_\R f(x)\,\sigma(x)\, dx\ = \ \int_\R g(x) \bigl((\sg x)^\eta \, |x|^{\alpha}\, (\log|x|)^{n}\,\sigma(x)\bigr) \,dx\,,
\end{equation}
by definition of the pairing (\ref{talpha2}). Since $\,\sigma$ vanishes to infinite order at the origin, so does the tempered distribution $\,(\sg x)^\eta  |x|^{\alpha} (\log|x|)^{n}\sigma(x)$. We can replace $\,\alpha$ by $\,\alpha-2m$, for any $\,m>0$, at the expense of making $\,g$ vanish to high order at the origin. By doing so we can make sure that $\,(\sg x)^\eta  |x|^{\alpha} (\log|x|)^{n}\sigma(x)$ has an extension across $\infty$ which vanishes there to order $\,k_\infty \geq 1$, in which case \lemref{mellem4} applies:
\begin{equation}
\label{adjeq3}
\begin{aligned}
&\int_\R g(x) \bigl((\sg x)^\eta \, |x|^{\alpha}\, (\log|x|)^{n}\,\sigma(x)\bigr) \,dx \ =
\\
&\ \ \ \,=\ \frac{1}{4\pi i}\int_{\operatorname{Re}(s)=s_0} \!\!\!\!\!\! M_{\delta+\eta} g\,(s)\,\,M_{\delta+\eta} \bigl((\sg x)^\eta  |x|^{\alpha} (\log|x|)^{n}\sigma\bigr)(1-s)\,ds \,,
\end{aligned}
\end{equation}
for $\,s_0 \gg 1$. Going back to the definition of $\,M_\delta$, one finds
\begin{equation}
\label{adjeq4}
M_{\delta+\eta} \bigl((\sg x)^\eta  |x|^{\alpha} (\log|x|)^{n}\sigma\bigr)(s)\ = \  \textstyle\frac{d^n\ }{ds^n}M_{\delta} \sigma(s+\alpha)\,,
\end{equation}
and the relation (\ref{adjeq1}) between $\,f$ and $\,g$ translates into the relation
\begin{equation}
\label{adjeq5}
M_\delta f\,(s)\ = \ \textstyle\frac{d^n\ }{ds^n} M_{\delta + \eta}g\,(s + \alpha)
\end{equation}
between their Mellin transforms. The identity asserted by the lemma follows from (\rangeref{adjeq2}{adjeq5}) when we translate the line of integration by $\,\alpha$ and perform an $n$-fold integration by parts on the right hand side of (\ref{adjeq3}).
\end{proof}

In proving \thmref{talphathm2}, we may as well suppose that each of $\,f$ and $\,\sigma$ is either even or odd. Since  $\,T_{\alpha,\eta}$ and $\,T^*_{\alpha,\eta}$ change the parity by $\,\eta$, the identity to be proved holds vacuously unless the parities of $\,f$ and $\,\sigma$ are related by $\,\eta$. Thus, from now on,
\begin{equation}
\label{adjeq6}
f \in \Sch_{\text{\rm{sis}}}(\R)\,,\ \ \ f(-x)\ = \ (-1)^{\delta+\eta} f(x)\,,\ \ \ \text{and}\ \ \ \sigma \in \Sch'_{\delta}(\R)\,.
\end{equation}
Letting $\,T_{\alpha,\eta}(f)$ play the role of $\,f$ in \lemref{adjointlem}, we find
\begin{equation}
\label{adjeq7}
\begin{aligned}
\int_\R T_{\alpha,\eta}(f)\, \,\sigma \,dx \ &=\ \frac{1}{4\pi i}\int_{\operatorname{Re}(s)=s_0} \!\!\!\! M_{\delta}\bigl( T_{\alpha,\eta}(f)\bigr)\,(s)\,\,M_{\delta}\sigma\,(1-s)\,ds
\\
&=\ \frac{(-1)^\delta}{4\pi i}\int_{\operatorname{Re}(s)=s_0} \!\!\!\!\!\!\!\!\!\! G_\delta(s)\,M_{\delta+\eta} f\,(s+\alpha)\,\,M_\delta \sigma\,(1-s)\,ds \,;
\end{aligned}
\end{equation}
at the second step we have used \lemref{talphamellem}. Since $\,T^*_{\alpha,\eta}\sigma\in \Sch'_{\alpha+\delta}(\R)$ is known to vanish to infinite order at the origin, we can also apply \lemref{adjointlem} with $\,T^*_{\alpha,\eta}\sigma$ in place of $\,\sigma\,$:
\begin{equation}
\label{adjeq8}
\int_\R f\, \,T^*_{\alpha,\eta}(\sigma) \,dx \ =\ \frac{1}{4\pi i}\int_{\operatorname{Re}(s)=s_0} \!\!\!\! M_{\delta+\eta}f(s)\,M_{\delta+\eta}\bigl( T^*_{\alpha,\eta}(\sigma)\bigr)\,(1-s)\,ds \,.
\end{equation}
The equations (\rangeref{adjeq7}{adjeq8}) reduce the assertion of the theorem to the identity
\begin{equation}
\label{adjeq9}
M_{\delta+\eta}\bigl(T^*_{\alpha,\eta}(\sigma)\bigr)\,(1-s) \ =\ (-1)^\delta \, G_{\delta}(s-\alpha)\, M_\delta \sigma(1-s+\alpha) \,,
\end{equation}
which would follow from \thmref{4nmelthm} if not only $\,\sigma$, but also $\,\widehat\sigma$ were to vanish to infinite order at $x=0$. Since we cannot make that assumption we must argue differently. The equations (\rangeref{adjeq7}{adjeq8}) hold in particular when $\,f$ is a Schwartz function which vanishes identically near the origin. In that case, the Fourier integral in the definition (\ref{talphadef}) of $\,T_{\alpha,\eta}f$ has a Schwartz function as argument. One can then prove \thmref{talphathm2}, for the particular choice of $\,f$, by direct computation, using only the Parseval identity (\ref{fourier3}) and the change of variables formula (\ref{changeofvariables}). We conclude that (\ref{adjeq9}) becomes valid when integrated against the Mellin transform of any $\,f\in\Sch_{\delta+\eta}(\R)$ which vanishes identically near $x=0$:
\begin{equation}
\label{adjeq10}
\begin{aligned}
&\int_{\operatorname{Re}(s)=s_0} M_{\delta+\eta}\bigl(T^*_{\alpha,\eta}(\sigma)\bigr)\,(1-s) \,\, M_{\delta+\eta}f\,(s)\,ds \ \ =
\\
&\ \ \ \ =\ \ (-1)^\delta \int_{\operatorname{Re}(s)=s_0}  G_{\delta}(s-\alpha)\,\, M_\delta \sigma(1-s+\alpha)\,\,M_{\delta+\eta}f\,(s)\,ds \,,
\end{aligned}
\end{equation}
provided $\,s_0\gg 0$. But $\,M_{\delta+\eta}f\,(s_0+iy)=2\mathcal F\bigl(e^{s_0 x}f(e^x)\bigr)(-y/2\pi)$, as can be seen by arguing as in (\ref{moderategrowth12}). Any $\,h\in C^\infty_c(\R)$ can play the role of $\,x\mapsto~e^{s_0 x}f(e^x)$. Thus, when we take the difference of the two expressions in (\ref{adjeq9}) and substitute $\,s=s_0+iy$, the resulting function of the variable $y$ -- which may be regarded as a tempered distri\-bution -- is perpendicular to $\,\mathcal F h$, for all test functions $\,h \in C^\infty_c(\R)$. That is possible only if the identity (\ref{adjeq9}) holds along the vertical line $\,\operatorname{Re}(s)=s_0$, or equivalently, for all $\,s\in \C$. The proof of \thmref{talphathm2} is now complete.

Our statement of the Voronoi summation formula for $GL(3)$ involves the integral transform operator
\begin{equation}
\label{talpha4}
\begin{aligned}
&(\sg x)^{\delta_3}\,|x|^{\lambda_3}\,\Sch(\R)\ \ni \ f \ \mapsto \ F \ \in \Sch_{\text{sis}}(\R)\,,
\\
&\ \ \ \ F\ =\ (\sg x)^{\delta_1}|x|^{1-\lambda_1} T_{\lambda_1-\lambda_2,\delta_3}\circ T_{\lambda_2-\lambda_3,\delta_1}\circ \mathcal F\bigl( (\sg x)^{\delta_3}|x|^{-\lambda_3} f\bigr)\,,
\end{aligned}
\end{equation}
which depends on the parameters $\,(\lambda_j,\,\delta_j)\in \C\times\Z/2\Z\,$, $\,1\leq j\leq 3\,$ \cite{MS2}. The passage from $\,f$ to $\,F$ does not affect the parity, so we may as well suppose that $\,f(-x)=(-1)^\eta f(x)$ and $\,F(-x)=(-1)^\eta F(x)$, in which case the Mellin transform $M_\eta F$ completely determines $\,F$. According to \lemref{talphamellem},
\begin{equation}
\label{talpha5}
\begin{aligned}
M_\eta F\,(s) \ = \ (-1)^{\delta_3}\,G_{\delta_1+\eta}(s-\lambda_1+1)\,G_{\delta_2+\eta}(s-\lambda_2+1)\ \times\qquad
\\
\times \ M_{\delta_3+\eta}h\,(s-\lambda_3+1)\,,
\end{aligned}
\end{equation}
where $\,h = \mathcal F\bigl( (\sg x)^{\delta_3}|x|^{-\lambda_3} f\bigr)\in\Sch_{\delta_3+\eta}(\R)$. Each of the three factors on the right has only first order poles, only at points in
\begin{equation}
\label{talpha6}
s \ \in \ \lambda_j \, -\, 1 \, + \,  (2\Z + \eta + \delta_j) \cap \Z_{\leq 0} \qquad \bigl(\, j=1,\,2,\,3\, \bigr)\,.
\end{equation}
Recall the definition (\ref{sis1}) of the partial order $\,\preceq\,$. In the generic situation, when no two of the pairs $(\l_j,\d_j)$ are related by the order, the poles of the three factors do not overlap. Hence, in view of corollary \ref{charschwartzcor},
\begin{equation}
\label{talpha7}
F \, \in \, {\sum}_{1\leq j\leq 3}\ (\sg x)^{\d_j}|x|^{1-\l_j}\Sch(\R) \ \ \ \text{if}\,\ (\l_i,\d_i)\npreceq (\l_j,\d_j)\,\ \text{for}\,\ i\neq j  \,,
\end{equation}
independently of the particular value of $\,\eta$. As far as the location of the poles is concerned, the three pairs $\,(\lambda_j,\,\delta_j)$ play completely symmetric roles in (\ref{talpha5}). Thus, if exactly two pairs are related by $\,\preceq\,$, we may as well suppose that $(\l_1,\d_1) \preceq (\l_2,\d_2)$. The poles at points $\,s\in\lambda_1 - 1  + (2\Z + \eta + \delta_1) \cap \Z_{\leq 0}$ can then have order two, but all other poles are still simple. In this situation,
\begin{equation}
\label{talpha8}
\begin{aligned}
F \ \in \ (\sg x)^{\d_1}|x|^{1-\l_1}\log{|x|}\,\Sch(\R)\ + \  {\sum}_{2\leq j\leq 3}\, (\sg x)^{\d_j}|x|^{1-\l_j}\Sch(\R)
\\
\text{if}\ \ (\l_1,\d_1) \preceq (\l_2,\d_2)\ \ \text{and}\ \ (\l_j,\d_j) \npreceq (\l_3,\d_3)\ \ \text{for}\ \ j=1,\,2\,.
\end{aligned}
\end{equation}
At first glance, $\,(\sg x)^{\d_1}|x|^{1-\l_1}\Sch(\R)\,$ may also contribute, but this space is contained in $\,(\sg x)^{\d_2}|x|^{1-\l_2}\Sch(\R)\,$ because $(\l_1,\d_1) \preceq (\l_2,\d_2)$. In the only remaining case the three pairs must be linearly ordered. Appealing to the symmetry among the $\,(\lambda_j,\,\delta_j)$ once again, we may suppose that the order increases with increasing $\,j$. Then
\begin{equation}
\label{talpha9}
\begin{aligned}
F \ \in \ {\sum}_{1\leq j\leq 3}\ (\sg x)^{\d_j}\,|x|^{1-\l_j}\,(\log{|x|})^{3-j}\,\Sch(\R)
\\
\text{if}\ \ (\l_1,\d_1) \preceq (\l_2,\d_2) \preceq (\l_3,\d_3)\,,
\end{aligned}
\end{equation}
since the poles at points $\,s\in\lambda_j - 1  + (2\Z + \eta + \delta_j) \cap \Z_{\leq 0}$, for $1\leq j\leq 3$, can have order up to $3-j$. The comment following (\ref{talpha8}) applies triply in the current setting.

The introduction to \cite{MS2} sketches a proof the Voronoi summation formula for $SL(2)$, which has a long history \cite{MS1}. Our formulation involves the $SL(2)$ analogue of the integral transform (\ref{talpha4}),
\begin{equation}
\label{talpha10}
|x|^{-\nu}\,\Sch(\R)\ \ni \ f \ \ \mapsto \ \ F\ = \ |x|^{1-\nu}\, T_{2\nu,0}\circ  \mathcal F\bigl(\, |x|^\nu f\bigr)\ \in \ \Sch_{\text{sis}}(\R)\,,
\end{equation}
with $\,\nu\in\C$. This case is simpler, of course. One can argue as before, to find
\begin{equation}
\label{talpha11}
F \ \in \
\begin{cases}
\ |x|^{1-\nu}\Sch(\R)\ +\ |x|^{1+\nu}\Sch(\R) \ &\text{if}\,\ \nu\notin\Z  \,,
\\
\ |x|^{1-\nu}\log|x|\,\Sch(\R)\ +\ |x|^{1+\nu}\Sch(\R) \ &\text{if}\,\ \nu\in\Z_{\leq 0}  \,.
\end{cases}
\end{equation}
As in the previous case, we can interchange $\,\nu$ and $-\nu$ in deriving (\ref{talpha11}), even though they do not occur symmetrically in the definition of the integral transform (\ref{talpha10}).

\section{The multi-variable case revisited}\label{multiagain}

We had remarked earlier that the distributions $\,\sigma\in C^{-\infty}(\R)$ which vanish to order $k\leq \infty$ at $x=0$ do not constitute a closed subspace, relative to the strong distribution topology. To put a useful topology on this space, one can use the methods of the previous section to translate the problem into a tractable problem in complex analysis. Alternatively, one can use the local description of distributions in definition \ref{inf_order} or \lemref{one_var_def} to define an appropriate topology. We shall pursue the latter strategy, which has the advantage of working just as well in the context of manifolds. Even though we shall state and prove certain results without mentioning the topology explicitly, the use of a topology will be visible in the background.

We begin with a version of definition \ref{inf_order} with dependence on parameters. Again $\,M$ denotes a $C^\infty$ manifold, and $\,S\subset M$ a locally closed submanifold. We consider a family of distributions $\,\sigma_n\in C^{-\infty}(M)$ indexed by $\,n=(n_1,\,n_2,\,\dots,\,n_d)\in \Z^d$, or more generally, indexed by $d$-tuples $\,n$ of integers ranging over some subset of $\,\Z^d$. For $\,n=(n_1,\,n_2,\,\dots,\,n_d)\in \Z^d$ and $\,y=(y_1,\,y_2,\,\dots,\,y_d)\in \R^d$, we let $\,ny$ denote the sum $\,\sum_j n_j y_j$.

\begin{definition}\label{multi_var_def}
The family $\,\sigma_n$ vanishes to order $\,k\geq 0$ along $\,S$, uniformly in $\,n$, if every $\,p\in S$ has a coordinate neighborhood $\,U_p$ in $\,M$ on which
\[
\sigma_n\ \ = \ \ {\sum}_{1\leq j\leq N}\ f_{n,j}\, D_{n,j}\, h_{n,j}\,,
\]
with $\,h_{n,j}\in L^\infty(U_p)$, with $\, f_{n,j}\in C^\infty(U_p)$ vanishing to order $\,k$ along $\,S\cap U_p$, and with differential operators $\,D_{n,j}$ on $\,U_p$ which are tangential to $\,S\cap U_p$, of order $\,r$ -- which may depend on $\,k$ and $\,U_p$, but not on $\,n$ -- such that the $L^\infty$ norms $\,\|h_{n,j}\|_\infty$ are bounded by a polynomial in $\,\|n\|$, and such that the coefficient functions of the $\,D_{n,j}$ as well as the $\,f_{n,j}$ are uniformly bounded, along with all their derivatives up to order $\,k+2r$. When this is the case for every $k\geq 0$, we say that the $\,\sigma_n$ vanish to infinite order along $\,S$, uniformly in the parameter $\,n$.
\end{definition}

The definition involves the choice of a coordinate system for the sole purpose of comparing the sizes of the $\,f_{n,j}$ and $\,D_{n,j}$ for various $\,n$. We shall soon argue that both the $\,f_{n,j}$ and $\,D_{n,j}$ can be made independent of $\,n$, so the particular choice of coordinate functions does not matter at all -- as can also be seen directly, of course.

\begin{lem}\label{multi_var_lem1}
A family $\,\sigma_n\,$, $\,n\in\Z^d$, vanishes to order $\,k\geq 0$ along $\,S$, uniformly in $\,n$, if and only if the series
\[
\tau(p,y) \ \ = \ \ {\sum}_{n\in\Z^d}\ \sigma_n(p)\,e(ny)\qquad \bigl(\,(p,y)\in M\times \R^d/\Z^d\,\bigr)
\]
converges on some open neighborhood of $\,S\times \R^d/\Z^d$ in $\,M\times \R^d/\Z^d$ in the strong distribution topology, to a distribution which vanishes to order $\,k$ along $\,S\times \R^d/\Z^d$.
\end{lem}

\begin{proof}
We shall argue locally, as we may. For the ``only if" direction, we change the coordinates on the coordinate neighborhood $\,U_p$ so that
\begin{equation}
\label{multi_var_1}
S \cap U_p\ \ = \ \ \{\,x=(x_1,\, x_2,\,\dots,\,x_m)\in U_p \mid x_1=\dots =x_c = 0\,\}\,.
\end{equation}
We then write
\begin{equation}
\label{multi_var_2}
f_{n,j}\, D_{n,j}\ \ = \ \ {\sum}_{}\ a_{n,j,I,J}\,\textstyle\frac{\partial^{|I|}\ }{\partial x_I}\,\frac{\partial^{|J|}\ }{\partial x_J}\,,
\end{equation}
with $\frac{\partial^{|I|}\ }{\partial x_I}$ running over all monomials in the $\,\frac{\partial\ }{\partial x_i}$, $\,1\leq i \leq c$, of degree $\,|I|\leq r$, and $\,\frac{\partial^{|J|}\ }{\partial x_J}$ over all monomials in the $\,\frac{\partial\ }{\partial x_j}$, $\,c < j \leq m\,$, also with $|J|\leq r$; moreover, $\,a_{n,j,I,J}=0\,$ unless $\,|I|+|J|\leq r$. Because of the hypotheses, the functions $\,a_{n,j,I,J}\,$ and all their partial derivatives up to order $\,k+2r$ are bounded independently of $\,n$, and each $\,a_{n,j,I,J}$ vanishes along $\,S\cap U_p$ to order $\,|I|+k$. Thus
\begin{equation}
\label{multi_var_3}
a_{n,j,I,J}\ \ = \ \ {\sum}_L\ x^L\,b_{n,j,I,J,L}\,,
\end{equation}
where now $\,x^L$ runs over all monomials in the $\,x_i$, $1\leq i \leq c$, of degree $\,|L| = |I|+ k$. In passing from the $\,a_{n,j,I,J}$ to the $\,b_{n,j,I,J,L}$, the bound on the partial derivatives gets weakened: the latter functions have partial derivatives bounded independently of $\,n$ up to order $\,2r-|I|\geq |I|+|J|$. When we substitute the expressions (\ref{multi_var_3}) into (\ref{multi_var_2}) and commute the $\,b_{n,j,I,J,L}$ across the operators $\,\frac{\partial^{|I|}\ }{\partial x_I}\,\frac{\partial^{|J|}\ }{\partial x_J}$, the result is an expression
\begin{equation}
\label{multi_var_4}
\sigma_n\ \ = \ \ {\sum}_{1\leq j\leq N}\ f_j\, D_j\, h_{n,j}
\end{equation}
as in definition \ref{multi_var_def}, but with functions $\,f_j$ and differential operators $\,D_j$ which no longer depend on $\,n$; the $\,h_{j,n}$ are $L^\infty$ functions whose norm still grows at most polynomially in $\,\|n\|$. Hence, for $\,s\in\N$ sufficiently large, the series
\begin{equation}
\label{multi_var_5}
\widetilde h_j(x,y)\ \ = \ \ {\sum}_{n\in\Z^d}\ (1 + 4\pi^2\|n\|^2)^{-s}\,h_{n,j}(x)\,e(ny)
\end{equation}
converges uniformly, to a bounded measurable function on $\,U_p\times \R^d/\Z^d$. By construction, the series
\begin{equation}
\label{multi_var_6}
\begin{aligned}
{\sum}_{n\in\Z^d}\ \sigma_n(x)\,e(ny)\ \ = \ \ {\sum}_{1\leq j\leq N}\ f_j\,\widetilde D_j\,\widetilde h_j(x,y)\,,
\\
\text{with}\ \ \ \ \widetilde D_j\ = \ D_j\,\bigl(1 - \sum_{i=1}^d\ \textstyle\frac{\partial^2\ }{\partial y_i^2}  \bigr)^s\,,
\end{aligned}
\end{equation}
converges to a distribution $\,\tau\in C^{-\infty}(U_p\times \R^d/\Z^d)\,$ in the strong distribution topology, and $\,\tau$ vanishes along $\,S\times \R^d/\Z^d$ to order $\,k$.

For the argument in the ``if" direction, we fix $\,p\in S$. Since $\,\R^d/\Z^d$ is compact, there exists an open neighborhood $\,U_p$ of $\,p$ in $M$ such that the open set on which the series for $\,\tau$ converges contains $\,U_p \times \R^d/\Z^d$. The local expressions for $\,\tau$ in definition \ref{inf_order} are only required to exist locally, but we can use a partition of unity to glue together such local expressions to get one that is valid on a neighborhood of $\,\{p\} \times \R^d/\Z^d$; equivalently, if we shrink $\,U_p$, we can get this type of expression globally on $\,U_p \times \R^d/\Z^d$. Shrinking $\,U_p$ further, if necessary, we may suppose that there exist coordinate functions $\,x_j$ on $\,U_p$, as in (\ref{multi_var_1}). We now argue as in the first half of the proof to put the local expression for $\,\tau$ into the following form:
\begin{equation}
\label{multi_var_7}
\tau(x,y)\ \ = \ \ {\sum}_{1\leq j\leq N}\ f_j(x)\, D'_j(x)\,D''_j(y)\, h_j(x,y)\,;
\end{equation}
here the $\,f_j\in C^\infty(U_p)$ vanish to order $k$ along $\,S\cap U_p$, the $\,D'_j(x)$ are differential operators on $\,U_p$ with polynomial coefficients, tangential to $\,S\cap U_p$, the $\,D''_j(y)$ are constant coefficient differential operators in the $\,y_i$, and the $\,h_j(x,y)$ are $L^\infty$ functions depending on both sets of variables. The torus $\,\R^d/\Z^d$ acts continuously on $\,C^{-\infty}(U_p\times \R^d/\Z^d)$. Taking Fourier coefficients, we find
\begin{equation}
\label{multi_var_7a}
\sigma_n(x)\ \ = \ \ {\sum}_{1\leq j\leq N}\ f_j(x)\, D'_j(x)\,\int_{\R^d/\Z^d} e(-ny)\,D''_j(y)\, h_j(x,y)\,dy\,.
\end{equation}
Integration by parts exhibits the integral as a bounded, measurable function of the $\,x_j$ whose $L^\infty$ norm is bounded by a polynomial in $\,\|n\|$. This is the kind of local expression for the $\,\sigma_n$ required by definition (\ref{multi_var_def}), with $\,f_{n,j}=f_j$ and $\,D_{n,j}=D'_j\,$ both independent of the parameter $\,n\in\Z^d$.
\end{proof}

For future reference, we record a fact which we just established in the course of the proof of \lemref{multi_var_lem1}:

\begin{cor}\label{multi_var_cor1}
In the setting of definition \ref{multi_var_def}, it is possible to choose the $\,f_{n,j}=f_j$ and $\,D_{n,j}=D_j$ independently of $\,n$. In terms of any local coordinate neighborhood $\,(U_p;x_1,\,x_2,\,\dots,\,x_m)$ such that
\[
S \cap U_p\ \ = \ \ \{\,x=(x_1,\, x_2,\,\dots,\,x_m)\in U_p \mid x_1=\dots =x_c = 0\,\}\,,
\]
the $\,f_j$ can be assumed to be polynomial functions and the $\,D_j$ differential operators with polynomial coefficients.
\end{cor}

The fact that the $\,f_{n,j}=f_j$ and $\,D_{n,j}=D_j$ can be chosen independently of $\,n$ makes it easy to extend our earlier results to families depending on parameters. We begin with a version of \propref{xalpha} for families of distributions; the earlier proof carries over almost word-for-word. As in \propref{xalpha}, we suppose that the submanifold $\,S \subset M$ has a global defining function $\,f\in C^\infty(\R)$ whose differential is non-zero at every point of $\,S$.

\begin{lem}\label{multi_var_lem2}
Let $\,\sigma_n\in C^{-\infty}(M)$, $\,n\in\Z^d$, be a family which vanishes along $\,S$ to order $\,0\leq k\leq\infty$, uniformly in $\,n$. If $\,\alpha, \beta \in \C$ and $\,\ell\geq 0$ satisfy the conditions in part a) of \propref{xalpha}, the family of distributions $\,(\sg f)^\delta |f|^\alpha (\log|f|)^\beta\sigma_n$ vanishes along $\,S$ to order $\,\ell$, also uniformly in $\,n$.
\end{lem}

Periodic distributions without constant term furnish the simplest example of distributions vanishing to infinite order. Such distributions can be represented as $k$-th derivatives of bounded continuous functions, for every $k\gg 0\,$ -- that is the crux of the proof of proposition \ref{oscillation}. In the setting of families, the same reasoning implies:

\begin{lem}\label{multi_var_lem3}
Let $\,\sigma_n\in C^{-\infty}(\R)$, $\,n\in\Z^d$, be a family of distributions which, for every sufficiently large $k\in\N$, can be expressed as $\,\sigma_n(x) = f_n^{(k)}(x)$, with $f_n$ continuous, bounded, and $\,\sup|f_n| = O(\|n\|^m)$ for some $m\in \N$. Then
the $\,\sigma_n$ have canonical extensions across $\,\infty\,$ which vanish there to infinite order, uniformly in $\,n$.
\end{lem}

Definition \ref{multi_var_def} imposes conditions on the $\,\sigma_n$ only near points of $\,S$. The Fourier transform of a tempered distribution is a global operation, so the extension of \thmref{fourier} to the present setting also requires a global hypothesis. We shall say that a family of tempered distributions $\,\sigma_n\in\Sch'(\R)$, $\,n\in\Z^d$, is bounded if there exist positive integers $m$, $k$, $\ell$, such that
\begin{equation}
\label{multi_var_8}
\begin{aligned}
\sigma_n(x)\ = \ \textstyle\frac{d^k\ }{dx^k}\,f_n(x)\,,\ \ \text{with}\ \ f_n\in C(\R)\ \ \text{and}
\\
{\sup}_{x\in\R}\bigl((1+x^2)^{-\ell}|f_n(x)|\bigr)\ = \ O(\|n\|^m)\,.
\end{aligned}
\end{equation}
If the $\,\sigma_n$ constitute a bounded family, then the family of Fourier transforms $\,\widehat \sigma_n$ is bounded, too.

\begin{lem}\label{multi_var_lem4}
If $\,\sigma_n\in\Sch'(\R)$, $\,n\in\Z^d$, is a bounded family which vanishes at $x=0$ to order $k\geq 0$, uniformly in $\,n$, the $\,\widehat\sigma_n$ extend to distributions on $\,\R\cup\{\infty\}$ which vanish at $\infty$ to order $k$, uniformly in $\,n$.
\end{lem}

The proof of \thmref{fourier} can easily be adapted to the present situation. This is also true for the proofs of lemmas \ref{one_var_def} and \ref{compactsupport}, which are used in the proof of \thmref{fourier}.

Other constructions for tempered distributions can also be carried for bounded families. We now state a number of results in this direction; all can be proved by keeping track of bounds in the analogous arguments for single distributions.

\begin{lem}\label{multi_var_lem5}
In the statements that follow $\,\,\sigma_n$, $\,n\in\Z^d$, is a bounded family of tempered distributions.\newline
\noindent {\rm a)\,} The family of distributions $\sigma_n(1/x)$ is bounded, provided the $\sigma_n$ have been extended to distributions on $\,\R\cup\{\infty\}$ which vanish at $\infty$ to order $k$ in the uniform sense, for some $\,k\geq 0$.\newline
\noindent {\rm b)\,} For $\alpha\in\C$ and $\delta\in\Z/2\Z$, the family $\,(\sg x)^\d |x|^\alpha \sigma_n(x)$ is bounded.\newline
\noindent {\rm c)\,} If $c_n\in \R^*$, $\,n\in\Z^d$, is a family of constants such that $|c_n| = O(\|n\|^m)$ and $|c_n|^{-1} = O(\|n\|^m)$ for some $m\in \N$, then $\,\sigma_n(c_n x)$ is a bounded family. If the original family vanishes to order $k\geq 0$ at the origin in the uniform sense, then so does the family $\,\sigma_n(c_n x)$.\newline
\noindent {\rm d)\,} If $f\in C^\infty(\R)$ and all of its derivatives grow at most polynomially as $|x|\to \infty$, the family $f(x)\sigma_n(x)$ is bounded. If the original family vanishes to order $k\geq 0$ at the origin in the uniform sense, then so does the family $\,f(x)\sigma_n(x)$.
\end{lem}

We conclude with a discussion of families of distributions arising from a geometric action. The submanifold $\,S\subset M$ will now assumed to be closed. We suppose that a Lie group $\,H$ acts smoothly on $M$, and that the action preserves $\,S\,$:
\begin{equation}
\label{action_1}
\ell_h\,S\, \subset \, S\ \ \ \text{for every}\ \ h \in H\,;
\end{equation}
here $\,\ell_h : M \to M$ denotes translation by $\,h$. We shall need to consider not only scalar distributions, but also distributions with values in an $\,H$-equivariant vector bundle $\,\mathcal E \to M$, i.e., a vector bundle to which the action of $\,H$ on $M$ lifts. Then $\,H$ acts on
\begin{equation}
\label{action_2}
C^{-\infty}(M,\mathcal E)\ \ =\ \ \text{space of $\,\mathcal E$-valued distributions on $M$}\,.
\end{equation}
Locally the datum of an $\,\mathcal E$-valued distribution amounts to an $\,r$-tuple of scalar distributions, with $\,r\,$=\,rank of $\,\mathcal E$. The notion of vanishing to order $\,k$ along $\,S$ therefore has meaning for $\,\mathcal E$-valued distributions. If $\,dh\,$ is a smooth measure on $\,H$ -- such as left or right Haar measure, for example -- and $\,\sigma$ an $\,\mathcal E$-valued distribution, the family $\,h\mapsto \ell_h \sigma$ can be integrated with respect to $\,dh\,$ over any compact measurable subset $\,\Omega\subset H$:
\begin{equation}
\label{action_3}
\int_\Omega \ell_h\,\sigma\,dh \ \in \ C^{-\infty}(M,\mathcal E)\,.
\end{equation}
Typically this type of integral arises when both $\,\sigma$ and $\,dh\,$ are invariant under a cocompact discrete subgroup $\,\Gamma \subset H$. In that case one may want to integrate $\,\ell_h \sigma$ over $\,\Gamma\backslash H$, or equivalently, over a fundamental domain $\,\Omega$ for the action of $\,\Gamma$ on $\,H$.

\begin{prop}\label{multi_var_prop}
In the situation {\rm (\rangeref{action_1}{action_3})}, if $\,\sigma \in C^{-\infty}(M,\mathcal E)$ vanishes to order $k\geq 0$ along $\,S$, then so does $\,\int_\Omega \ell_h\, \sigma\, dh\,$.
\end{prop}

\begin{proof}
We first give the argument for a scalar valued distribution $\,\sigma$. Partitions of unity for both $M$ and $\,H$ make it possible to reduce the problem to the following situation: the translates $\,\ell_h \sigma$, with $h\in\Omega$, all have compact support in a coordinate neighborhood $\,U$ as in (\ref{multi_var_1}). We choose a reference point $\,h_0\in\Omega$ and use corollary \ref{multi_var_cor1} -- for the ``trivial family" $\,\ell_{h_0}\sigma$ without dependence on a parameter $\,n$ -- to write
\begin{equation}
\label{action_4}
\ell_{h_0}\,\sigma(x)\ \ = \ \ {\sum}_{1\leq j\leq N}\ f_j(x)\,D_j(x)\,h_j(x)\,,
\end{equation}
in terms of polynomial functions $\,f_j$ which vanish on $\,S\cap U$ to order $\,k$ and differential operators $\,D_j$ with polynomial coefficients which are tangential to $\,S\cap U$. We enlarge the collection of $\,f_j$ and $\,D_j$ until we get finite generating sets over the polynomial algebra $\,\C[x]$ for the ideal of polynomials vanishing on $\,S\cap U$ to order $\,k$ and the space of differential operators $\,D_j$ tangential to $\,S\cap U$, of degree up to the maximum necessary in the expression (\ref{action_4}). The action of $\,H$ preserves the order of vanishing of functions along $\,S$ and the notion of tangentiality of a differential operator. It follows that there exist matrices of $C^\infty$ functions $a_{j,\ell}(h,x)$, $b_{j,\ell}(h,x)$, such that
\begin{equation}
\label{action_5}
\bigl(\ell_h f_j\bigr)(x)  = \,  {\sum}_{i}\, a_{j,i}(h,x)\,f_i(x)\,,\ \ \bigl(\ell_h D_j\bigr)(x)  = \,  {\sum}_{i}\, b_{j,\ell}(h,x)\,D_\ell(x)\,.
\end{equation}
Then
\begin{equation}
\label{action_6}
\ell_{h}\,\sigma(x)\ \ = \ \ {\sum}_{i,\,\ell,\,j}\ a_{j,i}(h,x)\,b_{j,\ell}(h,x)\,f_{i}(x)\,D_{\ell}(x)\,\bigl(\ell_{h}\,h_j\bigr)(x)\,,
\end{equation}
for $\,h\in\Omega$ and $\,x\in U$. We now move the $\,a_{j,i_1}$ and $\,b_{j,i_2}$ across the $\,D_\ell$. Arguing as in the proof of \lemref{one_var_def}, but in higher dimension, we can transform (\ref{action_6}) into an expression where the dependence on $\,h$ appears to the right of the differential operators:
\begin{equation}
\label{action_7}
\ell_{h}\,\sigma(x)\ \ = \ \ {\sum}_{i,\,\ell,\,j}\ f_{i}(x)\,D_{\ell}(x)\bigl(c_{i,\ell,j}(h,x)\,\bigl(\ell_{h}\,h_j\bigr)(x)\bigr)\,,
\end{equation}
with coefficient functions $\,c_{i,\ell,j}$ which are products of partial derivatives of the $\,a_{j,i}$ and of the $\,b_{j,\ell}$ and coordinate functions. Thus
\begin{equation}
\label{action_8}
\int_\Omega \ell_h\,\sigma\,dh \ \ = \ \ {\sum}_{i,\,\ell,\,j}\ f_{i}(x)\,D_{\ell}(x)\biggl(\int_\Omega c_{i,\ell,j}(h,x)\,\bigl(\ell_{h}\,h_j\bigr)(x)\,dh\biggr)
\end{equation}
does vanish to order $\,k$ along $\,S\cap U$, as was to be shown.

If $\,\sigma$ takes values in an $\,H$-equivariant vector bundle $\,\mathcal E$, we shrink the coordinate neighborhood $U$ so that $\,\mathcal E$ can be trivialized over $U$. We can then identify the $\mathcal E$-valued distribution $\,\sigma$ with an $\,r$-tuple of scalar distributions $(\sigma_1,\,\sigma_2,\,\dots,\,\sigma_r)$. The action of $\,H$, expressed in terms of the $\,r$-tuple, involves a matrix-valued factor of automorphy:
\begin{equation}
\label{action_9}
\begin{aligned}
&\ell_h(\sigma_1,\,\sigma_2,\,\dots,\,\sigma_r)(x) \ \ =
\\
&\ \ \ \ \ =\ \ \bigl(\,\textstyle{\sum}_{j}\, A_{1,j}(h,x)\,\ell_h\sigma_{j}(x)\,,\,\dots,\,{\sum}_{j}\, A_{r,j}(h,x)\,\ell_h\sigma_{j}(x)\,\bigr)\,,
\end{aligned}
\end{equation}
with $C^\infty$ coefficients $\,A_{i,j}(h,x)$. These must be moved across the $\,D_\ell$ along with the $\,a_{j,i}$ and $\,b_{j,\ell}$. Otherwise the argument remains the same.
\end{proof}
\bigskip

\begin{tabular}{lcl}
Stephen D. Miller                    & & Wilfried Schmid \\
Department of Mathematics            & & Department of Mathematics \\
Hill Center-Busch Campus             & & Harvard University \\
Rutgers University                        & & Cambridge, MA 02138 \\
110 Frelinghuysen Rd                 & & {\tt schmid@math.harvard.edu}\\
Piscataway, NJ 08854-8019            & & \\
{\tt miller@math.rutgers.edu}
\end{tabular}

\bibliographystyle{amsxport}

\end{document}